\documentclass[preprint,a4paper,10pt,oneside,onecolumn]{elsarticle}
\usepackage{graphicx}
\usepackage{latexsym,amsmath,amsfonts,amscd,amsthm,amssymb}
\usepackage{changebar}
\usepackage{color}
\usepackage{bm}
\usepackage{tikz}
\usepackage{multirow,multicol}
\usetikzlibrary{arrows,backgrounds,snakes}
\usepackage{epsfig,subfigure}
\usepackage{ulem}
\usepackage[english]{babel}

\newtheoremstyle{plainNoItalics}{}{}{\normalfont}{}{\bfseries}{.}{ }{}
\theoremstyle{plain}

\usepackage[linesnumbered,ruled,lined]{algorithm2e}
\usepackage{graphicx}
\usepackage[colorlinks=true, allcolors=blue]{hyperref}
\usepackage{indentfirst}
\usepackage{amssymb,amsmath,amsthm}
\usepackage{epsfig}
\usepackage{appendix}
\usepackage{enumerate}
\usepackage{comment}
\usepackage{color}
\usepackage{amsfonts}
\usepackage{mathrsfs}
\usepackage{tikz}
\usepackage{booktabs}
\usepackage{epsfig,graphics, color,bm}

\newtheorem{example}{\indent Example}[section]
\newtheorem{tableau}{\indent Tableau}[section]

\newtheorem{thm}{Theorem}

\newdefinition{rmk}{Remark}
\newproof{pf}{Proof}
\newproof{pot}{Proof of Theorem \ref{thm2}}

\newcommand{\bx}{{\bf x}}

\begin{document}
	
	\begin{frontmatter}
		\title{High order conservative LDG-IMEX methods for the degenerate nonlinear non-equilibrium radiation diffusion problems\tnoteref{t1}}
		\tnotetext[t1]{The first and the fourth authors are supported 
			by National Key R\&D Program of China No. 2022YFA1004500, NSFC No. 11971025 and No. 92270112, NSF of Fujian Province No. 2023J02003. The second author is partially
			supported by NSFC No. 12031013, Shanghai pilot innovation project No. 21JC1403500 and the Strategic Priority Research Program of Chinese Academy of Sciences Grant No. XDA25010401. The third author is supported by NSFC No. 12071214.}
		
		\author[1]{Shaoqin Zheng}
		\ead{sqzheng@stu.xmu.edu.cn}
		
		\author[2]{Min Tang}
		\ead{tangmin@sjtu.edu.cn}
		
		\author[3]{Qiang Zhang}
		\ead{qzh@nju.edu.cn}
		
		\author[1]{Tao Xiong\corref{cor1}}
		\ead{txiong@xmu.edu.cn}
		
		\cortext[cor1]{Corresponding author}
		
		\affiliation[1]{organization={School of Mathematical Sciences},
			addressline={Xiamen University},
			city={Xiamen},
			postcode={Fujian 361005},
			country={PR China}}
		\affiliation[2]{organization={Institute of Natural Sciences},
			addressline={Shanghai Jiao Tong University},
			city={Shanghai},
			postcode={200240},
			country={PR China}}
		\affiliation[3]{organization={Department of Mathematics},
			addressline={Nanjing University},
			city={Nanjing},
			postcode={210093},
			country={PR China}}
		\affiliation[4]{organization={School of Mathematical Sciences, Fujian Provincial Key Laboratory of Mathematical Modeling and High-Performance Scientific Computing},
			addressline={Xiamen University},
			city={Xiamen},
			postcode={Fujian 361005},
			country={PR China}}

	\begin{abstract}
        In this paper, we develop a class of high-order conservative methods for simulating non-equilibrium radiation diffusion problems. Numerically, this system poses significant challenges due to strong nonlinearity within the stiff source terms and the degeneracy of nonlinear diffusion terms. Explicit methods require impractically small time steps, while implicit methods, which offer stability, come with the challenge to guarantee the convergence of nonlinear iterative solvers. To overcome these challenges, we propose a predictor-corrector approach and design proper implicit-explicit time discretizations. In the predictor step, the system is reformulated into a nonconservative form and linear diffusion terms are introduced as a penalization to mitigate strong nonlinearities. We then employ a Picard iteration to secure convergence in handling the nonlinear aspects. The corrector step guarantees the conservation of total energy, which is vital for accurately simulating the speeds of propagating sharp fronts in this system. 

        For spatial approximations, we utilize local discontinuous Galerkin finite element methods, coupled with positive-preserving and TVB limiters. We validate the orders of accuracy, conservation properties, and suitability of using large time steps for our proposed methods, through numerical experiments conducted on one- and two-dimensional spatial problems. In both homogeneous and heterogeneous non-equilibrium radiation diffusion problems, we attain a time stability condition comparable to that of a fully implicit time discretization. Such an approach is also applicable to many other reaction-diffusion systems.
	\end{abstract}
	\begin{keyword}
		non-equilibrium radiation diffusion \sep predictor-corrector procedure \sep conservative \sep high order \sep local discontinuous Galerkin method \sep IMEX
		
	\end{keyword}
	
\end{frontmatter}
	
	\section{Introduction}
	\label{sec_intro}
	\setcounter{equation}{0}
	\setcounter{figure}{0}
	\setcounter{table}{0}
	
        In scenarios where radiation interacts dynamically with a material, but does not fully reach thermodynamic equilibrium, a commonly employed description involves a system of coupled time-dependent nonlinear diffusion equations. This system is known for its strong nonlinearity and close coupling, and it finds extensive applications across various fields, including inertial confinement fusion \cite{turner2001module}, astrophysics \cite{bowers1991numerical}, and Z-pinch experiments \cite{robinson2004three}. In this work, we specifically explore a two-temperature (2T) model. This model comprises a radiation diffusion equation within the framework of a gray approximation and a material energy balance equation \cite{bowers1991numerical,peterson1996two,winslow1995multifrequency}. The 2T model writes:
	\begin{equation}\label{RDE}
		\begin{cases}
			\dfrac{\partial E}{\partial t}-\nabla \cdot\left(D_r \nabla E\right)=\sigma\left(T^{4}-E\right),
			\\[5pt]
			\dfrac{\partial T}{\partial t}-\nabla \cdot\left(D_t \nabla T\right)=\sigma\left(E-T^{4}\right),
		\end{cases}
	\end{equation}
	where $E(\bx, t)$ is the radiation energy and $T(\bx, t)$ is the material temperature. 
        The energy exchange between materials and photons is controlled by the photon absorption cross-section $\sigma$, which is modeled by
	\begin{equation}\label{sigma}
		\sigma(T) = \frac{z(\bx)^3}{T^3},
	\end{equation}
	where $z(\bx)$ is a spatial dependent material coefficient that represents an atomic mass number. The following flux-limited energy radiation diffusion coefficient $D_r$ is widely used \cite{bowers1991numerical,knoll1999efficient}:
	
	\begin{equation}\label{Dr}
		D_r = \frac{1}{3\sigma+\frac{|\nabla E|}{E}}.
	\end{equation}
	Moreover, the material conduction coefficient has the following form \cite{spitzer1953transport}
	\begin{equation}\label{Dt}
		D_t = \kappa T^{\frac{5}{2}},
	\end{equation}
	where $\kappa$ is a constant. 
		
	Adding the two equations of \eqref{RDE} together, one gets
	\begin{equation}\label{conservation_quantity}
		\dfrac{\partial( E + T)}{\partial t}=\nabla \cdot\left(D_r \nabla E\right)+\nabla \cdot\left(D_t \nabla T\right),
	\end{equation}
	which indicates that the system \eqref{RDE} conserves the energy $E+T$ over the space. Energy conservation is not only physically important, but also crucial to capture a useful simplification of the 2T model. As $\sigma$ approaches $\infty$, the radiation energy tends to approximate the thermal equilibrium, namely $E\approx T^4$, \eqref{conservation_quantity} would lead to the gray radiative diffusion equation \cite{Larsen1983}
	\begin{equation}\label{diff_limit}
	\dfrac{\partial( T^4 + T)}{\partial t}=\nabla \cdot\left(D_r \nabla T^4\right)+\nabla \cdot\left(D_t \nabla T\right).
    \end{equation}	
    The gray radiation equation is an important approximation widely employed for studying diverse radiative heat transfer phenomena, including radiative transfer in stellar atmospheres.
    	
	Non-equilibrium radiation diffusion equations have attracted extensive research efforts. Marshak, for instance, developed a time-dependent radiative transfer model to investigate the impact of radiation on shock wave behavior \cite{marshak1958effect}. Analytical solutions have also been provided for specific Marshak wave problems in previous studies \cite{pomraning1979non, bingjing1996benchmark}. However, solving the system \eqref{RDE} numerically presents substantial challenges, mainly due to the following several reasons:
	\begin{itemize}
		\item The source terms might be very stiff. Specifically, the coefficient $\sigma$, as defined in \eqref{sigma}, tends to be very large for large values of $z$ or small values of $T$.
		This leads to pronounced energy exchanges and a close coupling between $E$ and $T^4$. Dealing with such stiff source terms when $\sigma$ is large requires implicit time discretizations, resulting in a highly nonlinear system with multiscale coefficients. 
		
		\item The diffusion terms $\nabla \cdot\left(D_r \nabla E\right)$ and $\nabla \cdot\left(D_t \nabla T\right)$ are nonlinear and degenerate. The coefficients $D_r$ and $D_t$ defined in \eqref{Dr}-\eqref{Dt} nonlinearly depend on $T$, and they may vary by several orders of magnitude when the temperature $T$ varies over a wide range. Moreover, when $T$ is small, both $D_r$ and $D_t $ are close to $0$, which result in degenerate diffusion. Consequently, solutions exhibit sharp fronts that propagate at finite speeds, akin to shock wave solutions in hyperbolic systems \cite{tang2021semi}. 
		
		\item A challenging condition number for a nonlinear iteration. When the temperature $T$ is low, energy $E$ and temperature $T$ operate at significantly different scales, with $E$ roughly proportionate to $T^4$. Consequently, using a fully implicit time discretization for both $E$ and $T$, combined with Newton iteration, can lead to a coefficient matrix with a challenging condition number. It is essential to devise suitable iterative solvers and efficient preconditioners to tackle this challenge effectively.
	\end{itemize}
	
	When solving \eqref{RDE} with explicit schemes, it is relatively straightforward to maintain the conservation of total energy $E+T$, however a significant drawback is the requirement of very small time steps. A parabolic time step $\Delta t=\mathcal{O}(h^2)$ with $\Delta t$ being the time step and $h$ being the mesh size is needed due to the nonlinear diffusion terms, and the stiff source terms lead to a severe time step constraint as well   \cite{yang2019moving}. Fully implicit schemes offer the advantage of using large time steps. However, the presence of strong nonlinearity and multiscale coefficients, when employing a large time step, the solutions from the previous time step may not serve as a suitable initial guess. As a result, iterative solvers used in fully implicit methods may struggle to converge effectively \cite{arbogast2019finite,ketcheson2009optimal}. 
    In past years, significant research efforts have been dedicated to developing efficient iterative solvers with effective preconditioners. One notable approach is the Jacobian-free Newton-Krylov method \cite{knoll1999efficient,knoll2001nonlinear}, which combines Newtonian external iterations and conjugate gradient-like (Krylov) internal iterations, resulting in superlinear convergence without the need to form Jacobian matrices. Additionally, a physical-based preconditioning Newton-Krylov method was explored in \cite{mousseau2000physics}, and an operator-split preconditioner was investigated in \cite{mousseau2003new}. Various other efficient time discretization methods have also been examined in \cite{knoll2003balanced,lowrie2004comparison,brown2005fully,mousseau2006temporal,ober2004studies,knoll2007numerical} and references therein.  
	 Two semi-implicit schemes allowing for large time steps for the gray radiation diffusion equation \eqref{diff_limit} have been developed in \cite{tang2021semi}.  
     Furthermore, radiation diffusion models are frequently coupled with hydrodynamic equations, and large deformations of complex fluid flows lead to mesh distortions. Consequently, for spatial discretization, there have been some works built upon distorted meshes \cite{kang2003p,sheng2009monotone,yu2019finite,zhao2013finite}, as well as moving meshes \cite{yang2015moving,yang2019moving}.  
	 
     Schemes that efficiently handle non-equilibrium radiation diffusion equations, while being robust and capable of accommodating large time steps for accurate sharp front capturing, remain relatively rare in the existing literature. To design such a scheme, several requirements have to be met:
        \begin{itemize}
        \item[I]  Preservation of the energy equation \eqref{conservation_quantity} at the discrete level. This indicates that the source terms in the two equations in \eqref{RDE} are in balance with each other. As far as \eqref{conservation_quantity} is preserved at the discrete level, when $\sigma $ is large and the system reaches the thermodynamic equilibrium $E=T^4$, the scheme's accuracy can be controlled by the gray radiation equation \eqref{diff_limit}; 
		\item[II] Conservation of energy. Because of the presence of the nonlinear degenerate diffusion terms, the solution exhibits hyperbolic properties at the front. As is well-established in the study of hyperbolic simulations, preserving the conservation properties of the scheme's discretization is of paramount importance to accurately capture finite-speed front propagations.
		\item[III]  Positivity Preserving. The system is only physically meaningful for positive radiation energy $E$ and positive material temperature $T$. The radiation temperature $T_r=E^{1/4}$ and all coefficients in \eqref{RDE} depend on $T$. Hence, it would be crucial to maintain positivity of $E$ and $T$ in order to ensure the robustness of the scheme and obtain physically meaningful solutions. 
		\item[IV] Efficient iterative solvers which can guarantee convergence and the providing of a good initial guess. Fully implicit time discretizations coupled with proper spatial discretizaitons may possibly meet the above mentioned requirements, but due to strong nonlinearities and multiscale variations of coefficients, one has to provide a good initial guess and a suitable preconditioner in order to make an iterative solver converging properly.   
	\end{itemize}
	
    We note that radiation diffusion equations belongs to reaction-diffusion equations which appear in many other fields, such as in material sciences \cite{Murray2003}, chemical reactions \cite{Epstein1998}, ecology and biology systems \cite{Jokisaari2017}, and many numerical approaches have been developed for reaction-diffusion equations, e.g. \cite{Zhu2009applications,zhang2014direct,liu2024primal} and reference therein. However, for the radiation diffusion equations \eqref{RDE}, the main challenges are the degeneracy in the nonlinear diffusion operators and stiffness in the source terms. Both are caused by the low material temperature $T$, making degeneracy, stiffness and nonlinearity strongly coupled, while standard reaction diffusion equations usually involve only linear diffusion but various nonlinear reactions. The design of an efficient iterative solver for radiation diffusion equations is very challenging. 
	In this paper, to address such difficulties, we present a novel and efficient method for \eqref{RDE} utilizing implicit-explicit (IMEX) time discretizations, designed to satisfy all four of the aforementioned criteria. Some techniques we develop can also be applied to other reaction diffusion systems. We introduce a predictor-corrector procedure as our proposed approach. In the predictor step, we multiply both sides of the second equation in \eqref{RDE} with $4T^3$, introduce $B=T^4$, and obtain a simpler system for $B$ and $E$ \cite{huang2016monotone}. This new system does not meet requirements I, II, and III. However, we can design an efficient iterative solver capable of ensuring convergence even for large time steps. The predictor step can provide a good initial guess. Then, in the corrector step we solve the energy conservative equation \eqref{conservation_quantity}, and substitute $E+T$ into the first equation of \eqref{RDE} to further update $E$. Requirements I and II are satisfied in the corrector step and the bad condition number for solving $E$ and $T$ together is avoided. Finally, requirement III can be attained through the utilization of suitable spatial discretizations.  

	For requirement IV, several strategies are employed to alleviate the challenges associated with solving a highly nonlinear system. First of all, implicit treatments of nonlinear diffusion terms are avoided by using the idea in \cite{xiong2022high,wang2020local}  from adding the same linear diffusion terms on both sides of those equations, which however are discretized explicitly and implicitly respectively. For both the predictor and corrector steps, only a mildly nonlinear system needs to be solved, and a simple Picard iteration is adopted. 
	Secondly, a local discontinuous Galerkin (LDG) finite element method is employed for spatial discretiztions, due to a shock wave-like behavior in the solutions of radiation \cite{marshak1958effect}. 
	The LDG method was introduced by Cockburn and Shu in \cite{cockburn1998local,cockburn2001runge} for convection-diffusion problems, which is well-suited for h-p adaptivity and is also very good at shock or sharp gradient capturing. For a first order in space LDG method, the iterative matrix for $E$ and $B$ can be shown to be an {\bf M}-matrix \cite{huang2016monotone}. This can guanrantee the convergence of the Picard iteration in the predictor step. Finally, for second and third orders of discretizations, appropriate spatial limiters are applied to control numerical oscillations for sharp gradient solutions and preserve the positivity of radiation energy $E$ and material temperature $T$. After spatial discretizations, the nonlinearity iteration can be written in the form of a linear system, with nonlinearities mainly appear in the diagonal part of the resulting system, a fast convergence can be obtained.

	The rest of this paper is organized as follows. In section \ref{sec_numer}, an efficient IMEX time discretization is introduced via a predictor-corrector procedure, followed by an LDG space discretization in section \ref{LDG}. In section \ref{sec_examp}, one and two spatial dimensional numerical examples are performed to verify the high order accuracy, conservation and large time step conditions, and good performances for capturing sharp fronts of the radiation energy in both homogeneous and heterogeneous media. Finally, a conclusion is drawn in section \ref{Conclusion}.

	\section{Time discretization}
	\label{sec_numer}
	\setcounter{equation}{0}
	\setcounter{figure}{0}
	\setcounter{table}{0}
	
	In this section, we present an efficient IMEX time discretization for \eqref{RDE}, employing a predictor-corrector procedure. The predictor step, though non-conservative, ensures fast and robust nonlinear convergence. In contrast, the corrector step maintains conservation, which is crucial for accurately propagating sharp fronts. This novel approach significantly mitigates nonlinearity as compared to a fully implicit time discretization, while maintaining a generous time step stability condition and ensuring fast and robust convergence. We begin with a first-order IMEX scheme, which serves as a general framework that can be readily extended to higher orders by integrating a multistage Runge-Kutta (RK) or multistep IMEX time-marching approach. Here, we adopt an IMEX RK scheme. The IMEX time discretization can be combined with any front-capturing spatial discretization. In this section, we keep space continuous, with a detailed description of the spatial discretization provided in the next section.
	
	\subsection{First order IMEX scheme}
	Let the discrete time steps be $t_n \, (n=0,1, \ldots)$ and the time step size be $\Delta t=t_{n+1}-t_n$. Instead of solving \eqref{RDE} directly, we propose a predictor-corrector procedure:
	\begin{itemize}
		\item {\bf The predictor step: } \\
		To mitigate the stiffness of the source terms caused by $\sigma=z^3/T^3$, as in \cite{huang2016monotone}, we multiply both sides of the second equation in \eqref{RDE} by $4T^3$. \eqref{RDE} can be rewritten as: 
		\begin{equation}
			\begin{cases}
				\label{RDE_2D}
				\dfrac{\partial E}{\partial t}-\nabla \cdot\left(D_r \nabla E\right)=\sigma \left(B-E\right),
				\\[10pt]
				\dfrac{\partial B}{\partial t}-4\kappa\left[\nabla \cdot\left(T^{\frac{11}{2}} \nabla T\right)-3T^{\frac{9}{2}}|\nabla T|^2\right]=4z^3\left(E-B\right),
								\\[10pt]
								T=B^{1/4}.
			\end{cases}
		\end{equation}
		As discussed in the introduction, $T^4$ in the source term of \eqref{RDE} needs to be treated implicitly to ensure thermal dynamic equilibrium ($E=T^4$) in the stiff regime when $T$ is small. To achieve this without resorting to nonlinear iterations for a given $\sigma$, we introduce a new variable, $B=T^4$, and update a system for both $E$ and $B$ in this predictor step. Utilizing the equation for $B$ not only avoids the stiffness associated with $\sigma$ in the second equation but also maintains balance between $E$ and $B$, ensuring they remain at the same scale when $T$ is small, thus preventing unfavorable condition numbers. Furthermore, in this equation, we continue to use the material temperature $T$ rather than $B$ for the diffusion terms, thereby preventing negative powers of $B$ which can lead to stiffness when $B=T^4$ is small. Once we have determined $B$, we can calculate $T$ using the relationship $B=T^4$. It is noteworthy that throughout this process and in the following, we consistently use $E$ and $T$ as our input and output variables.
		
		To update \eqref{RDE_2D} for $E$ and $B$, the diffusion terms remain nonlinear. An implicit treatment would lead to a complex nonlinear system. Following the approach presented in \cite{xiong2022high,wang2020local}, we introduce two linear diffusion terms $\alpha_0\Delta E$ and $\beta_0\Delta B$ with constant coefficients $\alpha_0$, $\beta_0$ on both sides of \eqref{RDE_2D} for $E$ and $B$, respectively. Starting from $E^n$ and $T^n$, a first-order IMEX scheme for updating predicted solutions $\tilde{E}^{n+1}$ and $\tilde{T}^{n+1}$ is defined as follows:
		\begin{equation}
			\begin{cases}
			\label{RDE_2Dplus}
			\frac{\tilde{E}^{n+1}-E^{n}}{\Delta t}-\nabla \cdot\left(D^n_r \nabla E^n\right)+\alpha_0\,\Delta E^n =\alpha_0\,\Delta \tilde{E}^{n+1}+\tilde\sigma^{n+1} \left(\tilde{B}^{n+1}-\tilde{E}^{n+1}\right),
				\\[10pt]
				\frac{\tilde{B}^{n+1}-B^{n}}{\Delta t}-\mathcal{H}^n+\beta_0\,\Delta B^n=\beta_0\,\Delta \tilde{B}^{n+1}+4z^3\left(\tilde{E}^{n+1}-\tilde{B}^{n+1}\right),
			\\[10pt]
			\tilde{T}^{n+1}=(\tilde{B}^{n+1})^{1/4},
			\end{cases}
		\end{equation}
	    where 
	    \begin{equation}
	    	\label{nonc_H}
	        \mathcal{H} = 4\kappa\left[\nabla \cdot\left(T^{\frac{11}{2}} \nabla T\right)-3T^{\frac{9}{2}}|\nabla T|^2\right].
	    \end{equation}
		Here, the superscript $n$ or $n+1$ denotes the corresponding values at time step $t_n$ or $t_{n+1}$ respectively, similarly in the following. As we can see, in each equation, two added identical terms are discretized differently, one explicitly and the other implicitly. This approach allows us to achieve time stability close to that of an implicit scheme while only needing to solve linearly implicit diffusion operators \cite{wang2020local, xiong2022high}. From the above, we observe that by introducing $B$ and adding the linear diffusion terms, when we solve \eqref{RDE_2Dplus} with an IMEX method, the only nonlinearity arises from the source term $\sigma(E-B)$ due to the coefficient $\sigma$ appearing in the first equation. If $\sigma$ is constant or a predefined function, \eqref{RDE_2Dplus} becomes a simple linear system. Additionally, if $\kappa=0$ for $\mathcal{H}$ in \eqref{nonc_H}, \eqref{RDE_2Dplus} is in a conservative form, and such a 2T model can be efficiently solved.
		
		\item {\bf The corrector step: } \\
		However, if $\kappa$ is not zero, the term $\mathcal{H}$ in \eqref{nonc_H} within the equation of $B$ is in a non-conservative form. This non-conservative form can result in incorrect sharp front propagation, as will be demonstrated in our numerical examples. To address this issue, a corrector step is required. To ensure the conservation of total energy, as shown in \eqref{conservation_quantity}, we begin by replacing the second equation in the original equation \eqref{RDE} with \eqref{conservation_quantity}, resulting in the following:
		\begin{equation}\label{RDE_2Dc}
			\begin{cases}
				\dfrac{\partial E}{\partial t}-\nabla \cdot\left(D_r \nabla E\right)=\sigma\left(T^{4}-E\right),
				\\[10pt]
				\dfrac{\partial Q}{\partial t}=\nabla \cdot\left(D_r \nabla E\right)+\nabla \cdot\left(D_t \nabla T\right),
				\\[10pt]
				T=Q-E.
			\end{cases}
		\end{equation}
		The equation for $Q$ is in a conservative form, allowing for easy conservation of total energy across space. Following \eqref{RDE_2Dplus}, we also introduce two linear diffusion terms, $\alpha_0\Delta E$ and $\gamma_0\Delta Q$, on both sides of \eqref{RDE_2Dc}. Using a first-order IMEX time discretization, we obtain:
		\begin{equation}\label{RDE_2Dplusc}
			\begin{cases}
				\dfrac{E^{n+1}-E^n}{\Delta t}-\nabla \cdot\left(D^n_r \nabla E^n\right)+\alpha_0\,\Delta E^n =\alpha_0\,\Delta E^{n+1} +\tilde{\sigma}^{n+1} \left((T^{n+1})^{4}-E^{n+1}\right),
				\\[10pt]
				\dfrac{Q^{n+1}-Q^n}{\Delta t}-\nabla \cdot\left(D^n_r \nabla E^n\right)-\nabla \cdot\left(D^n_t \nabla T^n\right)+\gamma_0\,\Delta Q^n=\gamma_0\,\Delta Q^{n+1},
				\\[10pt]
				T^{n+1}=Q^{n+1}-E^{n+1}.
			\end{cases}
		\end{equation}
		Solving for $E$ and $Q$ ensures good convergence and the conservation of $Q$.
        As we can see, in \eqref{RDE_2Dplusc}, adding these linear diffusion terms is crucial to simplify the nonlinear system. With this approach, we can first solve a linear system for $Q$ from the second equation. Then, we can use $T=Q-E$ to replace the $T^4$ term in the first equation, allowing us to solve this mildly nonlinear equation to further update $E$. The nonlinearity is simply $(Q-E)^4$ for $E$ and only appears in the diagonal part of the mass matrix after spatial discretization.
		
	\end{itemize}

    \begin{rmk}
    	\label{rem_1}
    	In the corrector step, we utilize the values obtained from the predictor step to determine the coefficients $\sigma$ in the source term and to provide an initial guess for the corrector step. One could contemplate simplifying the scheme by eliminating the predictor step and directly employing a Picard iteration in the corrector step to handle these coefficients. However, we would mention that $\sigma=z^3/T^3$ and $T^4=(Q-E)^4$ appears as a nonlinear term. In case of $T$ being small, with a stiff coefficient, the nonlinear iteration for solving $E$ in \eqref{RDE_2Dplusc} may not converge well, for example, the stand Marshak wave problem in Example \ref{Marshak_homo_problem}. Instead, in the predictor step, except $\sigma$, others linearly depend on $E$ and $B$. Even with a stiff $\sigma$, it appears in the diagonal part of the mass matrix, so that a fast and robust convergence can be obtained. Starting from a good initial guess provided by the predictor step, the corrector step can also converge well.
    \end{rmk}

	\subsection{High order IMEX scheme}
	\label{IMEX time discretization}
	
    The first-order IMEX scheme with a predictor-corrector procedure has offered a highly efficient and versatile framework for solving the 2T model \eqref{RDE}. To attain high-order accuracy in time while preserving these desirable properties, we employ a globally stiffly accurate IMEX RK time discretization \cite{boscarino2013implicit}. However, a multistep IMEX method can also be utilized \cite{Ascher1995}.
    
    We consider a system of additive ordinary differential equations:
    \begin{equation}
	\label{ODE_IMEX}
	\frac{\mathrm{d} \boldsymbol{y}}{\mathrm{d} t}=L(t, \boldsymbol{y})+N(t, \boldsymbol{y}), \quad \boldsymbol{y}\left(t_{0}\right)=\boldsymbol{y}_{0},
    \end{equation}
    where $\boldsymbol{y}= (y_{1}, y_2, \ldots, y_{d})^T$. $L(t, \boldsymbol{y})$ and $N(t, \boldsymbol{y})$ are linear and nonlinear operators, respectively. $L(t, \boldsymbol{y})$ will be discretized implicitly, while $N(t, \boldsymbol{y})$ will be discretized explicitly. An $s$-stage IMEX RK time discretization can be represented by a double Butcher tableau 
    \begin{equation*}\label{DBT}
	\begin{array}{c|c}
		\hat{c} & \hat{A}\\
		\hline
		\vspace{-0.25cm}
		\\
		& \hat{b}^T 
	\end{array}, \ \ \  \qquad
	\begin{array}{c|c}
		{c} & {A}\\
		\hline
		\vspace{-0.25cm}
		\\
		& {b^T} 
	\end{array}.
    \end{equation*}
    Here $A=\left(a_{i j}\right)$ and $\hat{A}=\left(\hat{a}_{i j}\right) \in \mathbf{R}^{s \times s}$. $\hat{A}$ is a strictly lower triangular matrix for explicit parts. For the implicit part, $A$ can be taken as a lower triangular matrix with a nonzero diagonal to get an efficient implementation, which is usually referred to as a diagonally implicit RK (DIRK) scheme. The vectors are $b^{T}=(b_{1},b_{2},\ldots, b_{s})$, $\hat{b}^{T}=(\hat{b}_{1},\hat{b}_{2}, \ldots, \hat{b}_{s})$, $c^{T}=(c_1, c_{2}, \ldots, c_{s})$, and $\hat{c}^{T}=(\hat{c}_1, \hat{c}_{2}, \ldots, \hat{c}_{s})$, where $c_{i}=\sum\limits_{j=1}^{i} a_{i j}$ and $\hat{c}_{i}=\sum\limits_{j=1}^{i-1} \hat{a}_{i j}$. Denoting $t_{n}^{(j)}=$ $t_{n}+c_{j} \Delta t$, $\hat{t}_{n}^{(j)}=$ $t_{n}+\hat{c}_{j} \Delta t$, the solution of \eqref{ODE_IMEX} can be updated from time level $t^{n}$ to $t^{n+1}$ in the following way:
    \begin{equation}
	\begin{cases}\label{IMEX_RK}
		\boldsymbol{Y}^{(1)}=\boldsymbol{y}_{n},
		\\[5pt]
		\boldsymbol{Y}^{(i)}=\boldsymbol{y}_{n}+\Delta t \sum\limits_{j=1}^{i-1} \hat{a}_{i j} N\left(\hat{t}_{n}^{(j)}, \boldsymbol{Y}^{(j)}\right)+\Delta t \sum\limits_{j=1}^{i} a_{i j} L\left(t_{n}^{(j)}, \boldsymbol{Y}^{(j)}\right), \quad 2 \leq i \leq s,
		\\[5pt]
		\boldsymbol{y}_{n+1}=\boldsymbol{y}_{n}+\Delta t \sum\limits_{i=1}^{s} \hat{b}_{i} N\left(\hat{t}_{n}^{(i)}, \boldsymbol{Y}^{(i)}\right)+\Delta t \sum\limits_{i=1}^{s} b_{i} L\left(t_{n}^{(i)}, \boldsymbol{Y}^{(i)}\right).
	\end{cases}
    \end{equation}
    The IMEX RK scheme is called to be globally stiffly accurate, if the coefficients satisfy \cite{boscarino2013implicit} are required to satisfy:
    $$\hat{c}_{s}=c_{s}=1,\text{ and } a_{s j}=b_j, \, \hat{a}_{s j}= \hat{b}_{j},\quad j=1,2,\ldots,s.$$
    With such an IMEX scheme, the final updating of $\boldsymbol{y}_{n+1}$ in \eqref{IMEX_RK} coincides with the last stage of updating $\boldsymbol{Y}^{(s)}$, so that we can take $\boldsymbol{y}_{n+1}=\boldsymbol{Y}^{(s)}$ and avoid the last cumulative step. In \ref{IMEX Butcher tableau}, Butcher tableaux from first order to third order, which are adopted in this work, are provided.

    If we choose $L$ and $N$ in \eqref{ODE_IMEX} based on the first order scheme \eqref{RDE_2Dplus} and \eqref{RDE_2Dplusc}, with \eqref{IMEX_RK}, the updating of the solutions at $t^{n+1}$ from $t^n$ can be presented as follows:
	\begin{subequations}
		\label{semi_PDE_final_cor}
		\begin{align}
			\frac{E^{n+1}-E^{n}}{\Delta t}&=\sum_{i=1}^{s-1} \hat{b}_{i}\left[\nabla\cdot(D^{(i)}_r\nabla E^{(i)})-\alpha_0\Delta E^{(i)}\right]+\sum_{i=1}^s b_{i}\left[\alpha_{0} \Delta E^{(i)} +{\sigma}^{(i)}(B^{(i)}-E^{(i)})\right],
			\\
			\frac{Q^{n+1}-Q^{n}}{\Delta t}&=\sum_{i=1}^{s-1}\hat{b}_{i}\left[\nabla \cdot\left(D^{(i)}_r \nabla E^{(i)}\right)+\nabla \cdot\left(D^{(i)}_t \nabla T^{(i)}\right)-\gamma_{0}\Delta Q^{(i)}\right] + \gamma_{0} \sum_{i=1}^s b_{i} \Delta Q^{(i)},
			\\
			T^{n+1}&=Q^{n+1}-E^{n+1}.
		\end{align}
	\end{subequations}
    Here similarly the superscript $(i)$ denotes variables at the time stage $t_n^{(i)}$ or $\hat{t}_n^{(i)}$, and the intermediate stage values for $2 \leq i \leq s$ are obtained from:
    \begin{itemize}
	\item {\bf the predictor step: } \\
	\begin{subequations}
		\label{semi_PDE_immediate}
		\begin{align}
			\frac{\tilde{E}^{(i)}-E^{n}}{\Delta t}&=\text{RHS}_E^{(i)}+a_{i i}\left[\alpha_{0} \Delta \tilde{E}^{(i)}+\tilde{\sigma}^{(i)}(\tilde{B}^{(i)}-\tilde{E}^{(i)})\right],\label{semi_PDE_immediate_1}
			\\
			\frac{\tilde{B}^{(i)}-B^{n}}{\Delta t}&=\text{RHS}_B^{(i)}+a_{i i}\left[\beta_{0} \Delta \tilde{B}^{(i)}+z^3(\tilde{E}^{(i)}-\tilde{B}^{(i)})\right],
			\label{semi_PDE_immediate_2}
			\\
			\tilde{T}^{(i)} &= (\tilde{B}^{(i)})^{1/4},
		\end{align}
    \end{subequations}
	\item {\bf the corrector step: } \\
	\begin{subequations}
		\label{semi_PDE_immediate_cor}
		\begin{align}
			\frac{E^{(i)}-E^{n}}{\Delta t}&=\text{RHS}_E^{(i)}+a_{i i}\left[\alpha_{0} \Delta E^{(i)}+\tilde{\sigma}^{(i)}(B^{(i)}-E^{(i)})\right],
			\\
			\frac{Q^{(i)}-Q^{n}}{\Delta t}&=\text{RHS}_Q^{(i)}+a_{i i}\gamma_{0} \Delta Q^{(i)},
			\\
			T^{(i)} & = Q^{(i)}-E^{(i)}.
		\end{align}
	\end{subequations}
    \end{itemize}

    The shorthand notations in \eqref{semi_PDE_immediate}-\eqref{semi_PDE_immediate_cor} are defined as:
    \begin{equation*}
	\begin{cases}
		\text{RHS}_E^{(i)}=\sum\limits_{j=1}^{i-1} \left[\hat{a}_{ij}\left(\nabla\cdot(D^{(j)}_r\nabla E^{(j)})-\alpha_0\Delta E^{(j)}\right)+ a_{i j}\left(\alpha_{0} \Delta E^{(j)}+\sigma^{(j)}(B^{(j)}-E^{(j)})\right)\right],
		\\[5pt]
		\text{RHS}_B^{(i)}=\sum\limits_{j=1}^{i-1} \left[\hat{a}_{ij}\left(\mathcal{H}^{(j)}-\beta_0 \Delta B^{(j)}\right)+ a_{i j}\left(\beta_{0} \Delta B^{(j)}+z^3(E^{(j)}-B^{(j)})\right)\right],
		\\[5pt]
		\text{RHS}_Q^{(i)}=\sum\limits_{j=1}^{i-1} \left[\hat{a}_{ij}\left(\nabla \cdot\left(D^{(j)}_r \nabla E^{(j)}\right)+\nabla \cdot\left(D^{(j)}_t \nabla T^{(j)}\right)-\gamma_{0}\Delta Q^{(j)}\right)+ a_{ ij}\gamma_{0} \Delta Q^{(j)}\right].
	\end{cases}
\end{equation*}
For the first stage, we take $\tilde{E}^{(1)}=E^n,\, \tilde{B}^{(1)}=(T^n)^4,\, E^{(1)} = E^n,\, Q^{(1)} = E^n + T^n$.

\subsection{Picard iteration}
For the first-order IMEX scheme \eqref{RDE_2Dplus} and \eqref{RDE_2Dplusc}, or the high-order IMEX scheme in the intermediate stages \eqref{semi_PDE_immediate}-\eqref{semi_PDE_immediate_cor}, each system is mildly nonlinear. Here we will describe how to solve those mildly nonlinear systems with a simple Picard iteration. The iteration does not rely on any specific spatial discretization, so we keep space continuous first. 

Taking the high-order IMEX scheme for the predictor step \eqref{semi_PDE_immediate} as an example, the updating $\tilde{E}$ and $\tilde{B}$ can be rewritten as:
\begin{subequations}\label{semi_Simplified form}
	\begin{align}
		&\left(\frac{1}{\Delta t}-a_{i i}\alpha_{0}\Delta +a_{i i}\tilde{\sigma}^{(i)}\right)\tilde{E}^{(i)}-a_{i i}\tilde{\sigma}^{(i)}\tilde{B}^{(i)}=\text{RHS}_E^{(i)}+\frac{1}{\Delta t}E^n,
		\label{semi_Simplified form_1}\\
		&\left(\frac{1}{\Delta t}-a_{i i}\beta_{0}\Delta +a_{i i}z^3\right)\tilde{B}^{(i)}-a_{i i}z^3 \tilde{E}^{(i)}=\text{RHS}_B^{(i)} + \frac{1}{\Delta t}B^n.
		\label{semi_Simplified form_2}
	\end{align}
\end{subequations} 
For a Picard iteration, starting from the iterative number $l=0$, we set $\tilde{E}^{i,0}=E^{(i-1)}$ and $\tilde{B}^{i,0}=B^{(i-1)}$, and update $\tilde{E}^{i,l+1}$ and $\tilde{B}^{i,l+1}$ from $\tilde{E}^{i,l}$ and $\tilde{B}^{i,l}$ iteratively as:
\begin{subequations}\label{semi_Simplified form_ite}
	\begin{align}
		&\left(\frac{1}{\Delta t}-a_{i i}\alpha_{0}\Delta +a_{i i}\tilde{\sigma}^{i,l}\right)\tilde{E}^{i,l+1}-a_{i i}\tilde{\sigma}^{i,l}\tilde{B}^{i,l+1}=\text{RHS}_E^{(i)}+\frac{1}{\Delta t}E^n,\\
		&\left(\frac{1}{\Delta t}-a_{i i}\beta_{0}\Delta +a_{i i}z^3\right)\tilde{B}^{i,l+1}-a_{i i}z^3 \tilde{E}^{i,l+1}=\text{RHS}_B^{(i)} + \frac{1}{\Delta t}B^n.
	\end{align}
\end{subequations}
As observed, when $\tilde{\sigma}^{i,l}$ is set based on the previous iterative step, \eqref{semi_Simplified form_ite} gives rise to a linear system for $\tilde{E}^{i,l+1}$ and $\tilde{B}^{i,l+1}$ that exhibits diagonal dominance. This property arises from the positivity of $\sigma$ and $z(\bx)$, resulting in rapid convergence when employing iterative methods to solve it.

Similar to \eqref{semi_Simplified form}, the corrector step \eqref{semi_PDE_immediate_cor} for $E$ and $Q$ can be rewritten as:
\begin{subequations}\label{semi_Simplified form_cor}
	\begin{align}
		&\left(\frac{1}{\Delta t}-a_{i i}\alpha_{0}\Delta +a_{i i}\tilde{\sigma}^{(i)}\right)E^{(i)}-a_{i i}\tilde{\sigma}^{(i)}B^{(i)}=\text{RHS}_E^{(i)}+\frac{1}{\Delta t}E^n,
		\label{semi_Simplified form_cor_1}\\
		&\left(\frac{1}{\Delta t}-a_{i i}\gamma_{0}\Delta \right)Q^{(i)}=\text{RHS}_Q^{(i)}+\frac{1}{\Delta t}Q^n.
		\label{semi_Simplified form_cor_2}
	\end{align}
\end{subequations}
Here, $Q^{(i)}$ can be readily obtained by solving the linear system \eqref{semi_Simplified form_cor_2}, which also possesses diagonal dominance. Subsequently, we solve \eqref{semi_Simplified form_cor_1} through a Picard iteration, with initial values $B^{i,0}=\tilde{B}^{(i)}=(\tilde{T}^{(i)})^4$, as follows:
\begin{subequations}\label{semi_iteration form_cor}
	\begin{align}
		&\left(\frac{1}{\Delta t}-a_{i i}\alpha_{0}\Delta +a_{i i}\tilde{\sigma}^{(i)}\right){E}^{i,l+1}=a_{i i}\tilde{\sigma}^{(i)}{B}^{i,l}+\text{RHS}_E^{(i)}+\frac{1}{\Delta t}E^n,\label{semi_iteration form_cor_1}\\
		&B^{i,l+1} = (Q^{(i)} - E^{i,l+1})^4. \label{semi_iteration form_cor_2}
	\end{align}
\end{subequations}
We solve the linear system \eqref{semi_iteration form_cor_1} to get $E^{i,l+1}$, and then update $B^{i,l+1}$ from \eqref{semi_iteration form_cor_2}. Since the initial values are set as $E^{i,0}=\tilde{E}^{(i)}$ and $B^{i,0}=(Q^{(i)}-E^{i,0})^4$ from the predictor step, a fast convergence can be obtained for the corrector step.

The above procedures are similar for the first order IMEX scheme \eqref{RDE_2Dplus} and \eqref{RDE_2Dplusc}, we omit them to save space.
	
	\section{LDG spatial discretization}
	\label{LDG}
	\setcounter{equation}{0}
	\setcounter{figure}{0}
	\setcounter{table}{0}
	
	For the first-order IMEX scheme \eqref{RDE_2Dplus}-\eqref{RDE_2Dplusc}, or the high-order IMEX scheme \eqref{semi_PDE_final_cor}-\eqref{semi_PDE_immediate_cor}, we can couple them with any front capturing spatial discretizations \cite{liu2011FDWENO,bessemoulin2012FV,arbogast2019finite,wang2020local,zhang2022high}. In this work, we utilize an LDG finite element method. The LDG method offers great flexibility for h-p adaptivity and excels at capturing sharp gradient propagations.
	
	\subsection{Some notations}
	For a computational domain $\Omega \subseteq \mathbb{R}^2$, we consider a partition $\mathcal{T}_h$ of $\Omega$ with a set of non-overlapping rectangular elements $\{I_{i,j}\}$, which can cover the whole domain $\Omega$. Here $I_{i j}=I_i \times I_j$, $I_i=\left[x_{i-\frac{1}{2}}, x_{i+\frac{1}{2}}\right]$ and $I_j=\left[y_{j-\frac{1}{2}}, y_{j+\frac{1}{2}}\right]$ for $i=1,2,\cdots,N_x,j=1,2,\cdots,N_y$. We denote the element length and width as $h_i^x=x_{i+\frac{1}{2}}- x_{i-\frac{1}{2}}$, $h_j^y=y_{j+\frac{1}{2}}- y_{j-\frac{1}{2}}$, respectively. $h=\max\limits_{i,j}\{h_i^x,h_j^y\}$ is the maximum edge size of these elements. The center of the element $I_{i,j}$ is $(x_i,y_j)$. We also assume that $\mathcal{T}_h$ is quasi-uniform, namely, $\max\limits_{i}\{h/h_i^x\}$ and $\max\limits_{j}\{h/h_j^y\}$ are upper bounded by a given positive constant.
	
	With the above partition, we follow \cite{xiong2022high} to give some notations which will be used in the following. Given any non-negative integer vector $\mathbf{k}=\left(k_1, k_2\right)$, we define a finite-dimensional discrete piecewise polynomial space as follows
	\begin{equation*}\label{piecewise polynomial space}
		W_{h}^{\mathbf{k}}=\{u \in L^2(\Omega):u|_K \in \mathcal{Q}^{\mathbf{k}}(K), \forall K \in \mathcal{T}_h\},
	\end{equation*}
	where $\mathcal{Q}^{\mathbf{k}}(K)$ consists of tensor product polynomials of degree not exceeding $k_\ell$ along the $\ell$-th direction on each element $K$, for $\ell=1, 2$. Besides, we denote ${\mathbf{W}}_h^{\mathbf{k}}=W_{h}^{\mathbf{k}}\times W_{h}^{\mathbf{k}}$ as a vector space, where each component belongs to $W_{h}^{\mathbf{k}}$.
	We define a unit normal vector $\mathbf{n}^e$ on each edge $e$ of $\mathcal{T}_h$ as follows: if $e \in \partial \Omega$, $\mathbf{n}^e$ is defined as the unit normal vector pointing outside of $\Omega$; for an interior edge $e=\partial K^{+} \cap \partial K^{-}$, the outward unit normal vectors of $e$ taken from the elements $K^{+}$ and $K^{-}$ are denoted by $\mathbf{n}^{+}$ and $\mathbf{n}^{-}$, respectively. Here we fix $\mathbf{n}^e$ as one of $\mathbf{n}^{\pm}$. If we denote $u^{+}$ and $u^{-}$ as the values of a function $u$ on $e$, taken from $K^{+}$ and $K^{-}$ respectively, then the jump $[[u]]$ over an edge $e$ for a scalar-valued function $u$ is defined as
	\begin{equation*}
		[[u]]|_e=-\left(u^{+} \mathbf{n}^{+}+u^{-} \mathbf{n}^{-}\right) \cdot \mathbf{n}^e.
	\end{equation*}
	For a vector-valued function $\mathbf{v}$, the jump $[[\mathbf{v} \cdot \mathbf{n}]]$ is defined as
	\begin{equation*}
		[[\mathbf{v} \cdot \mathbf{n}]]|_e=-\left(\mathbf{v}^{+} \cdot \mathbf{n}^{+}+\mathbf{v}^{-} \cdot \mathbf{n}^{-}\right) \mathbf{n}^e \cdot \mathbf{n}^{+} .
	\end{equation*}
	Accordingly, we express the averages of $u$ and $\mathbf{v} \cdot \mathbf{n}$ as
	\begin{equation*}
		\{\{u\}\}|_e=-\frac{1}{2}\left(u^{+} \mathbf{n}^{+}-u^{-} \mathbf{n}^{-}\right) \cdot \mathbf{n}^e,\quad\{\{\mathbf{v} \cdot \mathbf{n}\}\}|_e=-\frac{1}{2}\left(\mathbf{v}^{+} \cdot \mathbf{n}^{+}-\mathbf{v}^{-} \cdot \mathbf{n}^{-}\right) \mathbf{n}^e \cdot \mathbf{n}^{+}.
	\end{equation*}
	In this work, we take $\mathbf{n}^e=\mathbf{n}^{-}$, then
	\begin{equation*}
		[[u]]|_e=u^{+}-u^{-},\quad  \{\{u\}\}|_e=\frac{1}{2}\left(u^{+}+u^{-}\right),
	\end{equation*}
	and
	\begin{equation*}
		[[\mathbf{v} \cdot \mathbf{n}]]|_e=-(\mathbf{v}^{+} -\mathbf{v}^{-})\cdot\mathbf{n}^{-},\quad\{\{\mathbf{v} \cdot \mathbf{n}\}\}|_e=-\frac{1}{2}\left(\mathbf{v}^{+} +\mathbf{v}^{-}\right)\cdot\mathbf{n}^{-}.
	\end{equation*}
	
	In our implementation, we use an orthogonal basis of $W_h^{\boldsymbol{k}}$ with a uniform rectangular partition $h_x=h_x^i$ for $1\le i\le N_x$ and $h_y=h_y^j$ for $1\le j\le N_y$. In this case, the numerical solution can be expressed as
	\begin{equation*}\label{expand_of_solution}
		u_h(x,y)=\sum\limits_{i=1}^{N_x}\sum\limits_{j=1}^{N_y}\sum\limits_{m=1}^{k_1}\sum\limits_{n=1}^{k_2}u_{i,j}^{m,n}H_i^m(x)H_j^n(y), \quad(x,y)\in \Omega,
	\end{equation*}
	$u_h$ is $E_h$, $T_h$, or $B_h$ respectively. The local basis of $\mathcal{Q}^{\mathbf{k}}(K)$ on each element $K$ is denoted as $H_i^m(x)H_j^n(y)$ for $K=I_{i,j} \in \mathcal{T}_h$, and $\mathbf{k}=(k_1,k_2)$. For example, up to third order, the local bases are chosen as follows
	\begin{align*}
		&H_i^1(x)=1,\, H_i^2(x)=\frac{x-x_i}{h_x},\, H_i^3(x)=\left(\frac{x-x_i}{h_x}\right)^2-\frac{1}{12}, \quad x\in I_i,
		\\
		&H_j^1(y)=1,\, H_j^2(y)=\frac{y-y_j}{h_y},\, H_j^3(y)=\left(\frac{y-y_j}{h_y}\right)^2-\frac{1}{12}, \quad y\in I_j,
	\end{align*}
    with zero extension outside the cell $I_i$ or $I_j$ respectively.
	
	\subsection{First order IMEX-LDG scheme}
	With the above notations, a fully-discrete LDG scheme utilizing a first-order IMEX scheme \eqref{RDE_2Dplus} and \eqref{RDE_2Dplusc} is defined as follows. First, for the predictor step \eqref{RDE_2Dplus}, the scheme reads: given $E_h^{n},T_h^{n},B_h^{n}\in W_{h}^{\mathbf{k}}$ and $\mathbf{p}^n_h$, $\mathbf{q}^n_h$, $\mathbf{r}^n_h \in {\mathbf{W}}_{h}^{\mathbf{k}}$, we find $\tilde{E}_h^{n+1},\tilde{B}_h^{n+1}\in W_{h}^{\mathbf{k}}$ and $\tilde{\mathbf{p}}^{n+1}_h$, $\tilde{\mathbf{q}}^{n+1}_h \in {\mathbf{W}}_{h}^{\mathbf{k}}$, such that for any $\mu$, $\nu$, $\zeta\in W_{h}^{\mathbf{k}}$, and $\pmb{\upsilon}$, $\pmb{\xi} \in {\mathbf{W}}_{h}^{\mathbf{k}}$, such that
	\begin{subequations}
		\label{full_PDE_final_1st}
		\begin{align}
			\frac{1}{\Delta t}\left(\tilde{E}_h^{n+1}-E_h^{n}, \mu\right)&=\mathcal{G}_h^{n}(\mu)-\alpha_{0} \mathcal{L}^{n}_{h,\mathbf{p}}(\mu)+\alpha_{0} \mathcal{L}^{n+1}_{h,\tilde{\mathbf{p}}}(\mu)+\left(\tilde{\sigma}_h^{n+1} (\tilde{B}^{n+1}_h-\tilde{E}^{n+1}_h),\mu\right),
			\\
			\frac{1}{\Delta t}\left(\tilde{B}_h^{n+1}-B_h^{n}, \nu\right)&=\mathcal{H}_h^{n}(\nu)-\beta_{0} \mathcal{L}^{n}_{h,\mathbf{q}}(\nu)+\beta_{0} \mathcal{L}^{n+1}_{h,\tilde{\mathbf{q}}}(\nu)+\left(z^3 (\tilde{B}^{n+1}_h-\tilde{E}^{n+1}_h),\nu\right),
			\\
			(\tilde{T}_h^{n+1},\zeta) &= ((\tilde{B}_h^{n+1})^{1/4},\zeta),\\
			(\tilde{\mathbf{p}}_h^{n+1},\pmb{\upsilon}) &= \pmb{\mathcal{K}}_h(\tilde{E}_h^{n+1},\pmb{\upsilon}),
			\\
			(\tilde{\mathbf{q}}_h^{n+1},\pmb{\xi}) &= \pmb{\mathcal{K}}_h(\tilde{B}_h^{n+1},\pmb{\xi}),
		\end{align}
	\end{subequations}
	with
	\begin{equation*}
		\mathcal{G}_h^{n}(\mu)=\pmb{\mathcal{C}}_h((D_r \mathbf{p}_h)^{n},\mu),
		\quad
		\mathcal{H}_h^{n}(\nu)=\pmb{\mathcal{D}}_h(T_h^{n},\mathbf{r}_h^{n},\nu),
		\quad
	\end{equation*}
	where
	\begin{equation*}
		\mathcal{L}^{n+1}_{h,\tilde{\mathbf{p}}}(\mu) = \pmb{\mathcal{L}}_h(\tilde{\mathbf{p}}_h^{n+1},\mu),
		\quad
		\mathcal{L}^{n+1}_{h,\tilde{\mathbf{q}}}(\nu) = \pmb{\mathcal{L}}_h(\tilde{\mathbf{q}}_h^{n+1},\nu), \quad
		\mathcal{L}^{n}_{h,\mathbf{p}}(\mu) = \pmb{\mathcal{L}}_h(\mathbf{p}_h^{n},\mu),
		\quad
		\mathcal{L}^{n}_{h,\mathbf{q}}(\nu) = \pmb{\mathcal{L}}_h(\mathbf{q}_h^{n},\nu).
	\end{equation*}
	Correspondingly, for the corrector step \eqref{RDE_2Dplusc}, the scheme is defined as: we look for $E_h^{n+1},Q_h^{n+1}\in W_{h}^{\mathbf{k}}$ and $\mathbf{p}^{n+1}_h$, $\mathbf{w}^{n+1}_h \in {\mathbf{W}}_{h}^{\mathbf{k}}$, for any $\mu$, $\phi\in W_{h}^{\mathbf{k}},$ and $\pmb{\upsilon}$, $\pmb{\eta}\in {\mathbf{W}}_{h}^{\mathbf{k}}$, such that
	\begin{subequations}
		\label{full_PDE_final_cor_1st}
		\begin{align}
			\frac{1}{\Delta t}\left(E_h^{n+1}-E_h^{n}, \mu\right)&=\mathcal{G}_h^{n}(\mu)-\alpha_{0} \mathcal{L}^{n}_{h,\mathbf{p}}(\mu)+\alpha_{0} \mathcal{L}^{n+1}_{h,\mathbf{p}}(\mu)+\left(\tilde{\sigma}_h^{n+1} (B^{n+1}_h-E^{n+1}_h),\mu\right) ,
			\\
			\frac{1}{\Delta t}\left(Q_h^{n+1}-Q_h^{n}, \phi\right)&=\mathcal{E}_h^{n}(\phi)-\gamma_{0} \mathcal{L}^{n}_{h,\mathbf{w}}(\phi)+\gamma_{0} \mathcal{L}^{n+1}_{h,\mathbf{w}}(\phi),
			\\
			T^{n+1}_h & = Q^{n+1}_h - E^{n+1}_h, \\
			(\mathbf{p}_h^{n+1},\pmb{\upsilon}) & = \pmb{\mathcal{K}}_h(E_h^{n+1},\pmb{\upsilon}),
			\\  
			(\mathbf{w}_h^{n+1},\pmb{\eta}) &= \pmb{\mathcal{K}}_h(Q_h^{n+1},\pmb{\eta}),
		\end{align}
	\end{subequations}
	with	
	\begin{equation*}
		\mathcal{E}_h^{n}(\phi)=\pmb{\mathcal{F}}_h((D_t \mathbf{r})_h^{n},\phi)+\pmb{\mathcal{C}}_h((D_r \mathbf{p}_h)^{n},\phi),
	\end{equation*}
	where 
	\begin{equation*}
		\mathcal{L}^{n}_{h,\mathbf{p}}(\mu) = \pmb{\mathcal{L}}_h(\mathbf{p}_h^{n},\mu),
		\quad
		\mathcal{L}^{n}_{h,\mathbf{w}}(\phi) = \pmb{\mathcal{L}}_h(\mathbf{w}_h^{n},\phi).
	\end{equation*}
    
	Those operators in \eqref{full_PDE_final_1st}-\eqref{full_PDE_final_cor_1st} are defined as:
	\begin{align*}
		\pmb{\mathcal{L}}_h\left(\mathbf{p}_h,\mu\right)=&-\sum\limits_{K\in\mathcal{T}_h}\int_K \mathbf{p}_h\cdot\nabla\mu d\bx+\sum\limits_{K\in\mathcal{T}_h}\int_{\partial K} \mathbf{n}^e\cdot\hat{\mathbf{p}}_h\mu d s,
		\\[5pt]
		\pmb{\mathcal{K}}_h\left(E_h,\pmb{\upsilon}\right)=&-\sum\limits_{K\in\mathcal{T}_h}\int_K E_h \nabla\cdot\pmb{\upsilon} d\bx+\sum\limits_{K\in\mathcal{T}_h}\int_{\partial K} \mathbf{n}^e\cdot\pmb{\upsilon}\hat{E}_h d s,
		\\[5pt]
		\pmb{\mathcal{C}}_h(D_r \mathbf{p}_h,\mu)=&-\sum\limits_{K\in\mathcal{T}_h}\int_K D_{r,h}\mathbf{p}_h\cdot\nabla\mu d\bx+\sum\limits_{K\in\mathcal{T}_h}\int_{\partial K} \hat{D}_{r,h}\mathbf{n}^e\cdot\hat{\mathbf{p}}_h\mu d s,
		\\[10pt]
		\pmb{\mathcal{F}}_h(D_t \mathbf{r}_h,\phi)=&-\sum\limits_{K\in\mathcal{T}_h}\int_K D_{t,h}\mathbf{r}_h\cdot\nabla\phi d\bx+\sum\limits_{K\in\mathcal{T}_h}\int_{\partial K} \hat{D}_{t,h}\mathbf{n}^e\cdot\hat{\mathbf{r}}_h\phi d s.
		\\[10pt]
		\pmb{\mathcal{D}}_h(T_h,\mathbf{r}_h,\nu)=&-4\kappa\sum\limits_{K\in\mathcal{T}_h}\int_K \left(T^{\frac{11}{2}}_h\mathbf{r}_h\cdot\nabla\nu+3T^{\frac{9}{2}}_h|\mathbf{r}_h|^2\nu \right)d\bx \\ &+4\kappa\sum\limits_{K\in\mathcal{T}_h}\int_{\partial K} \hat{d}_h(T_h;\mathbf{n}^e)\mathbf{n}^e\cdot\hat{\mathbf{r}}_h\nu d s.
	\end{align*}
	In the above and below, $(\cdot,\cdot)$ denotes the usual $L^2$ inner product on $\Omega$. We can find that $\mathbf{p}_h$, $\mathbf{q}_h$, $\mathbf{r}_h$, and $\mathbf{w}_h$ approximate $\nabla E$, $\nabla B$, $\nabla T$, and $\nabla Q$ respectively. 
	
	The variables with $\hat{\cdot}$ are numerical fluxes. For such a diffusion system, we adopt an alternating left-right flux, namely
	\begin{align*}
		&\hat{u}_h = u_h^-,\quad \hspace{14mm} \text{for a scalar-valued function }\, u=E,B,T,Q,
		\\[5pt]
		&\mathbf{n}^e\cdot\hat{\mathbf{v}}_h = \mathbf{n}^-\cdot\mathbf{v}_h^+,\quad \text{for a vector-valued function }\, \mathbf{v}=\mathbf{p},\mathbf{q},\mathbf{r},\mathbf{w}.
	\end{align*}
	The alternating right-left flux or central fluxes as in \cite{castillo2002optimal,jang2014analysis,wang2015stability,xiong2022high} can also be used. As used in \cite{cockburn1998local,xu2010local,yan2002local}, $\hat{D}_{r,h}=\{\{D_{r,h}\}\}|_e$ is defined as a central numerical flux in approximating of $D_{r,h}$, and $\hat{d}_h(T_h;\mathbf{n}^e)$ in $\pmb{\mathcal{D}}_h(T_h,\mathbf{r}_h,\nu)$ is an approximation to $T^{\frac{11}{2}}$,
	\begin{equation*}
		\hat{d}_h(T_h;\mathbf{n}^e)=\left\{\begin{array}{ll}
			[[\frac{2}{13}T_h^{\frac{13}{2}}]]|_e / [[T_h]]|_e, & \text { if } [[T_h]]|_e \neq 0 ,\\[10pt]
			(T_h^{\frac{11}{2}})^{-}, & \text { otherwise. }  
		\end{array}\right.
	\end{equation*}
	$\hat{D}_{t,h}$ is chosen to approximate $D_t$ in the same manner as $\hat{d}_h(T_h;\mathbf{n}^e)$.
	
	\subsection{High order IMEX-LDG scheme}
	By applying a high-order IMEX RK scheme \eqref{IMEX_RK} to \eqref{full_PDE_final_1st}-\eqref{full_PDE_final_cor_1st}, similar to the first-order scheme, the updating of the solutions at $t^{n+1}$ from $t^n$ is given as follows:
		\begin{subequations}
			\label{full_PDE_final_cor}
			\begin{align}
				\frac{1}{\Delta t}\left(E_h^{n+1}-E_h^{n}, \mu\right)&=\sum\limits_{i=1}^{s} \left[\hat{b}_{i}\left(\mathcal{G}_h^{(i)}(\mu)-\alpha_{0} \mathcal{L}^{(i)}_{h,\mathbf{p}}(\mu)\right)+ b_{i }\left(\alpha_{0} \mathcal{L}^{(i)}_{h,\mathbf{p}}(\mu)+\left({\sigma}_h^{(i)} (B^{(i)}_h-E^{(i)}_h),\mu\right)\right)\right],
				\\
				\frac{1}{\Delta t}\left(Q_h^{n+1}-Q_h^{n}, \phi\right)&=\sum\limits_{i=1}^{s} \left[\hat{b}_{i}\left(\mathcal{E}_h^{(i)}(\phi)-\gamma_{0} \mathcal{L}^{(i)}_{h,\mathbf{w}}(\phi)\right)+ b_{i}\gamma_{0} \mathcal{L}^{(i)}_{h,\mathbf{w}}(\phi)\right],				
				\\
				T^{n+1}_h & = Q^{n+1}_h - E^{n+1}_h, \\
				(\mathbf{p}_h^{n+1},\pmb{\upsilon}) &= \pmb{\mathcal{K}}_h(E_h^{n+1},\pmb{\upsilon}),
				\\
				(\mathbf{w}_h^{n+1},\pmb{\eta}) &= \pmb{\mathcal{K}}_h(Q_h^{n+1},\pmb{\eta}).
			\end{align}
		\end{subequations}
	The intermediate stage values for $2\leq i \leq s$ are obtained from:
	\begin{itemize}
		\item {\bf the predictor step: } \\
		\begin{subequations}
			\label{full_PDE_immediate}
			\begin{align}
				\frac{1}{\Delta t}\left(\tilde{E}_h^{(i)}-E_h^{n}, \mu\right)&=\text{RHS}_E^{(i)}(\mu)+a_{i i}\left[\alpha_{0} \mathcal{L}^{(i)}_{h,\tilde{\mathbf{p}}}(\mu)+\left(\tilde{\sigma}_h^{(i)} (\tilde{B}^{(i)}_h-\tilde{E}^{(i)}_h),\mu\right)\right],\label{full_PDE_immediate_1}
				\\
				\frac{1}{\Delta t}\left(\tilde{B}_h^{(i)}-B_h^{n}, \nu\right)&=\text{RHS}_B^{(i)}(\nu)+a_{i i}\left[\beta_{0} \mathcal{L}^{(i)}_{h,\tilde{\mathbf{q}}}(\nu)+\left(z^3 (\tilde{B}^{(i)}_h-\tilde{E}^{(i)}_h),\nu\right)\right],\label{full_PDE_immediate_2}
							\\
				(\tilde{T}_h^{(i)},\zeta) &= ((\tilde{B}_h^{(i)})^{1/4},\zeta), \\
				(\tilde{\mathbf{p}}_h^{(i)},\pmb{\upsilon}) &= \pmb{\mathcal{K}}_h(\tilde{E}_h^{(i)},\pmb{\upsilon}),
				\\
				(\tilde{\mathbf{q}}_h^{(i)},\pmb{\xi}) &= \pmb{\mathcal{K}}_h(\tilde{B}_h^{(i)},\pmb{\xi});
			\end{align}
		\end{subequations}
		\item {\bf the corrector step: } \\
		\begin{subequations}
			\label{full_PDE_immediate_cor}
			\begin{align}
				\frac{1}{\Delta t}\left(E_h^{(i)}-E_h^{n}, \mu\right)&=\text{RHS}_E^{(i)}(\mu)+a_{i i}\left[\alpha_{0} \mathcal{L}^{(i)}_{h,\tilde{\mathbf{p}}}(\mu)+\left(\tilde{\sigma}_h^{(i)} (B^{(i)}_h-E^{(i)}_h),\mu\right)\right],
				\\
				\frac{1}{\Delta t}\left(Q_h^{(i)}-Q_h^{n}, \phi\right)&=\text{RHS}_Q^{(i)}(\phi)+a_{i i}\gamma_{0} \mathcal{L}^{(i)}_{h,\mathbf{w}}(\phi),
				\\
				T^{(i)}_h & = Q^{(i)}_h - E^{(i)}_h, \\
				(\mathbf{p}_h^{(i)},\pmb{\upsilon}) &= \pmb{\mathcal{K}}_h(E_h^{(i)},\pmb{\upsilon}),
				\\
				(\mathbf{w}_h^{(i)},\pmb{\eta}) &= \pmb{\mathcal{K}}_h(Q_h^{(i)},\pmb{\eta}).
			\end{align}
		\end{subequations}
	\end{itemize}
	
	The shorthand notations in \eqref{full_PDE_immediate}-\eqref{full_PDE_immediate_cor} are defined as:
	\begin{equation*}
		\begin{cases}
			\text{RHS}_E^{(i)}(\mu)=\sum\limits_{j=1}^{i-1} \left[\hat{a}_{ij}\left(\mathcal{G}_h^{(j)}(\mu)-\alpha_{0} \mathcal{L}^{(j)}_{h,\mathbf{p}}(\mu)\right)+ a_{i j}\left(\alpha_{0} \mathcal{L}^{(j)}_{h,\mathbf{p}}(\mu)+\left(\tilde{\sigma}_h^{(j)} (B^{(j)}_h-E^{(j)}_h),\mu\right)\right)\right],
			\\[5pt]
			\text{RHS}_B^{(i)}(\nu)=\sum\limits_{j=1}^{i-1} \left[\hat{a}_{ij}\left(\mathcal{H}_h^{(j)}(\nu)-\beta_{0} \mathcal{L}^{(j)}_{h,\mathbf{q}}(\nu)\right)+ a_{i j}\left(\beta_{0} \mathcal{L}^{(j)}_{h,\mathbf{q}}(\nu)+\left(z^3 (\tilde{B}^{(j)}_h-\tilde{E}^{(j)}_h),\nu\right)\right)\right],
			\\[5pt]
			\text{RHS}_Q^{(i)}(\phi)=\sum\limits_{j=1}^{i-1} \left[\hat{a}_{ij}\left(\mathcal{E}_h^{(j)}(\phi)-\gamma_{0} \mathcal{L}^{(j)}_{h,\mathbf{w}}(\phi)\right)+ a_{ ij}\gamma_{0} \mathcal{L}^{(j)}_{h,\mathbf{w}}(\phi)\right].
		\end{cases}
	\end{equation*}
    
    As we can observe, each of the intermediate steps \eqref{full_PDE_immediate}-\eqref{full_PDE_immediate_cor} in a high-order scheme shares the same structure as the first-order scheme \eqref{full_PDE_final_1st}-\eqref{full_PDE_final_cor_1st}. Thus, the first-order scheme provides a general framework for coupling with a multi-stage IMEX RK method, making it convenient to extend to higher orders. Furthermore, the choice of operators and numerical fluxes remains consistent with the first-order scheme.

	\subsection{Picard iteration for a full scheme}
	To represent our scheme as a mildly nonlinear system more effectively, we introduce notations for matrices and vectors. We define a long vector $\boldsymbol{u_c}$ to represent the coefficients of a two-dimensional numerical solution $u_h$, as follows:
	\begin{equation}\label{coe_of_basis}
		\boldsymbol{u_c} = \left(u_{1,1}^{1,1},\cdots,u_{1,1}^{k_1,k_2},u_{2,1}^{1,1},\cdots,u_{2,1}^{k_1,k_2},\cdots,u_{N_x,1}^{1,1},\cdots,u_{N_x,1}^{k_1,k_2},u_{1,2}^{1,1},\cdots,u_{N_x,N_y}^{k_1,k_2}\right)^T.
	\end{equation}
	Such a vector $\boldsymbol{u_c}$ is the value to be updated by a numerical method in the LDG framework.
	For instance, when taking the nonlinear integrals $\left(\sigma_h E_h,\mu\right)$ and $\left(\sigma_h B_h,\mu\right)$ in \eqref{full_PDE_final_1st} and \eqref{full_PDE_final_cor_1st}, practical approximations are employed. These approximations involve Gaussian quadrature integration along each direction, facilitated by the use of an interpolation operator $\mathcal{I}_h$ within each cell $K=I_i \times I_j$, i.e.
	\begin{equation*}
		\left(\sigma_h E_h,\mu\right)_K\approx\int_K\mathcal{I}_h\left(\sigma_h E_h\mu\right)d\bx=\dfrac{h_x h_y}{4}\sum\limits_{G=1}^{k_1+1}\sum\limits_{G'=1}^{k_2+1} \hat{\omega}_G\hat{\omega}_{G'}(\sigma_h E_h\mu)\left(\frac{h_x}{2}\eta_G+x_i,\frac{h_y}{2}\eta_{G'}+y_j\right),
	\end{equation*}
	where $\eta_G$ are Gaussian quadrature points on the interval $[-1,1]$ with corresponding weights $\hat{\omega}_G$. We define $ k_{1,2}=(k_1+1)(k_2+1)$, therefore, by taking $\mu$ over the local bases $H^m_i(x)H^n_j(y)$ in each cell, $(\sigma_h E_h,\mu)$ can be formatted as a matrix $\Lambda_{h,T}$ multiplied by a long vector $\boldsymbol{E_c}$, that is $\Lambda_{h,T} \boldsymbol{E_c}$, 
	where $\Lambda_{h,T}$ is a $k_{1,2}$-block at most tri-diagonal matrix in the size of $N_x N_y$. Here $\sigma$ is $T$-dependent, so that $\Lambda_{h,T}$ depends on $T$, and $\Lambda_{h,T} \boldsymbol{E_c}$ is nonlinear.

	With similar matrix-by-vector notations, the high order fully-discrete IMEX-LDG scheme \eqref{full_PDE_immediate}, corresponding to \eqref{semi_Simplified form}, can be written in the following form
	\begin{subequations}\label{Simplified form}
		\begin{align}
			&\left(\frac{1}{\Delta t}M_h-a_{i i}\alpha_{0}S_h+a_{i i}\Lambda_{h,\tilde{T}}^{(i)}\right)\tilde{\boldsymbol{E}}_{\boldsymbol{c}}^{(i)}-a_{i i}\Lambda_{h,\tilde{T}}^{(i)}\tilde{\boldsymbol{B}}_{\boldsymbol{c}}^{(i)}=F_E^{(i)},
			\label{Simplified form_1}\\
			&\left(\frac{1}{\Delta t}M_h+a_{i i}Z_h-a_{i i}\beta_0S_h\right)\tilde{\boldsymbol{B}}_{\boldsymbol{c}}^{(i)}-a_{i i}Z_h \tilde{\boldsymbol{E}}_{\boldsymbol{c}}^{(i)}=F_B^{(i)},
			\label{Simplified form_2}
		\end{align}
	\end{subequations}
	where $M_h\boldsymbol{E_c}$, $M_h\boldsymbol{B_c}$, $S_h\boldsymbol{E_c}$, $S_h\boldsymbol{B_c}$, $Z_h\boldsymbol{E_c}$ and $Z_h\boldsymbol{B_c}$ are matrix-by-vector forms of $(E_h,\mu)$, $(B_h,\nu)$, $(\mathcal{L}_{h,\mathbf{p}},\mu)$, $(\mathcal{L}_{h,\mathbf{q}},\nu)$, $(z^3 E_h,\nu)$ and $(z^3 B_h, \nu)$, respectively.
	The right-hand-side long vectors $F_E^{(i)}$ and $F_B^{(i)}$ are corresponding to $(\text{RHS}_E^{(i)},\mu)+(E_h,\mu)/\Delta t$ and $(\text{RHS}_B^{(i)},\nu)+(B_h,\nu)/\Delta t$ respectively.

	Similar to \eqref{Simplified form}, the corrector step \eqref{full_PDE_immediate_cor}, corresponding to \eqref{semi_Simplified form_cor}, can be written as
	\begin{subequations}\label{Simplified form_cor}
		\begin{align}
			&\left(\frac{1}{\Delta t}M_h-a_{i i}\alpha_{0}S_h+a_{i i}\Lambda_{h,\tilde{T}}^{(i)}\right)\boldsymbol{E_c}^{(i)}-a_{i i}\Lambda_{h,\tilde{T}}^{(i)}\boldsymbol{B_c}^{(i)}=F_E^{(i)},
			\label{Simplified form_cor_1}\\
			&\left(\frac{1}{\Delta t}M_h-a_{i i}\gamma_{0}S_h\right)\boldsymbol{Q_c}^{(i)}=F^{(i)}_Q+\frac{1}{\Delta t}M_h\boldsymbol{Q_c}^n.
			\label{Simplified form_cor_2}
		\end{align}
	\end{subequations}
	
	For the predictor step, \eqref{Simplified form} can be expressed as a nonlinear algebraic system
	\begin{equation}\label{algebraic system}
		A(\tilde{U}^{(i)})\tilde{U}^{(i)} = F(U^n,U^{(1)},U^{(2)},\ldots,U^{(i-1)}),
	\end{equation}
	with 
	\begin{equation*}
		U=\left [
		\begin {array} {c}
		\boldsymbol{E_c}\\
		\boldsymbol{B_c}
		\end {array}
		\right ],\quad
		F=\left [
		\begin {array} {c}
		F_E^{(i)}\\
		F_B^{(i)}
		\end {array}
		\right ],\quad
		A(U)=\left [
		\begin {array} {cc}
		M_1(U) & M_2(U)\\
		M_3 & M_4
		\end {array}
		\right ].
	\end{equation*}
	The corresponding submatrices in $A(U)$ are given by
	\begin{align*}
		&M_1(U)=\frac{1}{\Delta t}M_h-a_{i i}\alpha_{0}S_h+a_{i i}\Lambda_{h,T},\quad M_2(U)=-a_{i i}\Lambda_{h,T},\\
		&M_3=-a_{i i}Z_h,\quad M_4=\frac{1}{\Delta t}M_h+a_{i i}Z_h-a_{i i}\beta_0S_h.
	\end{align*}
	As we can see, only $M_1(U)$ and $M_2(U)$ mildly depend on $U$ due to $\Lambda_{h,T}$. Hence an iteration procedure, corresponding to \eqref{semi_Simplified form_ite}, can be written as
	\begin{equation}\label{simple Picard iteration}
		A(U^{i,l})U^{i,l+1} = F(U^n,U^{(1)},U^{(2)},\ldots,U^{(i-1)}),
	\end{equation}
	with the initial value taking to be $U^{i,0}=U^{(i-1)}$.
	
	A detailed procedure of the Picard iteration for solving \eqref{simple Picard iteration} is as follows: starting from an iterative number $l=0$, we set $E_h^{i,0}=E_h^{(i-1)}$, $B_h^{i,0}=B_h^{(i-1)}$, and update the unknowns $B_h^{i,l+1}$ and $E_h^{i,l+1}$ iteratively through the following three steps:
	\begin{itemize}
		\item \textbf{Step 1}: compute $A(U^{i,l})$
		\begin{equation*}
			A(U^{i,l})=\left [
			\begin {array} {cc}
			M_1(U^{i,l}) & M_2(U^{i,l})\\
			M_3 & M_4
			\end {array}
			\right ];
		\end{equation*}
		\item \textbf{Step 2}: if $z(\bx)$ is piecewise constant, $M_3$ is an easy invertible diagonal matrix. From \eqref{simple Picard iteration}, by Gauss elimination, we update $E_h^{i,l+1}$ and $B_h^{i,l+1}$ with their long vectors by
		\begin{equation}
			\label{a18}
			\begin{cases}
				\left(M_1{(M_3)}^{-1}M_4-M_2\right)\boldsymbol{B_c}^{i,l+1}=M_1{(M_3)}^{-1}F_B^{(i)}-F_E^{(i)},\\
				M_3 \boldsymbol{E_c}^{i,l+1}=F_B^{(i)}-M_4\boldsymbol{B_c}^{i,l+1}.
			\end{cases}
		\end{equation}
		Otherwise, we solve \eqref{simple Picard iteration} with the GMRES method \cite{saad1986gmres}.
		\item \textbf{Step 3}: we take $\tilde{E}_h^{(i)} = E_h^{i,l+1}$ and $\tilde{B}_h^{(i)} = B_h^{i,l+1}$ when iteration stops. The stop criteria is
		\begin{equation*}
			\|\boldsymbol{E_c}^{i,l+1}-\boldsymbol{E_c}^{i,l}\| < \delta.
		\end{equation*}
		In our numerical tests, we take an $L^2$ norm and a threshold $\delta=10^{-8}$.
	\end{itemize}

	For the corrector step \eqref{Simplified form_cor}, $\boldsymbol{Q_c}^{(i)}$ is solved directly from \eqref{Simplified form_cor_2}, and then obtain $\boldsymbol{E_c}^{(i)}$ and $\boldsymbol{B_c}^{(i)}$ iteratively by
	\begin{subequations}\label{iteration form_cor}
		\begin{align}
			&\left(\frac{1}{\Delta t}M_h-a_{i i}\alpha_{0}S_h+a_{i i}\Lambda_{h,\tilde{T}}^{(i)}\right)\boldsymbol{E_c}^{i,l+1}=a_{i i}\Lambda_{h,\tilde{T}}^{(i)}\boldsymbol{B_c}^{i,l}+F_E^{(i)},\label{iteration form_cor_1}\\
			&\boldsymbol{B_c}^{i,l+1} = (\boldsymbol{Q_c}^{(i)} - \boldsymbol{E_c}^{i,l+1})^4.
		\end{align}
	\end{subequations}
	with the initial values $B_h^{i,0}=\tilde{B}_h^{(i)}$.
	The stop criteria is the same as Step 3 in the predictor step.
	
	For the Picard iteration in the predictor step, we can show that, if $E_h$ and $B_h\in W_h^{(0,0)}$, and the implicit part of the IMEX RK time discretizations has a nonnegative diagonal, the resulting matrix $A(\tilde{U}^{(i)})$ in \eqref{algebraic system} is an $\mathbf{M}$ matrix if $\tilde{U}^{(i)} \geq 0$: 
	\begin{thm}\label{M_matrix}
		{\rm The matrix $A(\tilde{U}^{(i)})$ in \eqref{algebraic system} is an $\mathbf{M}$ matrix for piecewise constant $\mathcal{Q}^0$ finite elements if $\tilde{U}^{(i)} \geq 0$.}
	\end{thm}
	We can easily show that the mass matrix $A(\tilde{U}^{(i)})$ has a nonnegative diagonal, while off-diagonal entries are all non-positive. Besides, $A(U^{(i)})$ is diagonally dominant, so $A(U^{(i)})$ is an $\mathbf{M}$ matrix, namely its inverse has all nonnegative arguments \cite{berman1994nonnegative}. Thanks to this property, rapid convergence can be assured. A similar fast convergence is also observed in the iteration for the correction step \eqref{iteration form_cor}, with $\Lambda_{h,\tilde{T}}$ obtained from the predictor step and initial values taken from $\tilde{B}$. 
	
	\subsection{Positivity preserving and TVB limiters}
	For the solutions of \eqref{RDE}, both the radiation energy $E(\bx,t)$ and the material temperature $T(\bx,t)$ remain positive at all times \cite{guangwei2009progress}. Unfortunately, the schemes described above cannot preserve the positivity of the solutions when they are close to zero. For solutions in $W_{h}^{(0,0)}$, if the right-hand side terms $F_E^{(i)}$ and $F_B^{(i)}$ in \eqref{Simplified form} are non-negative, the solutions without correction preserve positivity. However, this requires a time step size on the order of $\mathcal{O}(h^2)$. Furthermore, an M-matrix cannot be guaranteed for higher-order schemes. In this work, to develop an efficient scheme with larger time steps, we rely on the following procedure to preserve positivity and control numerical oscillations. 
	
	Here we employ a positivity preserving limiter denoted as $P\Pi_h$ from \cite{zhang2009numerical}. The $P\Pi_h$ limiter has also been applied to porous medium equations in \cite{wang2020local}. Additionally, TVB limiters \cite{cockburn1989tvb} are utilized to control numerical oscillations in the case of non-smooth solutions for second and third-order schemes. As mentioned earlier, even for a first-order scheme with solutions in $W^{0,0}_h$, the cell averages may not necessarily remain positive when using large time step sizes. For diffusion systems, a cut-off limiter is commonly employed to ensure a lower bound on the solutions \cite{shi2018functional,xiong2022high,yang2019moving}. 
	
	For a one-dimensional case, we denote $\bar{u}_j$ as the cell average of the numerical solution $u_h$ in the cell $I_j$ and $u_{\text{cold}}$ as the minimum value of $u(\bx,0)$. First, for the limiter $P\Pi_h$, it is applied as follows:
	\begin{enumerate}
		\item  Check the cell average $\bar{u}_j$ in each cell $I_j$, if it is less than $u_{\text{cold}}$, we set 
		\begin{equation*}
			P \Pi_{h} u_{h}= u_{\text {cold }}; 
		\end{equation*}
		
		\item Then if $\bar{u}_j\geq u_{\text {cold}}$, but at least one value is less than $u_{\text {cold}}$ appearing at two endpoints of the cell $I_j$ or at any
		Gaussian point inside the cell $I_j$, we set $P \Pi_{h} u_{h}=P \Pi^1_{h} u_{h}^{1}$, where $u_{h}^{1}$ is a local $L^{2}$-projection of the solution $u_{h}$ to a linear polynomial within this cell, where
		\begin{equation*}
			P \Pi^1_{h}u_{h}^{1}= \begin{cases}
				\left[1-2 h_{j}^{-1}\left(x-x_{j}\right)\right] \bar{u}_{j}+2u_{\text {cold }}h_{j}^{-1}\left(x-x_{j}\right) ,& \text { if } u_{h, j+\frac{1}{2}}^{1,-}<u_{\text {cold}}, \\ \,\\
				\left[1+2 h_{j}^{-1}\left(x-x_{j}\right)\right] \bar{u}_{j}-2u_{\text {cold }}h_{j}^{-1}\left(x-x_{j}\right) ,& \text { if } u_{h, j-\frac{1}{2}}^{1,+}<u_{\text {cold}};
			\end{cases}
		\end{equation*}
		
		\item Otherwise, we set $P \Pi_{h} u_{h}|_{I_{j}}=u_{h}|_{I_{j}}$.
	\end{enumerate}
	
	After applying the PP limiter $P\Pi_h$, a generalized TVB limiter is further applied to control numerical oscillations for second and third order schemes:
	\begin{align*}
		& \begin{cases}u_{j-\frac{1}{2}}^{+,(\bmod)}= \text {TVB minmod }\left(\bar{u}_{j}-u_{j-\frac{1}{2}}^{+}, \frac{\theta}{2}\left(\bar{u}_{j}-\bar{u}_{j-1}\right), \frac{\theta}{2}\left(\bar{u}_{j+1}-\bar{u}_{j}\right),\frac{1}{4}\left(\bar{u}_{j+1}-\bar{u}_{j-1}\right)\right), \\
			u_{j+\frac{1}{2}}^{-,(\bmod)}= \text {TVB minmod }\left(u_{j+\frac{1}{2}}^{-}-\bar{u}_{j}, \frac{\theta}{2}\left(\bar{u}_{j}-\bar{u}_{j-1}\right), \frac{\theta}{2}\left(\bar{u}_{j+1}-\bar{u}_{j}\right),\frac{1}{4}\left(\bar{u}_{j+1}-\bar{u}_{j-1}\right)\right),\end{cases} 
	\end{align*}
	where
	\begin{align*}
		&\text {TVB minmod }(a, b, c, d)= \begin{cases}a, & \text { if }|a| \leq M h^{2}, \\
			\operatorname{minmod}(a, b, c, d), & \text { otherwise, }\end{cases} 
	\end{align*}
	and
	\begin{align*}
		&\operatorname{minmod}(a, b, c, d)= \begin{cases}\operatorname{sign}(a) \min (|a|,|b|,|c|), & \text { if } a b>0, a c>0, \\
			0, & \text { otherwise. }\end{cases}
	\end{align*}
	If $u_{j-\frac{1}{2}}^{+,(\bmod)}$ is different from $u_{j-\frac{1}{2}}^{+}$, or $u_{j+\frac{1}{2}}^{-,(\bmod)}$ is different from $u_{j+\frac{1}{2}}^{-}$, the cell $I_{j}$ is identified as a troubled cell, and the polynomial in this cell is modified to be
	$$
	u_h(x)=\bar{u}_{j}+\left(u_{j+\frac{1}{2}}^{-,(\bmod)}+u_{j-\frac{1}{2}}^{+,(\bmod)}\right) h_{j}^{-1}\left(x-x_{j}\right).
	$$
	In the generalized TVB limiter, we take the parameters $M=0.1$ and $\theta=1.5$ in our numerical tests. 
	
	For the two-dimensional case, we denote $\bar{u}_{i,j}$ as the cell average of the numerical solution $u_h$ in the cell $I_{ij}$ and $u_{\text{cold}}$ as the minimum value of $u_h(\bx,0)$. $u_{i\pm\frac12,j\pm\frac12}=u_h(x_{i\pm\frac12},y_{j\pm\frac12})$ are denoted the values of the endpoints within the cell $I_{i,j}$. The PP limiter $P\Pi_h$ is applied as follows:
	
	\begin{enumerate}
		\item  Check the cell average $\bar{u}_{i,j}$ in each cell $I_{i,j}$, if it is less than $u_{\text{cold}}$, we set 
		\begin{equation*}
			P \Pi_{h} u_{h}= u_{\text {cold }};
		\end{equation*}
		
		\item Then if $\bar{u}_{i,j}\geq u_{\text {cold}}$, but at least one value is less than $u_{\text {cold}}$ appearing at those endpoints of the cell $I_{ij}$ or at any
		Gaussian point inside the cell $I_{ij}$, we set $P \Pi_{h} u_{h}=P \Pi_{h}^1 u_{h}^{1}$, where \[ u_{h}^{1}=\bar{u}_{i,j}+ \alpha_1\psi_x^i+ \alpha_2\psi_y^j, \quad \psi_x^i = \frac{x-x_i}{h_x^i},\,\psi_y^j = \frac{y-y_j}{h_y^j}. \]
		$u_h^1$ is a local $L^{2}$-projection of the solution $u_{h}$ to $\mathcal{P}^1$ within this cell $I_{i,j}$. For $u^1_h$, the extremum occurs only at four endpoints of the cell $I_{i,j}$, and due to $\bar{u}_{i,j}\ge u_{\text{cold}}$ at most two adjacent points are less than $u_{\text {cold}}$. We set those end-point values to be $u_{\text {cold}}$ if they are less than $u_{\text {cold}}$. Next we only need to determine the parameters $\alpha_1$, $\alpha_2$ by using two end-point values, and we take two minimum point values at those end points. They are determined as follows:
		\begin{itemize}
			\item if only one end-point value at $(x_{i-\frac12},y_{j-\frac12})$ is modified, we have
			\begin{align*}
			P\Pi_{h}^1u^1_{h} =&\bar{u}_{i,j}+\left(\bar{u}_{i,j}-u_{\text {cold}}+\frac{u^1_{i+\frac12,j-\frac12}-u^1_{i-\frac12,j+\frac12}}{2}\right)\psi_x^i \\ &+\left(\bar{u}_{i,j}-u_{\text {cold}}-\frac{u^1_{i+\frac12,j-\frac12}-u^1_{i-\frac12,j+\frac12}}{2}\right)\psi_y^j,
		    \end{align*}
			\item otherwise if two end-point values at $(x_{i-\frac12},y_{j-\frac12})$ and $(x_{i-\frac12},y_{j+\frac12})$ are changed, we take
			$$P\Pi_{h}^1u^1_{h}=\bar{u}_{i,j}+2(\bar{u}_{i,j}-u_{\text {cold}})\psi_x^i.$$
		\end{itemize}
	    Other cases can be determined similarly, we omit them to save space;
		
		\item Otherwise, we set $P \Pi_{h} u_{h}|_{I_{ij}}=u_{h}|_{I_{ij}}$.
	\end{enumerate}
	
	For the two-dimensional generalized TVB limiter, we define $w^1(y)=u_h(x_i,y)$ in the cell $I_{i,j}$ and then modify $w^1(y)$ along $x=x_i$ as in the one-dimensional case, that is,
	
	\begin{align*}
		& \begin{cases}u_{j-\frac{1}{2}}^{+,(\bmod)}= \text {TVB minmod}\left(\bar{u}_{i,j}-w^{1}(y_{j-\frac12}^+), \frac{\theta}{2}\left(\bar{u}_{i,j}-\bar{u}_{i,j-1}\right), \frac{\theta}{2}\left(\bar{u}_{i,j+1}-\bar{u}_{i,j}\right),\frac{1}{4}\left(\bar{u}_{i,j+1}-\bar{u}_{i,j-1}\right)\right), \\
			u_{j+\frac{1}{2}}^{-,(\bmod)}= \text {TVB minmod}\left(w^{1}(y_{j+\frac12}^-)-\bar{u}_{i,j}, \frac{\theta}{2}\left(\bar{u}_{i,j}-\bar{u}_{i,j-1}\right), \frac{\theta}{2}\left(\bar{u}_{i,j+1}-\bar{u}_{i,j}\right),\frac{1}{4}\left(\bar{u}_{i,j+1}-\bar{u}_{i,j-1}\right)\right).\end{cases} 
	\end{align*}
	If $u_{j-\frac{1}{2}}^{+,(\bmod)}$ is different from $w^{1}(y_{j-\frac12}^+)$, or $u_{j+\frac{1}{2}}^{-,(\bmod )}$ is different from $w^{1}(y_{j+\frac12}^-)$, the cell $I_{i,j}$ is identified as a troubled cell along $x=x_i$ and $w^{1}(y)$ is modified to be
	$$
	w^{1}(y)=\bar{u}_{i,j}+\left(u_{j+\frac{1}{2}}^{-,(\bmod)}+u_{j-\frac{1}{2}}^{+,(\bmod)}\right) h_{j}^{-1}\left(y-y_{j}\right).
	$$
	Similarly, along $y=y_j$ we can also define $w^2(x)=u_h(x,y_j)$, and
	\begin{align*}
		& \begin{cases}u_{i-\frac{1}{2}}^{+,(\bmod)}= \text {TVB minmod}\left(\bar{u}_{i,j}-w^{2}(x_{i-\frac12}^+), \frac{\theta}{2}\left(\bar{u}_{i,j}-\bar{u}_{i-1,j}\right), \frac{\theta}{2}\left(\bar{u}_{i+1,j}-\bar{u}_{i,j}\right),\frac{1}{4}\left(\bar{u}_{i+1,j}-\bar{u}_{i-1,j}\right)\right), \\
			u_{i+\frac{1}{2}}^{-,(\bmod)}= \text {TVB minmod}\left(w^{2}(x_{i+\frac12}^-)-\bar{u}_{i,j}, \frac{\theta}{2}\left(\bar{u}_{i,j}-\bar{u}_{i-1,j}\right), \frac{\theta}{2}\left(\bar{u}_{i+1,j}-\bar{u}_{i,j}\right),\frac{1}{4}\left(\bar{u}_{i+1,j}-\bar{u}_{i-1,j}\right)\right).\end{cases} 
	\end{align*}
	If $u_{i-\frac{1}{2}}^{+,(\bmod)}\neq w^{2}(x_{i-\frac12}^+)$, or $u_{i+\frac{1}{2}}^{-,(\bmod )}\neq w^{2}(x_{i+\frac12}^-)$, the cell $I_{i,j}$ is identified as a troubled cell along $y=y_j$ and $w^{2}(x)$ is modified to be
	$$
	w^{2}(x)=\bar{u}_{i,j}+\left(u_{i+\frac{1}{2}}^{-,(\bmod)}+u_{i-\frac{1}{2}}^{+,(\bmod)}\right) h_{i}^{-1}\left(x-x_{i}\right).
	$$
	
	If $I_{i,j}$ is a troubled cell along either the $x$ or the $y$ direction, we set $u_h(x,y)$ in the cell $I_{i,j}$ as $u_h(x,y)=w^1(y)+w^2(x)-\bar{u}_{i,j}$.
	
	\subsection{Algorithm flowchart}
	Finally we present our high order IMEX-LDG scheme updating procedure from time level $t_n$ to $t_{n+1}$ in Algorithm \ref{algorithm}.
	\begin{algorithm}
		\caption{The updating from $t^n$ to $t^{n+1}$ with an $s$-stage $(k+1)$-th order IMEX method}\label{algorithm}
		\KwIn{$E^n$, $T^n$}
		\KwOut{$E^{n+1}$, $T^{n+1}$}
		$E^{(1)} = E^n$, $B^{(1)} = (T^n)^4$\;
		\For{$i=2$ \KwTo $s$}{$E^{i,0}=E^{(i-1)}$, $B^{i,0}=B^{(i-1)}$, $l=0$\;
			$M_3=-a_{i i}Z_h,\quad M_4=\frac{1}{\Delta t}M_h+a_{i i}Z_h-a_{i i}\beta_0S_h$\;
			\While{$\|E^{i,l+1}-E^{i,l}\| \ge \delta$ or $l=0$}{$M_1=\frac{1}{\Delta t}M_h-a_{i i}\alpha_{0}S_h+a_{i i}\Lambda_{h,T},\quad M_2=-a_{i i}\Lambda_{h,T}$\;
				update $E^{i,l+1}$ and $B^{i,l+1}$ by \eqref{a18}\;
				apply the limiters to $E^{i,l+1}$ and $B^{i,l+1}$\;
				$l=l+1$\;}
			$\tilde{E}^{(i)} = E^{i,l+1}$, $\tilde{B}^{(i)} = B^{i,l+1}$, $\tilde{T}^{(i)} = (\tilde{B}^{(i)})^{1/4}$  \;
			obtain $Q^{(i)}$ by \eqref{Simplified form_cor_2}\;			
			$E^{i,0}=\tilde{E}^{(i)}$, $l=0$\;
			\While{$\|E^{i,l+1}-E^{i,l}\| \ge \delta$ or $l=0$}{
				update $E^{i,l+1}$ by \eqref{iteration form_cor_1}\;
				apply the limiters to $E^{i,l+1}$\;
				$B^{i,l+1} = (Q^{(i)} - E^{i,l+1})^4$\;
				apply the limiters to $B^{i,l+1}$\;
				$l=l+1$\;
			}
			$E^{(i)} = E^{i,l+1}$, $T^{(i)} = Q^{(i)}-E^{(i)}$\;
		}
		$E^{n+1} = E^{(s)}$, $T^{n+1} = T^{(s)}$\;
		
	\end{algorithm}
	
	\begin{rmk}
		{\rm (Numerical boundary treatment) There are mainly three types of boundary conditions in our numerical experiments: periodic, Neumann, and mixed boundary conditions. The numerical fluxes at the boundary under periodic or Neumann boundary conditions are chosen as in \cite{cockburn1998local}. 
			
		Taking 1D with $\Omega=[x_L,x_R]$ as an example, for a periodic boundary condition, we take $\phi_L^{-}=\phi_R^-$, $\phi_R^+=\phi_L^+$ where $\phi=u$ or $u_x$. For a Neumann boundary condition, we set $\hat{u}_L=u_L^+$, $\hat{u}_R=u_R^-$, while $\hat{v}_L=\hat{v}_R=0$ where $v=u_x$. Here $\phi^{\pm}_L$, $\phi^{\pm}_R$, are the left and right limits at $x_L$ and $x_R$ respectively. Similarly for $u^+_L$ and $u^-_R$.
			
		For a mixed boundary condition, e.g. in the Marshak wave problem in section \ref{sec_examp}, the physical boundary and initial conditions on the domain $x\in\Omega=[x_L,x_R]=[0,1]$ are given as 
			\begin{equation}
				\begin{cases}\label{Marshak_1D_IBC}
					\dfrac{E}{4}-\dfrac{1}{6\sigma}\dfrac{\partial E}{\partial x}=1,\quad  \dfrac{\partial T}{\partial x} = 0, \quad x=0,
					\\[5pt]
					\dfrac{E}{4}+\dfrac{1}{6\sigma}\dfrac{\partial E}{\partial x}=0,\quad \dfrac{\partial T}{\partial x} = 0, \quad x=1,
					\\[5pt]
					E|_{t=0}=1.0\times 10^{-5},\quad T|_{t=0} = E^{1/4}|_{t=0}.
				\end{cases}
			\end{equation}
		Taking the left boundary $x_L=x_\frac12=0$ as an example, we omit the artificial term $\frac{|\nabla E|}{E}$ in $D_r$ and the boundary condition is simplified to \cite{rongpei2012discontinuous,kang2003p,huang2016monotone} 	
			\begin{equation*}
				\dfrac{E}{4}-\dfrac{1}{2}D_r\frac{\partial E}{\partial x}=1.
			\end{equation*}
		Hence, in both the predictor and corrector steps, we set the numerical fluxes as
			\begin{equation*}
				\hat{E}_h = E_h^{+}, \quad\hat{D}_{r,h}\hat{p}_h = \dfrac{\hat{E}_h}{2}-2.
			\end{equation*}
			Similarly at the right boundary $x_{N+\frac{1}{2}}=1$, we take
			\begin{equation*}
				\hat{E}_h = E_h^{-}, \quad\hat{D}_{r,h}\hat{p}_h = -\dfrac{\hat{E}_h}{2}.
			\end{equation*}
			In the corrector step, we set 
			\begin{equation*}
				\hat{Q}_h =\hat{E}_h + \hat{T}_h , \quad\hat{w}_h = 3\sigma(\hat{T}_h)\hat{D}_{r,h}\hat{p}_h.
			\end{equation*}
		} 
	\end{rmk}
	
	\section{Numerical Examples}
	\label{sec_examp}
	\setcounter{equation}{0}
	\setcounter{figure}{0}
	\setcounter{table}{0}
	
	In this section, we conduct numerical experiments to validate the high-order accuracy, conservation properties, suitability for large time steps, and effectiveness in capturing sharp fronts in both homogeneous and heterogeneous media using our proposed schemes.
	For the added diffusion terms, we take the coefficients as $ \alpha_0 = \tau\max\{\frac{1}{3\sigma}\}$, $\beta_0 = \tau\max\{D_t\}$ and $\gamma_0 = \alpha_0 + \beta_0$ with $\tau=0.6$ \cite{wang2020local}. In space, the $\mathcal{Q}^{\boldsymbol{k}}$ basis with $k$-th piecewise polynomial in each direction for $k = 0, 1, 2$ is taken. Correspondingly, an $s$-stage $(k+1)$-th order globally stiffly accurate IMEX RK time discretization is employed in time with the double Butcher tableau given in \ref{IMEX Butcher tableau}. In the following, our schemes are denoted as $(k+1)$-th order methods for $k = 0, 1, 2$ with $s=2, 3, 5$, respectively.
		
	\begin{example}
		\label{test1D_order}
		{(\bf Accuracy test in 1D)} {\rm First we consider a 1D example with smooth initial values and periodic boundary conditions at the equilibrium, which are given by}
		\begin{equation}
			\label{a23}
			\begin{cases}
				T(x,0)=0.8+0.1\sin(x),
				\\[5pt]
				E(x,0)=T(x,0)^4, 
			\end{cases}
		\end{equation}
	{\rm
		on the computational domain $[-\pi,\pi]$. We take $\kappa=0.1$ and a homogeneous medium $z(x)=1$. The problem is run to time $t = 5$ using the 1st, 2nd, and 3rd order methods, respectively. Since the exact solution is not available, we compute the numerical errors by comparing numerical solutions with a reference solution. In this case, the source term is not stiff as $T$ is away from $0$, so we compute the reference solution by a 3rd order LDG method with a 3rd order explicit strong-stability-preserving RK time discretization \cite{gottlieb1998total}, on a much refined mesh $N = 1024$. In Tables \ref{table_1}-\ref{table_3}, we show the numerical $L^2$ errors and orders of accuracy for the 1st, 2nd, and 3rd order schemes with different time steps, respectively. From these tables, we can see that our methods achieve the corresponding orders of accuracy when the time step is $\Delta t = \mathcal{O}(h).$ From the numerical results, we find that larger $\Delta t$ lead to larger errors. In Table \ref{table_3}, for the 3rd order method, an order reduction can be observed for a large ratio of $\Delta t/h$. The order increases with further mesh refinements. Such a phenomeno may be due to a high order IMEX time discretization. }
		\begin{table}[!htbp]
			\centering
			\caption{The numerical $L^2$ errors and orders of accuracy for the 1st order scheme with different time steps for Example \ref{test1D_order}. $t = 5$. }\label{table_1}
			\begin{tabular}{cccccccc}
				\toprule
				\multirow{2}{*}{N} & \multirow{2}{*}{$\Delta t$}           & \multicolumn{1}{c}{$L^2$ error} & order                     & \multicolumn{1}{c}{$L^2$ error} & order                     & \multicolumn{1}{c}{$L^2$ error} & order                     \\ \cline{3-8} 
				&                                       & \multicolumn{2}{c}{E}                                       & \multicolumn{2}{c}{T}                                       & \multicolumn{2}{c}{B}                                       \\ \hline
				4   & \multirow{6}{*}{$\frac{h}{2}$} & 3.46e-2 & -    & 1.74e-2 & -    & 3.51e-2 & -    \\
				8   &                      & 1.73e-2 & 1.00 & 8.65e-3 & 1.01 & 1.76e-2 & 0.99 \\
				16  &                      & 8.68e-3 & 1.00 & 4.34e-3 & 0.99 & 8.84e-3 & 1.00 \\
				32  &                      & 4.34e-3 & 1.00 & 2.17e-3 & 1.00 & 4.42e-3 & 1.00 \\
				64  &                      & 2.17e-3 & 1.00 & 1.09e-3 & 1.00 & 2.21e-3 & 1.00 \\
				128 &                      & 1.08e-3 & 1.00 & 5.43e-4 & 1.00 & 1.10e-3 & 1.00 \\ \hline
				4   & \multirow{6}{*}{$h$} & 3.47e-2 & -    & 1.75e-2 & -    & 3.52e-2 & -    \\
				8   &                      & 1.74e-2 & 1.00 & 8.68e-3 & 1.01 & 1.77e-2 & 0.99 \\
				16  &                      & 8.70e-3 & 1.00 & 4.36e-3 & 1.00 & 8.86e-3 & 1.00 \\
				32  &                      & 4.35e-3 & 1.00 & 2.18e-3 & 1.00 & 4.43e-3 & 1.00 \\
				64  &                      & 2.17e-3 & 1.00 & 1.09e-3 & 1.00 & 2.21e-3 & 1.00 \\
				128 &                      & 1.09e-3 & 1.00 & 5.44e-4 & 1.00 & 1.11e-3 & 1.00 \\ \hline
				4   & \multirow{6}{*}{$3h$} & 3.52e-2 & -    & 1.77e-2 & -    & 3.57e-2 & -    \\
				8   &                      & 1.77e-2 & 0.99 & 8.85e-3 & 1.00 & 1.80e-2 & 0.99 \\
				16  &                      & 8.86e-3 & 1.00 & 4.44e-3 & 1.00 & 9.02e-3 & 1.00 \\
				32  &                      & 4.43e-3 & 1.00 & 2.22e-3 & 1.00 & 4.51e-3 & 1.00 \\
				64  &                      & 2.21e-3 & 1.00 & 1.11e-3 & 1.00 & 2.25e-3 & 1.00 \\
				128 &                      & 1.11e-3 & 1.00 & 5.55e-4 & 1.00 & 1.13e-3 & 1.00 \\ \hline
				4   & \multirow{6}{*}{$5h$} & 3.53e-2 & -    & 1.78e-2 & -    & 3.58e-2 & -    \\
				8   &                      & 1.79e-2 & 0.98 & 9.00e-3 & 0.98 & 1.82e-2 & 0.97 \\
				16  &                      & 9.04e-3 & 0.99 & 4.55e-3 & 0.98 & 9.21e-3 & 0.99 \\
				32  &                      & 4.56e-3 & 0.99 & 2.30e-3 & 0.99 & 4.64e-3 & 0.99 \\
				64  &                      & 2.28e-3 & 1.00 & 1.15e-3 & 1.00 & 2.32e-3 & 1.00 \\
				128 &                      & 1.14e-3 & 1.00 & 5.74e-4 & 1.00 & 1.16e-3 & 1.00 \\ \bottomrule
			\end{tabular}
		\end{table}
		
		\begin{table}[!htbp]
			\centering
			\caption{The numerical $L^2$ errors and orders of accuracy for the 2nd order scheme with different time steps for Example \ref{test1D_order}. $t=5.$}\label{table_2}
			\begin{tabular}{cccccccc}
				\toprule
				\multirow{2}{*}{N} & \multirow{2}{*}{$\Delta t$}           & \multicolumn{1}{c}{$L^2$ error} & order & \multicolumn{1}{c}{$L^2$ error} & order & \multicolumn{1}{c}{$L^2$ error} & order \\ \cline{3-8} 
				&                                       & \multicolumn{2}{c}{E}                   & \multicolumn{2}{c}{T}                   & \multicolumn{2}{c}{B}                   \\ \hline
				4   & \multirow{6}{*}{$\frac{h}{2}$} & 1.25e-2 & -    & 6.18e-3 & -    & 1.27e-2 & -    \\
				8   &                       & 3.00e-3 & 2.06 & 1.57e-3 & 1.98 & 3.04e-3 & 2.06 \\
				16  &                       & 7.36e-4 & 2.02 & 3.86e-4 & 2.02 & 7.48e-4 & 2.02 \\
				32  &                       & 1.83e-4 & 2.01 & 9.60e-5 & 2.01 & 1.86e-4 & 2.01 \\
				64  &                       & 4.58e-5 & 2.00 & 2.40e-5 & 2.00 & 4.66e-5 & 2.00 \\
				128 &                       & 1.14e-5 & 2.00 & 5.99e-6 & 2.00 & 1.16e-5 & 2.00  \\ \hline
				4   & \multirow{6}{*}{$h$} & 1.25e-2 & -    & 6.19e-3 & -    & 1.27e-2 & -    \\
				8   &                      & 3.00e-3 & 2.06 & 1.57e-3 & 1.98 & 3.04e-3 & 2.06 \\
				16  &                      & 7.36e-4 & 2.02 & 3.86e-4 & 2.02 & 7.48e-4 & 2.02 \\
				32  &                      & 1.83e-4 & 2.01 & 9.60e-5 & 2.01 & 1.86e-4 & 2.01 \\
				64  &                      & 4.58e-5 & 2.00 & 2.40e-5 & 2.00 & 4.66e-5 & 2.00 \\
				128 &                      & 1.14e-5 & 2.00 & 5.99e-6 & 2.00 & 1.16e-5 & 2.00  \\ \hline
				4   & \multirow{6}{*}{$3h$} & 1.25e-2 & -    & 6.29e-3 & -    & 1.27e-2 & -    \\
				8   &                      & 3.03e-3 & 2.05 & 1.60e-3 & 1.97 & 3.10e-3 & 2.04 \\
				16  &                      & 7.41e-4 & 2.03 & 3.90e-4 & 2.04 & 7.58e-4 & 2.03 \\
				32  &                      & 1.85e-4 & 2.01 & 9.63e-5 & 2.02 & 1.87e-4 & 2.02 \\
				64  &                      & 4.61e-5 & 2.00 & 2.41e-5 & 2.00 & 4.68e-5 & 2.00 \\
				128 &                      & 1.15e-5 & 2.00 & 6.03e-6 & 2.00 & 1.17e-5 & 2.00  \\ \hline
				4   & \multirow{6}{*}{$5h$} & 1.30e-2 & -    & 6.29e-3 & -    & 1.27e-2 & -    \\
				8   &                      & 3.18e-3 & 2.03 & 1.93e-3 & 1.71 & 3.81e-3 & 1.74 \\
				16  &                      & 9.11e-4 & 1.80 & 4.37e-4 & 2.14 & 8.77e-4 & 2.12 \\
				32  &                      & 2.00e-4 & 2.19 & 9.99e-5 & 2.13 & 1.93e-4 & 2.18 \\
				64  &                      & 4.81e-5 & 2.05 & 2.52e-5 & 1.99 & 4.92e-5 & 1.97 \\
				128 &                      & 1.20e-5 & 2.01 & 6.28e-6 & 2.00 & 1.22e-5 & 2.00  \\ \bottomrule
			\end{tabular}
		\end{table}
		
		\begin{table}[!htbp]
			\centering
			\caption{The numerical $L^2$ errors and orders of accuracy for the 3rd order scheme with different time steps for Example \ref{test1D_order}. $t=5.$}\label{table_3}
			\begin{tabular}{cccccccc}
				\toprule
				\multirow{2}{*}{N} & \multirow{2}{*}{$\Delta t$}           & \multicolumn{1}{c}{$L^2$ error} & order & \multicolumn{1}{c}{$L^2$ error} & order & \multicolumn{1}{c}{$L^2$ error} & order \\ \cline{3-8} 
				&                                       & \multicolumn{2}{c}{E}                   & \multicolumn{2}{c}{T}                   & \multicolumn{2}{c}{B}                   \\ \hline
				4   & \multirow{6}{*}{$\frac{h}{2}$} & 1.71e-3 & -    & 1.03e-3 & -    & 1.78e-3 & -    \\
				8   &                      & 1.97e-4 & 3.12 & 1.14e-4 & 3.18 & 2.04e-4 & 3.12 \\
				16  &                      & 2.49e-5 & 2.99 & 1.47e-5 & 2.96 & 2.61e-5 & 2.97 \\
				32  &                      & 3.11e-6 & 3.00 & 1.86e-6 & 2.98 & 3.31e-6 & 2.98 \\
				64  &                      & 3.89e-7 & 3.00 & 2.34e-7 & 2.99 & 4.15e-7 & 2.99 \\
				128 &                      & 4.86e-8 & 3.00 & 2.93e-8 & 3.00 & 5.20e-8 & 3.00  \\ \hline
				4   & \multirow{6}{*}{$h$} & 1.72e-3 & -    & 1.04e-3 & -    & 1.78e-3 & -    \\
				8   &                      & 1.97e-4 & 3.12 & 1.14e-4 & 3.18 & 2.04e-4 & 3.12 \\
				16  &                      & 2.49e-5 & 2.99 & 1.47e-5 & 2.96 & 2.62e-5 & 2.96 \\
				32  &                      & 3.12e-6 & 3.00 & 1.87e-6 & 2.97 & 3.33e-6 & 2.98 \\
				64  &                      & 3.91e-7 & 3.00 & 2.36e-7 & 2.98 & 4.21e-7 & 2.98 \\
				128 &                      & 4.90e-8 & 3.00 & 2.97e-8 & 2.99 & 5.29e-8 & 2.99  \\ \hline
				4   & \multirow{8}{*}{3$h$} & 1.98e-3 & -    & 1.15e-3 & -    & 2.06e-3 & -    \\
				8   &                       & 2.21e-4 & 3.17 & 1.28e-4 & 3.17 & 2.36e-4 & 3.12 \\
				16  &                       & 2.99e-5 & 2.89 & 1.78e-5 & 2.85 & 3.30e-5 & 2.84 \\
				32  &                       & 4.28e-6 & 2.80 & 2.81e-6 & 2.66 & 5.30e-6 & 2.64 \\
				64  &                       & 7.37e-7 & 2.54 & 6.50e-7 & 2.11 & 1.31e-6 & 2.01 \\
				128 &                       & 1.25e-7 & 2.56 & 1.07e-7 & 2.61 & 2.16e-7 & 2.60 \\
				256 &                       & 1.97e-8 & 2.66 & 1.60e-8 & 2.74 & 3.23e-8 & 2.74 \\
				512 &                       & 2.85e-9 & 2.79 & 2.21e-9 & 2.85 & 4.48e-9 & 2.85  \\ \hline
				4   & \multirow{8}{*}{5$h$} & 2.13e-3 & -    & 1.19e-3 & -    & 2.20e-3 & -    \\
				8   &                       & 4.29e-4 & 2.31 & 2.37e-4 & 2.32 & 4.82e-4 & 2.19 \\
				16  &                       & 6.43e-5 & 2.74 & 3.65e-5 & 2.70 & 7.39e-5 & 2.71 \\
				32  &                       & 1.06e-5 & 2.61 & 8.27e-6 & 2.14 & 1.68e-5 & 2.13 \\
				64  &                       & 1.94e-6 & 2.44 & 1.64e-6 & 2.33 & 3.31e-6 & 2.35 \\
				128 &                       & 3.70e-7 & 2.39 & 3.09e-7 & 2.41 & 6.18e-7 & 2.42 \\
				256 &                       & 6.82e-8 & 2.44 & 5.69e-8 & 2.44 & 1.15e-7 & 2.43 \\
				512 &                       & 1.11e-8 & 2.61 & 8.99e-9 & 2.66 & 1.82e-8 & 2.66  \\ \bottomrule
			\end{tabular}
		\end{table}
		
	\end{example}

\begin{figure*}[!htbp]
		\begin{minipage}[t]{0.49\linewidth}
			\centerline{\includegraphics[scale=0.57]{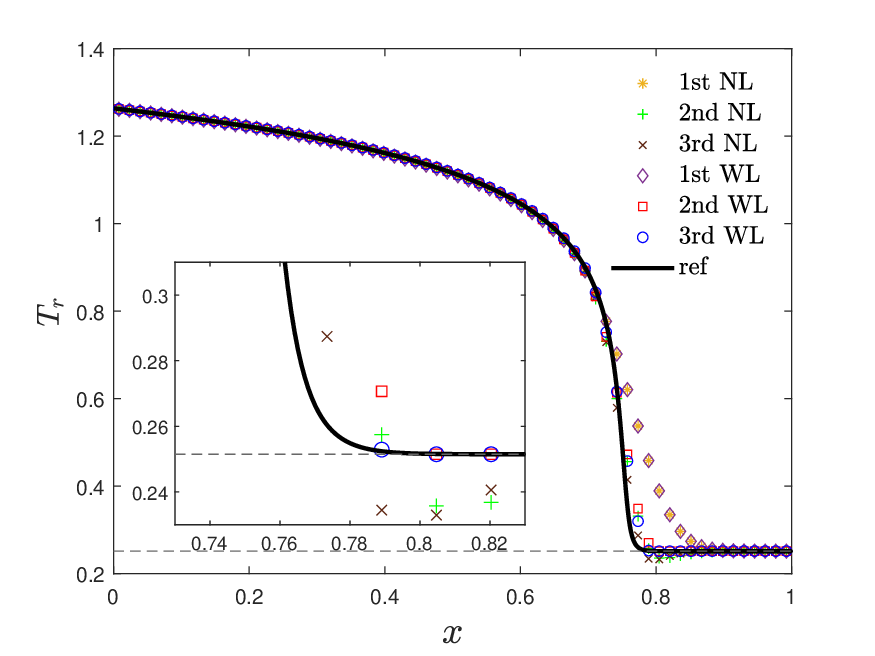}}
		\end{minipage}
		\begin{minipage}[t]{0.49\linewidth}
			\centerline{\includegraphics[scale=0.57]{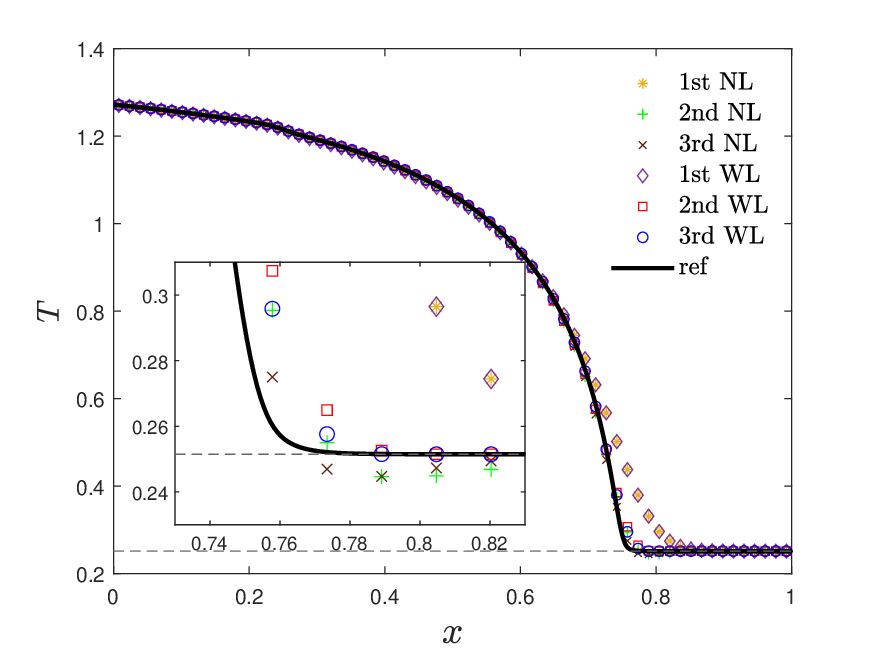}}
		\end{minipage}
		\vskip 0.1mm
		\begin{minipage}[t]{0.49\linewidth}
			\centerline{\includegraphics[scale=0.57]{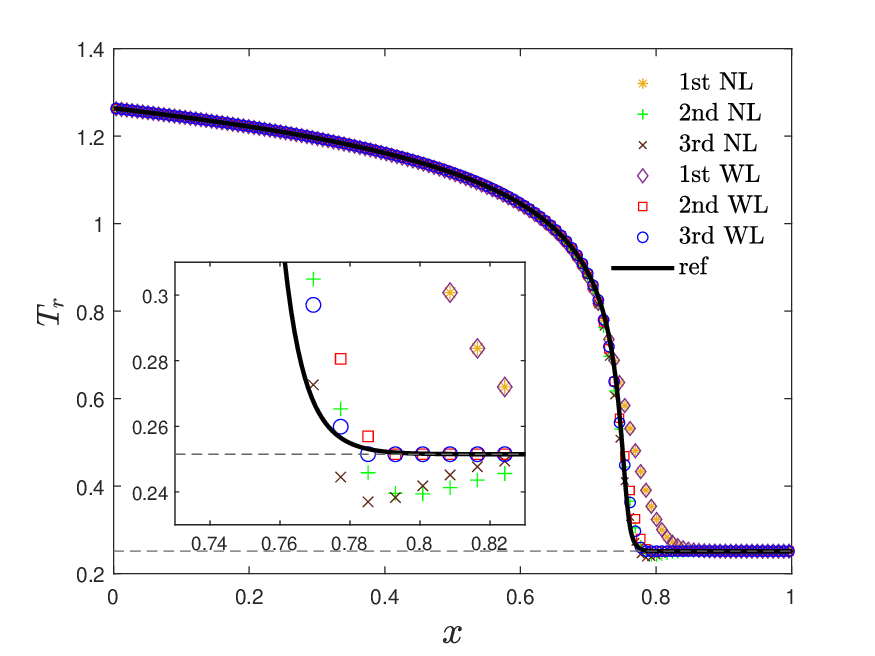}}
		\end{minipage}
		\begin{minipage}[t]{0.49\linewidth}
			\centerline{\includegraphics[scale=0.57]{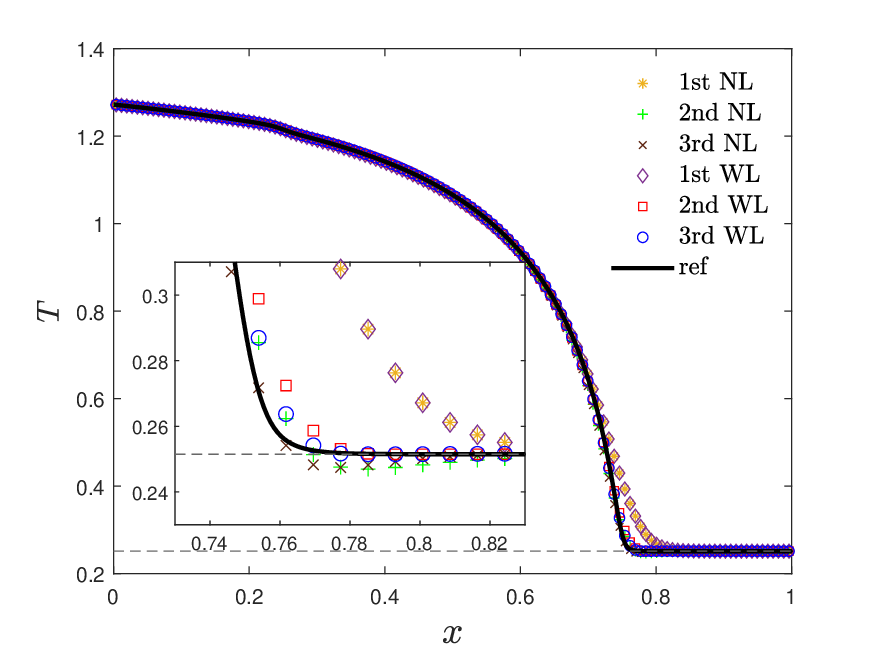}}
		\end{minipage}
		\caption{The numerical results for Example \ref{test_limiters} with $\kappa=0$. 1st, 2nd, and 3rd order schemes are used. ``WL’': with limiters, ``NL'’: without limiters. ``ref'' represents the reference solutions. Top row: $N=64$; Bottom row: $N=128$. Left column: radiation temperature $T_r$; Right column: material temperature $T$. $t=1$, $\Delta t = \frac{h}{5}.$
		}\label{fig1_1}
	\end{figure*}

\begin{figure*}[!htbp]
		\begin{minipage}[t]{0.49\linewidth}
			\centerline{\includegraphics[scale=0.57]{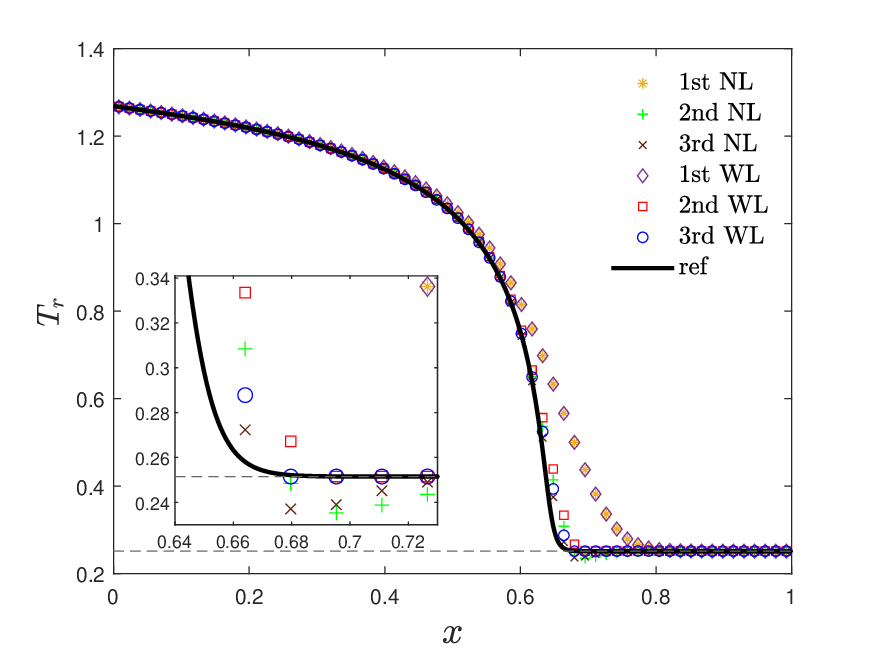}}
		\end{minipage}
		\begin{minipage}[t]{0.49\linewidth}
			\centerline{\includegraphics[scale=0.57]{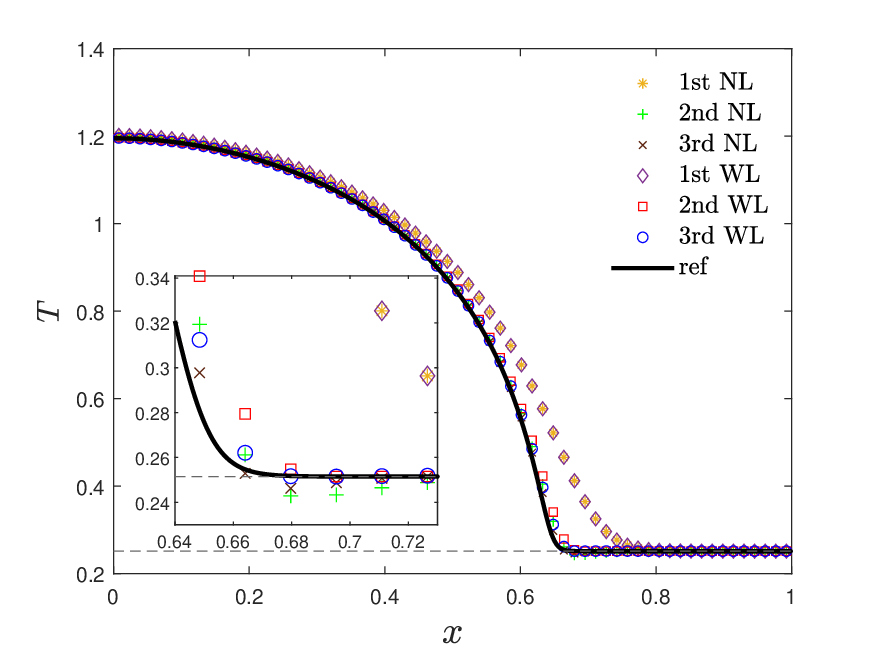}}
		\end{minipage}
		\vskip 0.1mm
		\begin{minipage}[t]{0.49\linewidth}
			\centerline{\includegraphics[scale=0.57]{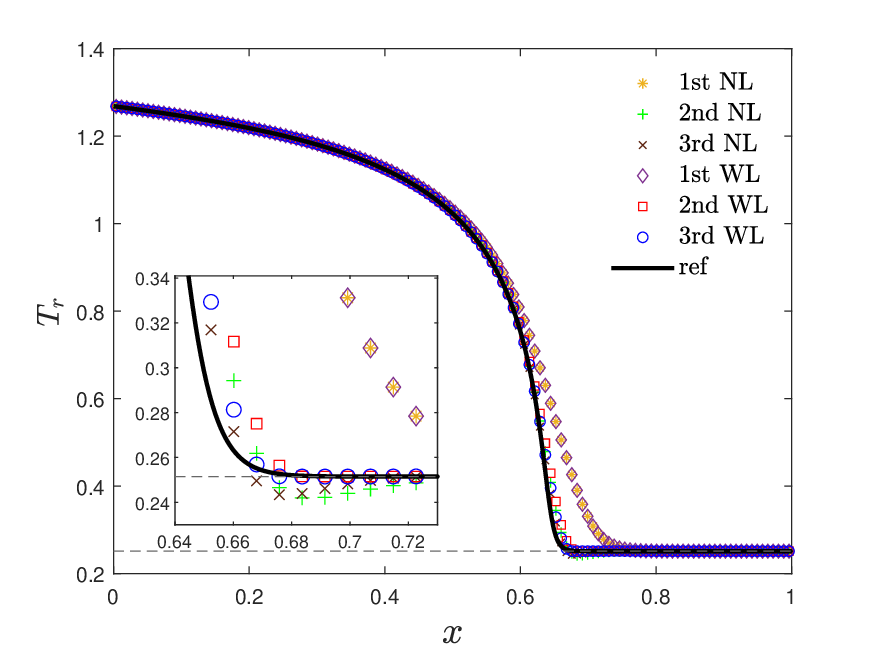}}
		\end{minipage}
		\begin{minipage}[t]{0.49\linewidth}
			\centerline{\includegraphics[scale=0.57]{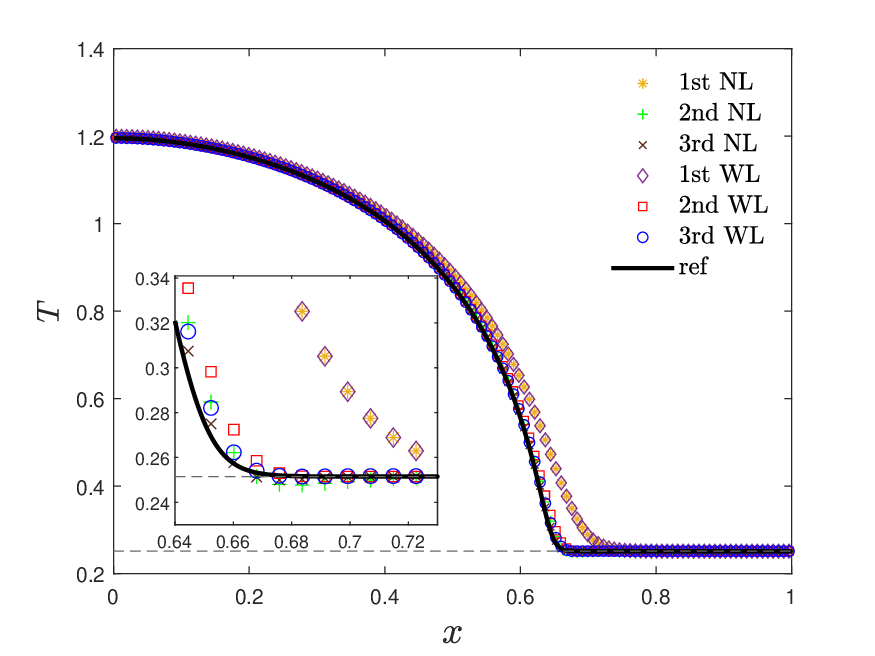}}
		\end{minipage}
		\caption{The numerical results for Example \ref{test_limiters} with $\kappa=0.1,$ and $N=128$. 1st, 2nd, and 3rd order schemes are used. ``WL’' : with limiters, ``NL'’: without limiters. ``ref'' represents the reference solutions. Top row: $N=64$; Bottom row: $N=128$. Left column: radiation temperature $T_r$; Right column: material temperature $T$. $t=0.5$, $\Delta t = \frac{h}{5}.$
		}\label{fig1_2}
	\end{figure*}

	\begin{example}\label{test_limiters} 
		{\rm Next we consider the following initial boundary data \cite{lowrie2004comparison}, with a sharp transition in the initial values}
		\begin{equation}
			\begin{cases}
				E(x,0)=E_L+(E_R-E_L)\dfrac{1+\tanh [50(x-0.25)]}{2},
				\\[5pt]
				T(x,0) = E(x,0)^{1/4},
				\\[5pt]
				\dfrac{E}{4}-\dfrac{1}{6\sigma}\dfrac{\partial E}{\partial x}=1 , \, \, \, \, \dfrac{\partial T}{\partial x} = 0, x=0,
				\\[5pt]
				\dfrac{E}{4}+\dfrac{1}{6\sigma}\dfrac{\partial E}{\partial x}=V_R , \dfrac{\partial T}{\partial x} = 0, x=1,
			\end{cases}
		\end{equation}
		{\rm where $V_R=1\times 10^{-3}$ , $E_L=4$ and $E_R=4\times 10^{-3}.$
		We consider two cases $\kappa = 0$ and $\kappa =0.1$ in $D_t$. For the case of $\kappa=0$, the system of $E$ and $B$ itself is in a conservative form. We take $z(x)=1$ and a time step size $\Delta t = \frac{1}{5}h$. The results at $t=1$ for $\kappa=0$ and $t=0.5$ for $\kappa=0.1$ are shown in Fig.~\ref{fig1_1} and \ref{fig1_2}, respectively. Reference solutions are obtained by a 1st order explicit RK LDG method on $1024$ elements with a small enough time step, which is denoted as ``ref''. The left column is the radiation temperature $T_r$ and the right column is the material temperature $T$. Two different mesh sizes are considered, $N=64$ and $N=128$. ``WL'' refers to numerical solutions obtained with limiters, and correspondingly ``NL'' refers to numerical solutions without limiters. It can be observed that the higher the order, the closer the numerical solutions are, as compared to the reference solutions. In addition, the 1st order solutions perform well without limiters, but the 2nd order and 3rd order solutions have small oscillations in front of the sharp gradient. After applying limiters, these oscillations can be well controlled. As we can see from Fig.~\ref{fig1_1} and Fig.~\ref{fig1_2}, refining the mesh from $N = 64$ to $N = 128$, the numerical solutions match the reference solutions better, no matter with or without limiters. Especially, deviations due to the application of the limiters are also reduced with the mesh refinement. 
		}
	\end{example}
	
\begin{figure*}[!htbp]
		\begin{minipage}[t]{0.49\linewidth}
			\centerline{\includegraphics[scale=0.57]{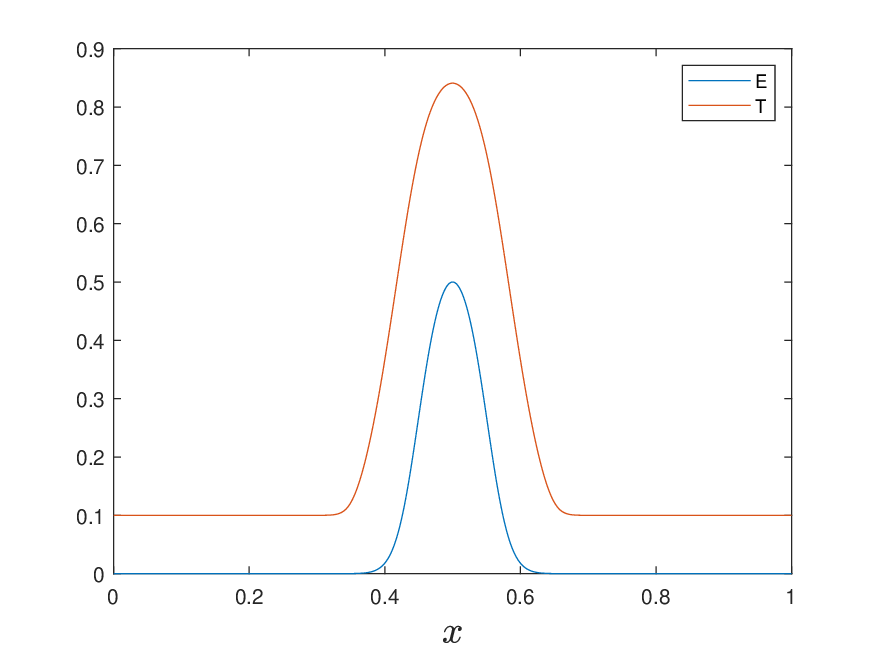}}
		\end{minipage}
		\begin{minipage}[t]{0.49\linewidth}
			\centerline{\includegraphics[scale=0.57]{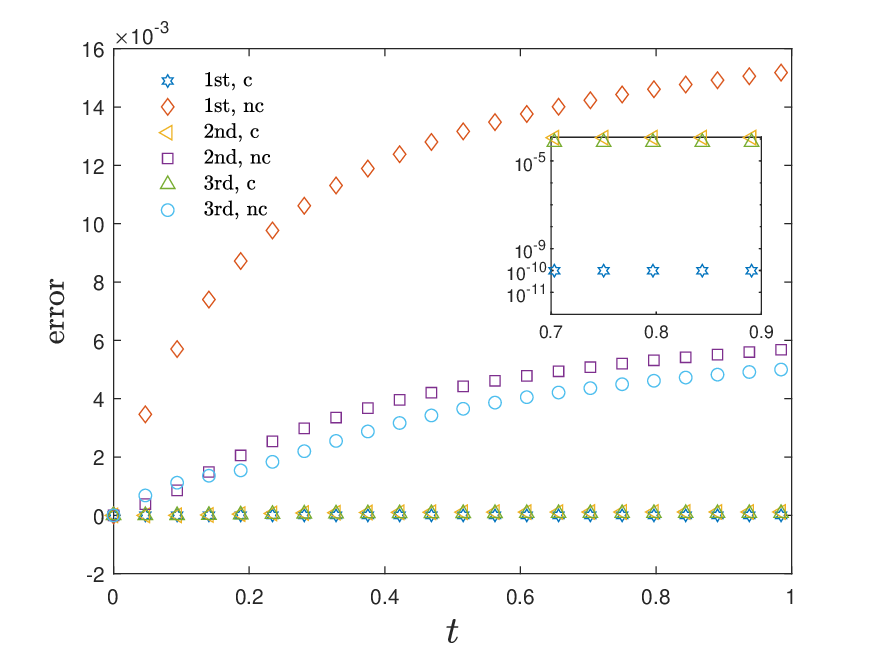}}
		\end{minipage}
		\caption{The numerical results for Example \ref{exp_conservation} with $\kappa=0.5$ and $N=128$. 1st, 2nd, and 3rd order schemes are used. ``c'': with a conservation corrector step; ``nc'': without a corrector step. Left: initial values; Right: time evolution of conservation errors. $\Delta t = \frac{1}{30}h$ for ``2nd, nc'', $\Delta t = \frac{1}{50}h$ for ``3rd, nc'' and $\Delta t = \frac{1}{5}h$ for all others. 
		}\label{figcon_ini}
	\end{figure*}

\begin{figure*}[!htbp]
		\begin{minipage}[t]{0.49\linewidth}
			\centerline{\includegraphics[scale=0.57]{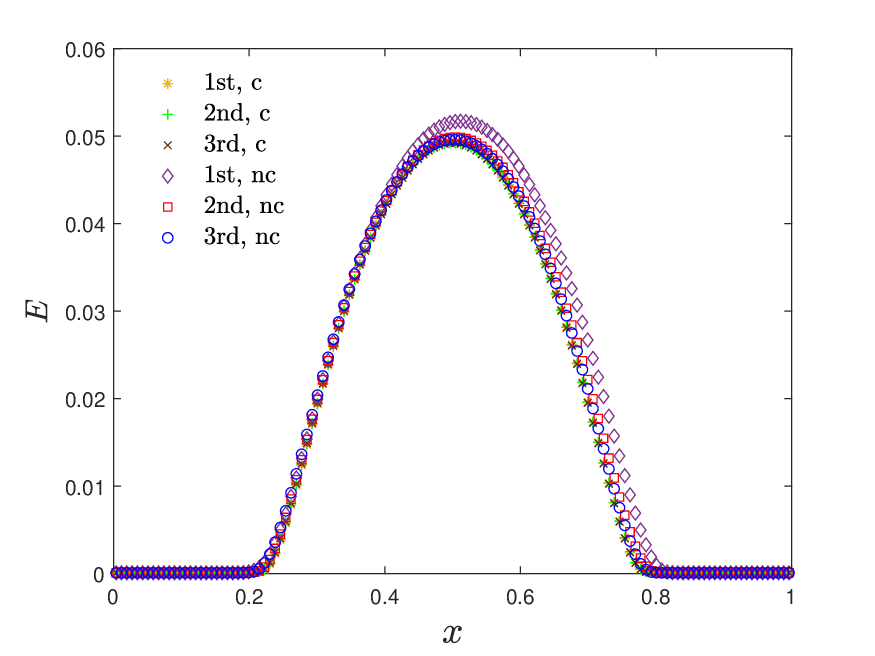}}
		\end{minipage}
		\begin{minipage}[t]{0.49\linewidth}
			\centerline{\includegraphics[scale=0.57]{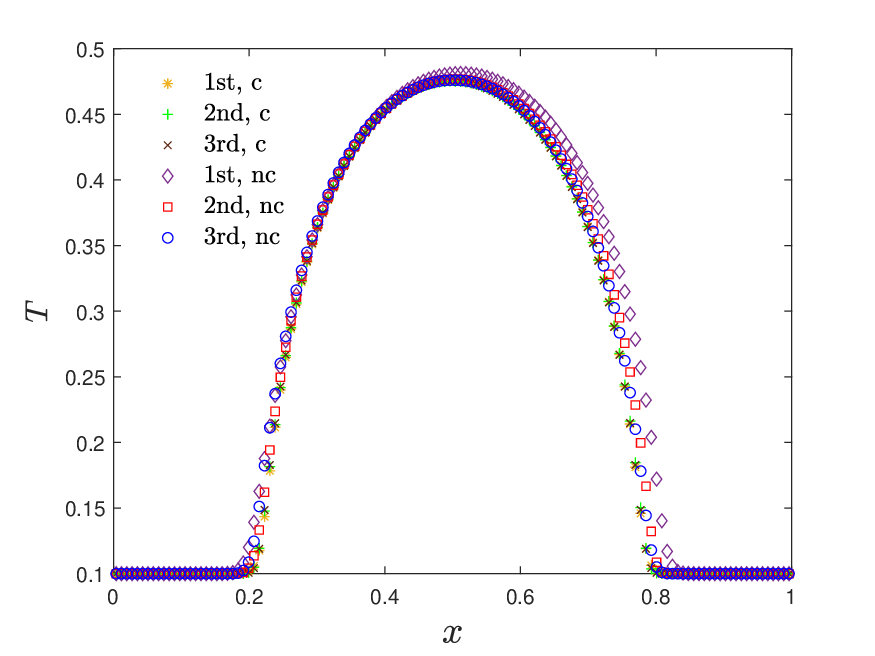}}
		\end{minipage}
		\caption{The numerical results for Example \ref{exp_conservation} with $\kappa=0.5$ and $N=128$. 1st, 2nd, and 3rd order schemes are used. ``c'': with a conservation corrector step; ``nc'': without a corrector step. Left: radiation energy $E$; Right: material temperature $T$. $\Delta t = \frac{1}{30}h$ for ``2nd, nc'', $\Delta t = \frac{1}{50}h$ for ``3rd, nc'' and $\Delta t = \frac{1}{5}h$ for all others.
		}\label{figcon_sol}
	\end{figure*}

	\begin{example}\label{exp_conservation}
		{\rm We now consider an example with a periodic boundary condition to verify the conservation errors. The initial values are taken to be } 
		\begin{equation}
			\begin{cases}\label{con_ini}
				E(x,0)=E_L+(E_R-E_L)\dfrac{1+\tanh [200(x-0.5)^2]}{2},
				\\[5pt]
				T(x,0) = E(x,0)^{1/4},
			\end{cases}
		\end{equation}
		{\rm with $E_L=1, E_R=0.0001$, which are shown in the left column of Fig. \ref{figcon_ini}. 
		We run the solution with $\kappa=0.5$ up to time $t = 1$. We take $z(x)=1$ and mesh numbers $N=128$. 1st, 2nd, and 3rd order schemes are used. ``c'' denotes numerical results obtained with the conservation corrector step, while ``nc'' is without the corrector step.
		As we have observed numerically, without the corrector step, a smaller time step is needed for the convergence of iteration. We take $\Delta t = \frac{1}{30}h$ for the 2nd method, and $\Delta t = \frac{1}{50}h$ for the 3rd order method, without a corrector step, and $\Delta t = \frac{1}{5}h$ for all others. In Fig.~\ref{figcon_ini}, on the right column, we show the time evolution of conservation errors for the total energy $Q=E+T$. We can observe that with a corrector step, the errors are much smaller than those without a corrector step. In Fig. \ref{figcon_sol}, we show the numerical solutions for different orders with or without a corrector step. As we can see, with a corrector step, all results match each other well. Without a corrector step, we can clearly observe a deviation, especially for first and second order methods. However, we would note that with a corrector step, the method itself is conservative for the total energy, although limiters for second and third order methods would slightly destroy such a conservation. }
		
	\end{example}
	
\begin{figure*}[!htbp]
		\begin{minipage}[t]{0.49\linewidth}
			\centerline{\includegraphics[scale=0.57]{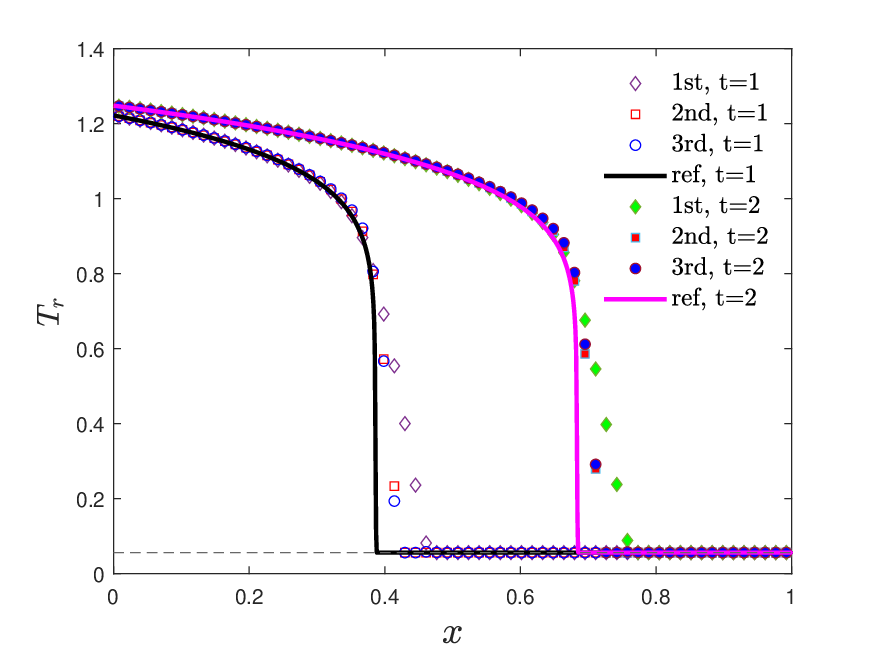}}
		\end{minipage}
		\begin{minipage}[t]{0.49\linewidth}
			\centerline{\includegraphics[scale=0.57]{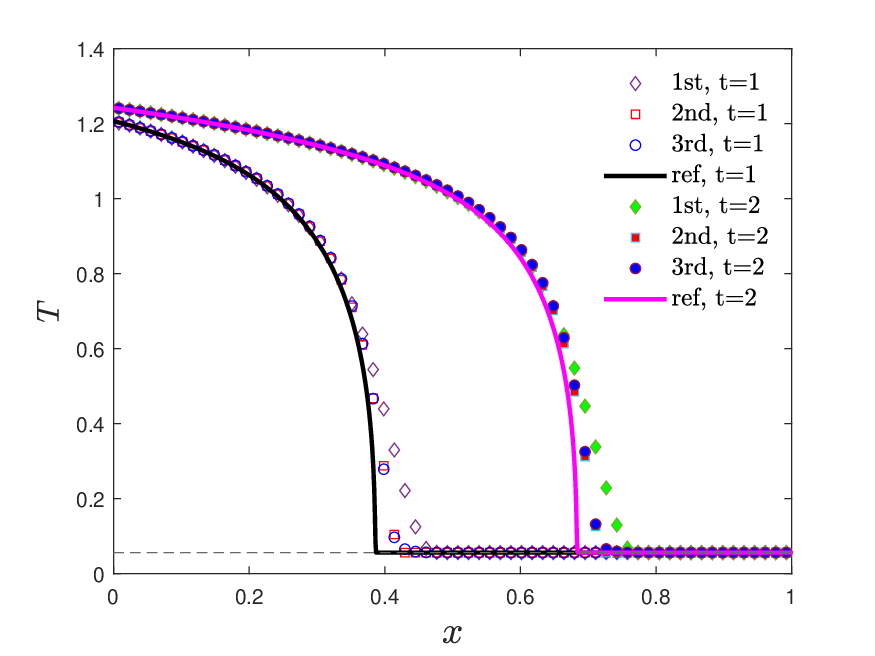}}
		\end{minipage}
		\vskip 0.1mm
		\begin{minipage}[t]{0.49\linewidth}
			\centerline{\includegraphics[scale=0.57]{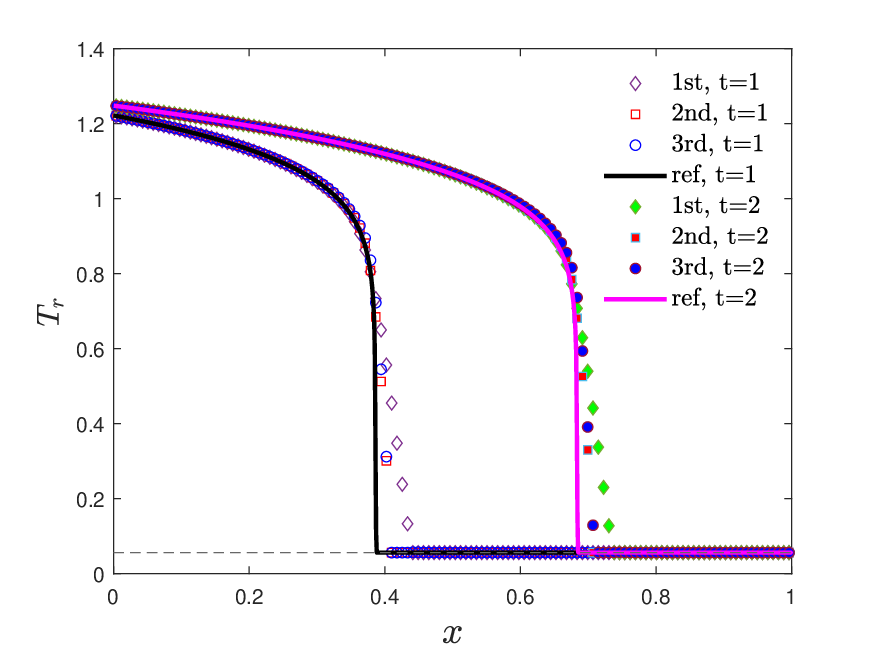}}
		\end{minipage}
		\begin{minipage}[t]{0.49\linewidth}
			\centerline{\includegraphics[scale=0.57]{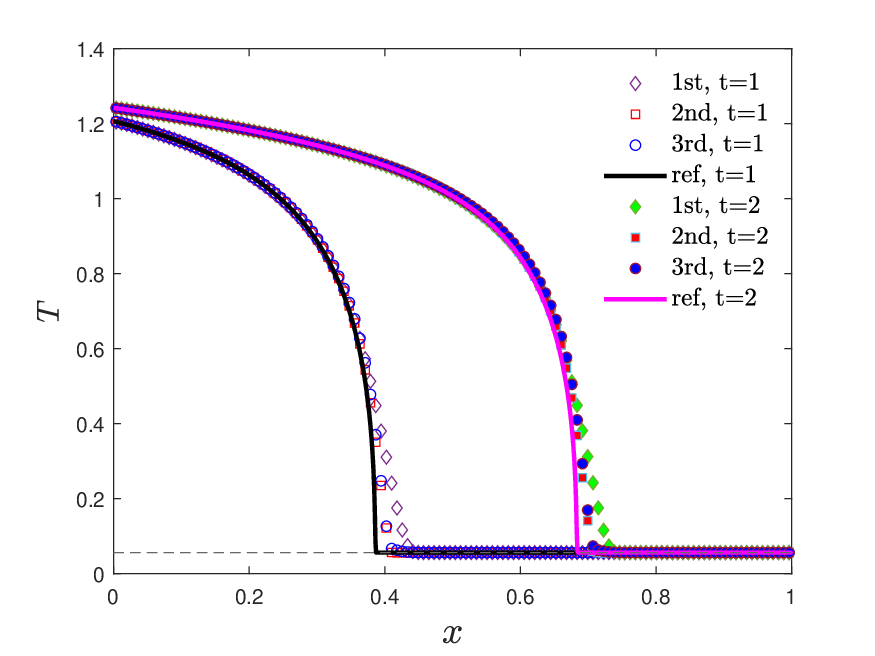}}
		\end{minipage}
		\caption{The numerical results of 1st, 2nd, and 3rd order schemes for Example \ref{Marshak_homo_problem} with $\kappa=0,$  $\Delta t = \frac{h}{5},$ time $t=1$ and $t=2.$ ``ref'' represents the reference solutions. Top: $N=64$; Bottom: $N=128$. Left column: radiation temperature $T_r$; Right column: material temperature $T$.
		}\label{Marshak_homo_1}
	\end{figure*}

\begin{figure*}[!htbp]
		\begin{minipage}[t]{0.49\linewidth}
			\centerline{\includegraphics[scale=0.57]{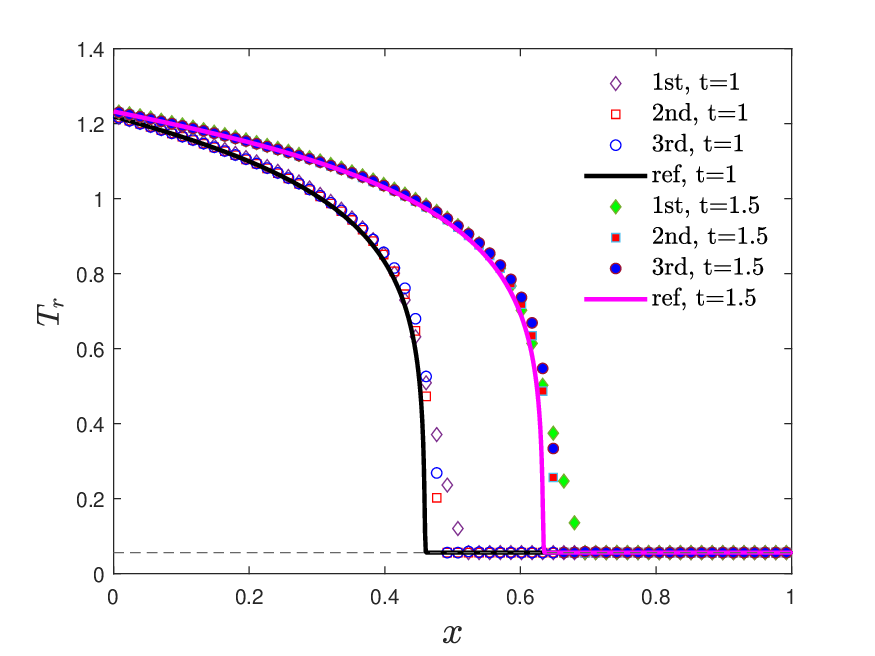}}
		\end{minipage}
		\begin{minipage}[t]{0.49\linewidth}
			\centerline{\includegraphics[scale=0.57]{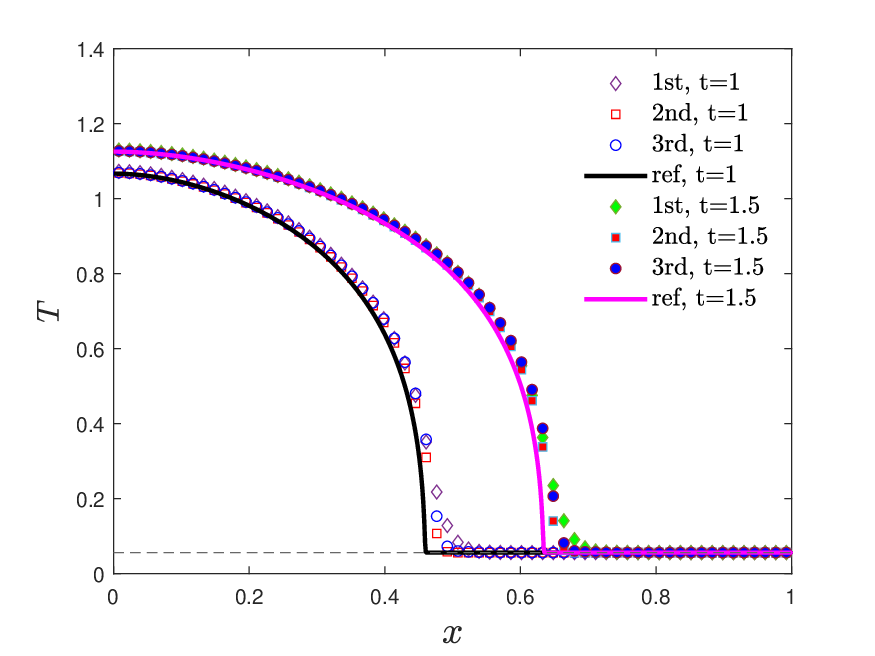}}
		\end{minipage}
		\vskip 0.1mm
		\begin{minipage}[t]{0.49\linewidth}
			\centerline{\includegraphics[scale=0.57]{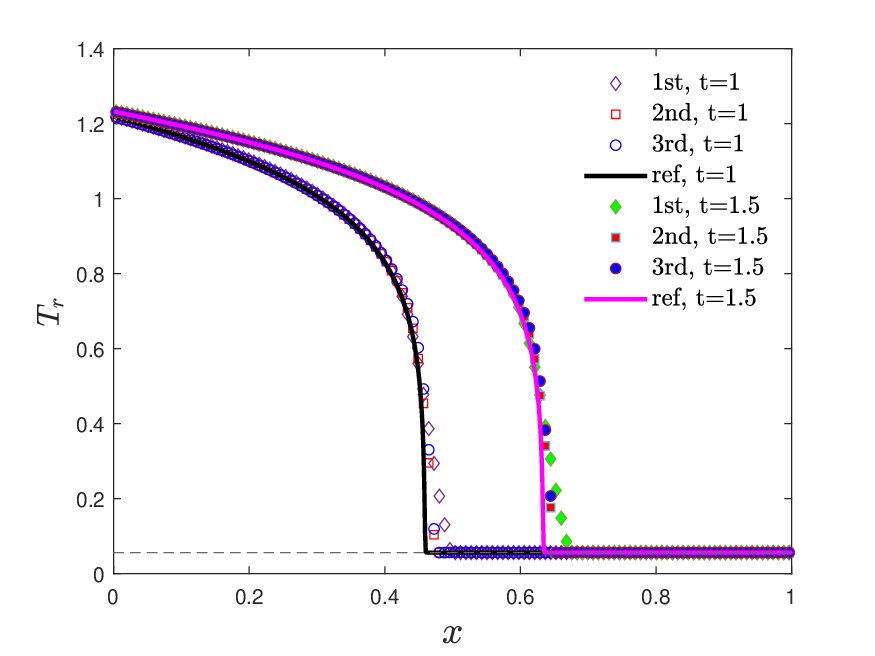}}
		\end{minipage}
		\begin{minipage}[t]{0.49\linewidth}
			\centerline{\includegraphics[scale=0.57]{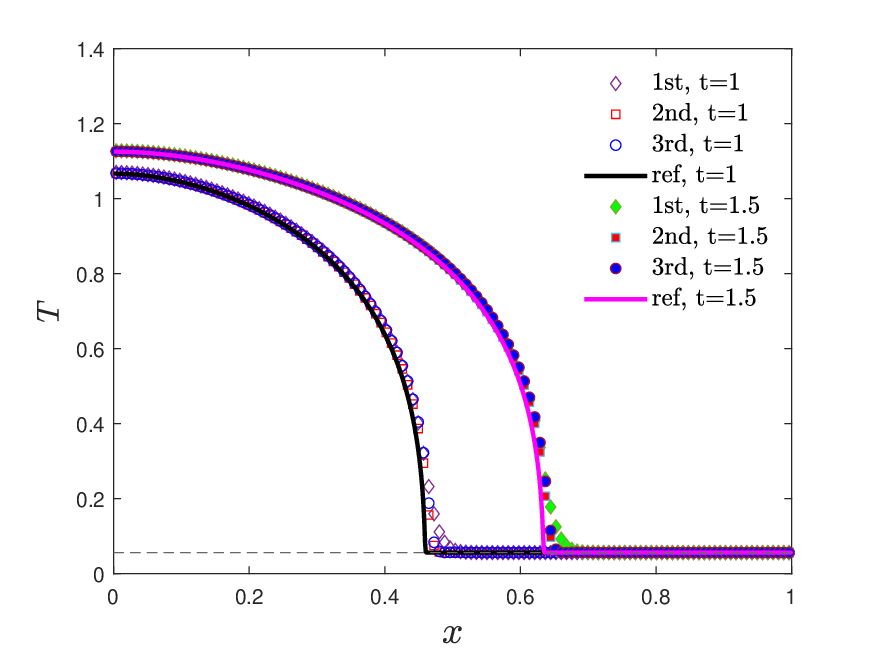}}
		\end{minipage}
		\caption{The numerical results of 1st, 2nd, and 3rd order schemes for Example \ref{Marshak_homo_problem} with $\kappa=0.1,$ $\Delta t = \frac{h}{5},$ time $t=1$ and $t=1.5.$ ``ref'' represents the reference solutions. Top: $N=64$; Bottom: $N=128$. Left column: radiation temperature $T_r$; Right column: material temperature $T$.
		}\label{Marshak_homo_2}
	\end{figure*}

\begin{figure*}[!htbp]
		\begin{minipage}[t]{0.49\linewidth}
			\centerline{\includegraphics[scale=0.57]{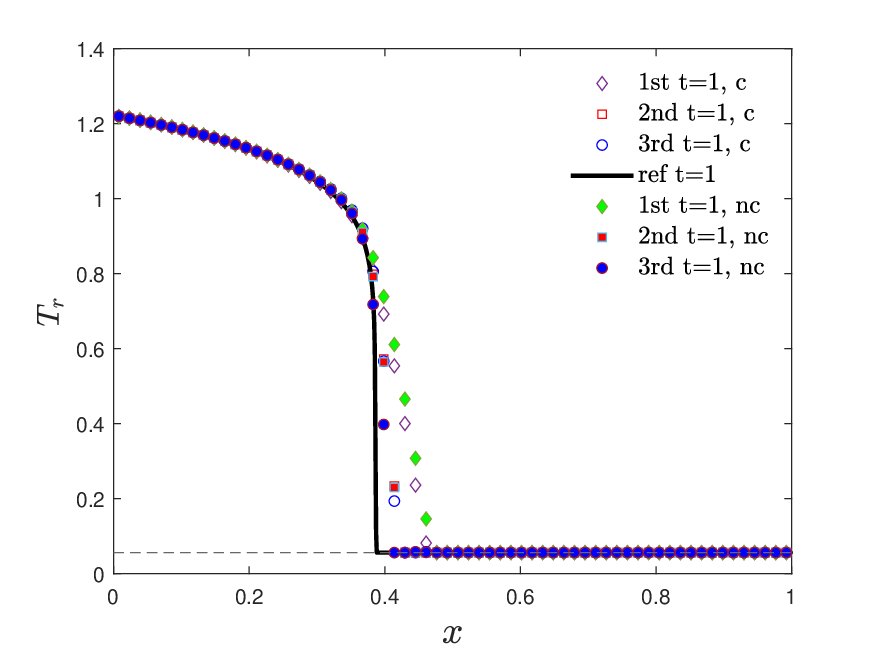}}
		\end{minipage}
		\begin{minipage}[t]{0.49\linewidth}
			\centerline{\includegraphics[scale=0.57]{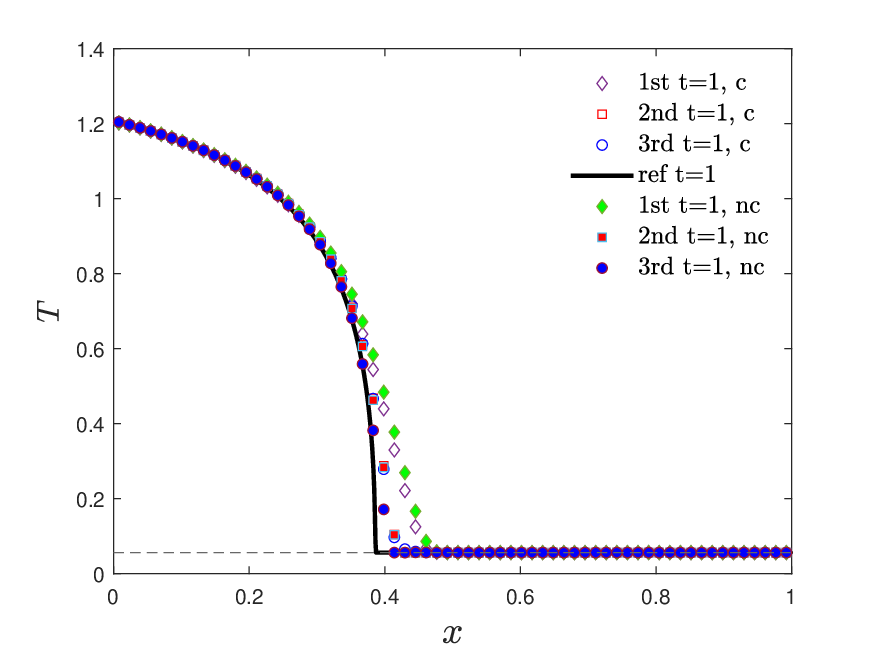}}
		\end{minipage}
		\vskip 0.1mm
		\begin{minipage}[t]{0.49\linewidth}
			\centerline{\includegraphics[scale=0.57]{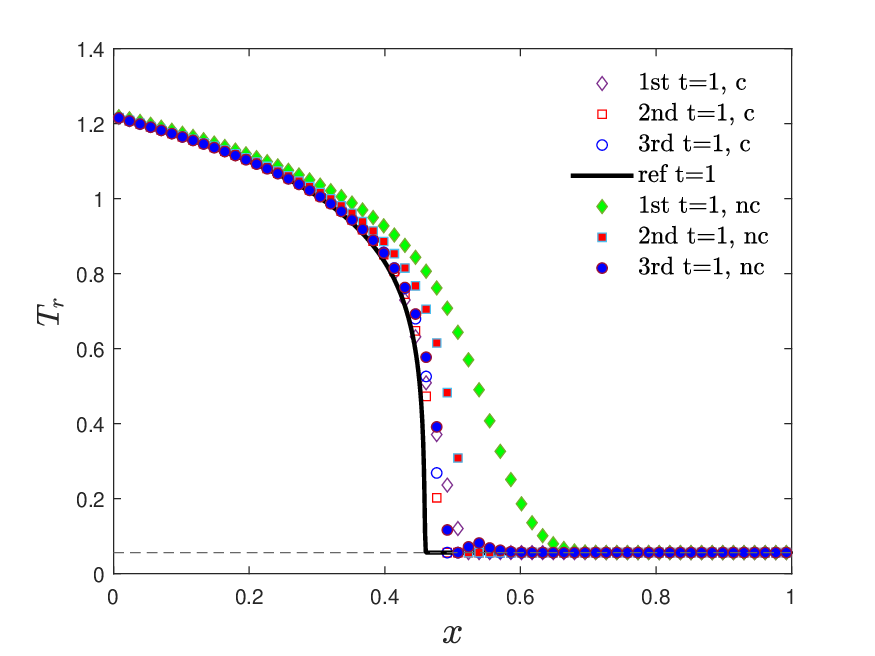}}
		\end{minipage}
		\begin{minipage}[t]{0.49\linewidth}
			\centerline{\includegraphics[scale=0.57]{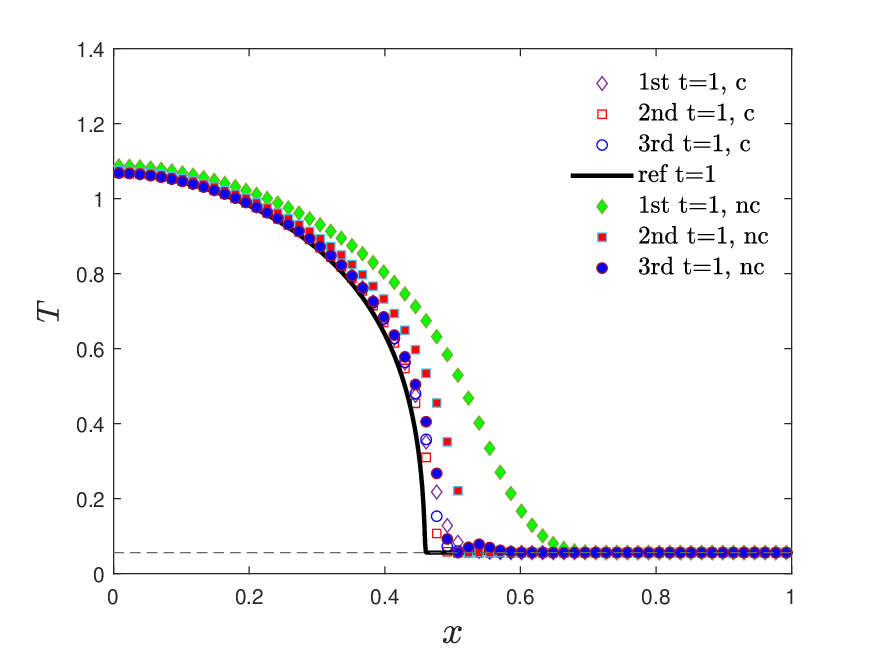}}
		\end{minipage}
		\caption{The numerical results of 1st, 2nd and 3rd order schemes for Example \ref{Marshak_homo_problem} with $t=1,$ $\Delta t = \frac{h}{5},$ $N=64$. ``ref'' means the reference solutions. Left column: radiation temperature $T_r$; Right column: material temperature $T$. Top: $\kappa=0$; Bottom: $\kappa=0.1$. ``c'' and ``nc'' denote numerical results with or without a corrector step, respectively.
		}\label{Marshak_homo_1_cnc}
	\end{figure*}

	\begin{example}\label{Marshak_homo_problem} 
	{\rm Here we consider the standard homogeneous Marshak wave problem \cite{marshak1958effect,tang2021accurate,zhang2009numerical} with the initial and boundary conditions \eqref{Marshak_1D_IBC}, where the atomic mass number $z(x)=1$. This is a benchmark problem for radiation diffusion problems. It is very challenging since the left inflow boundary value is not consistent with the initial datum. Similarly, we consider both $\kappa=0$ and $\kappa=0.1$ in $D_t$, respectively. The reference solutions are computed by a 1st order explicit RK LDG method on $1024$ elements with a small enough time step. {In Fig.~\ref{Marshak_homo_1} and Fig.~\ref{Marshak_homo_2}, we show the results at different times $t=1,2$ and $t=1,1.5$, respectively.} The left column is for the radiation temperature $T_r$ and the right column is for the material temperature $T$. A mesh refinement from $N=64$ to $N=128$ is also considered. As we can see that the 1st order numerical results are very deviated from the reference solutions, second and third order methods capture the sharp fronts more accurately. With mesh refinement, all results are getting closer to the reference solutions, indicating that our methods are convergent with mesh refinement. The results are also consistent with those in \cite{mousseau2000physics,knoll2007numerical}. As compared to \cite{rongpei2012discontinuous,huang2016monotone}, a relatively larger time step size $\Delta t=\frac15h$ can be used, {where $\Delta t=\frac{4}{125}h$ with $h=\frac{1}{80}$ in \cite{rongpei2012discontinuous} and $\Delta t=\frac18h^2$ with $h=\frac{1}{16}$ in \cite{huang2016monotone}}. In Fig~\ref{Marshak_homo_1_cnc}, we compare the results with and without a corrector step. In the case of $\kappa=0$, the system $E$ and $B$ in \eqref{RDE_2D} itself is in a conservative form, although it only conserves the total energy $Q=E+T$ up to an error of numerical precision, the results with or without a corrector step are almost the same. For $\kappa=0.1$, the results with a corrector step clearly match the reference solutions better, especially for first and second order methods. This has demonstrated that the corrector step is very necessary.

	}
    \end{example}

    \begin{rmk}
	As discussed in Remark \ref{rem_1}, if we eliminate the predictor step and directly employ a Picard iteration in the corrector step, it may not work well for some challenge problems. For this standard Marshak wave problem, if we consider $E(x,0)=10^{-7}$, $N=64$ and $\Delta t=\frac{h}{5}$, even a first order scheme does not converge well. Such an approach with higher orders work even worse. However, our predictor-corrector procedure works well for these test cases. With the results in the previous Example \ref{exp_conservation}, we have shown that our methods can ensure both conservation and robustness.
    \end{rmk}
	
\begin{figure*}[!htbp]
		\begin{minipage}[t]{0.49\linewidth}
			\centerline{\includegraphics[scale=0.55]{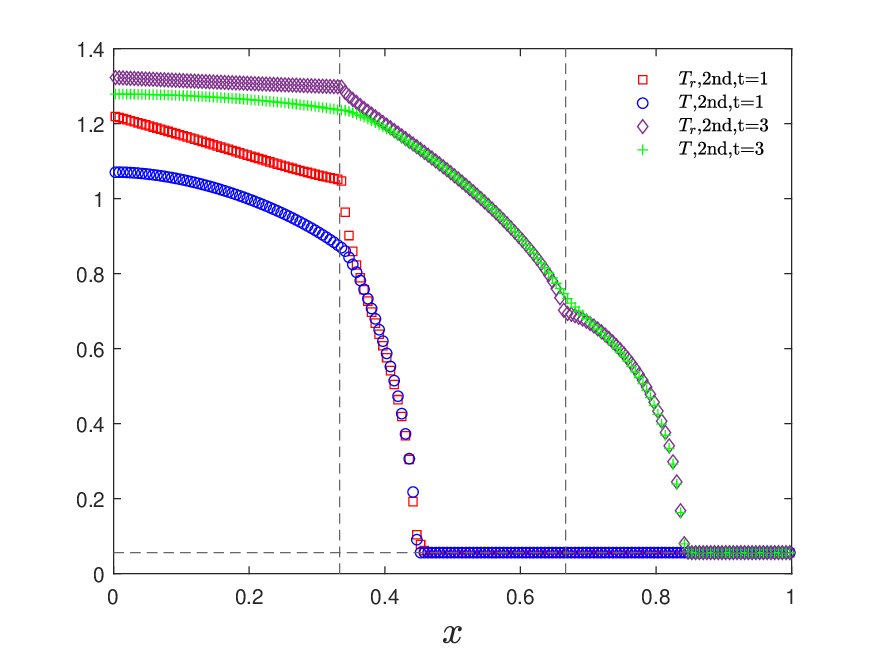}}
		\end{minipage}
		\begin{minipage}[t]{0.49\linewidth}
			\centerline{\includegraphics[scale=0.55]{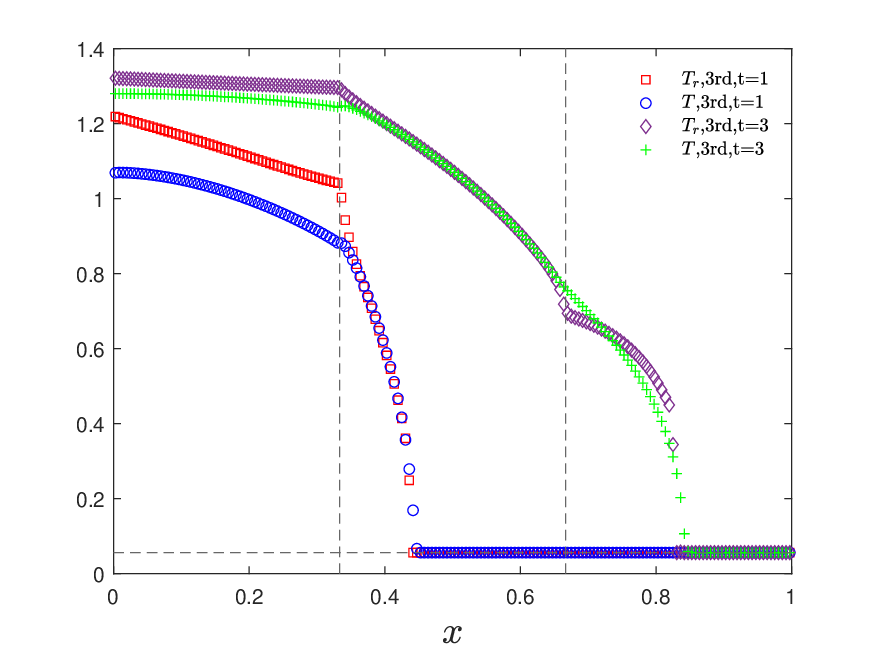}}
		\end{minipage}
		\caption{The numerical results of 2nd and 3rd order schemes for Example \ref{Marshak_hete_problem} at time $t=1$ and $t=3$. $N=180$ and $\Delta t=\frac{h}{15}$. Left: 2nd order; Right: 3rd order.
		}\label{Marshak_hete_1}
	\end{figure*}

\begin{figure*}[!htbp]
		\begin{minipage}[t]{0.49\linewidth}
			\centerline{\includegraphics[scale=0.55]{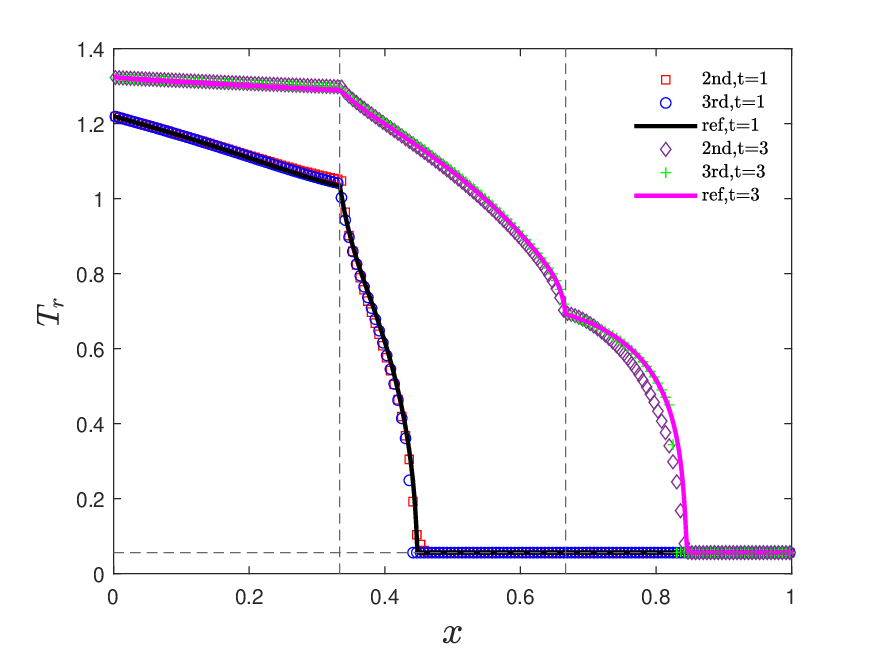}}
		\end{minipage}
		\begin{minipage}[t]{0.49\linewidth}
			\centerline{\includegraphics[scale=0.55]{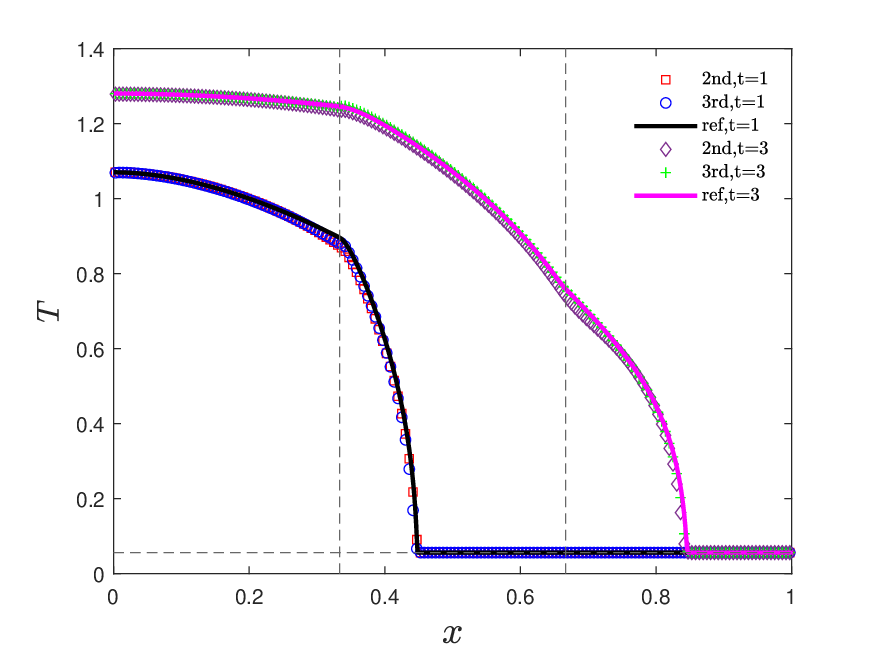}}
		\end{minipage}
		\caption{The comparison between the results of 2nd and 3rd order schemes for Example \ref{Marshak_hete_problem} at time $t=1$ and $t=3$. Left: radiation temperature $T_r$; Right: material temperature $T$. $N=180$ and $\Delta t=\frac{h}{15}$.
		}\label{Marshak_hete_2}
	\end{figure*}

	\begin{example}\label{Marshak_hete_problem} 
		{\rm In this example, we consider a heterogeneous Marshak problem \cite{marshak1958effect,tang2021accurate,zhang2009numerical} with the initial and boundary conditions \eqref{Marshak_1D_IBC} with $z = 3$ inside the interval $[\frac13, \frac23]$ and $z = 1$ elsewhere. 
			We take $\kappa=0.1$ and two different times $t=1,3$. In Fig.~\ref{Marshak_hete_1}, we show the results of 2nd and 3rd order methods, on a mesh $N=180$ and the time step is taken to be $\Delta t=\frac{1}{15}h$. The 1st order results are not presented since they deviate a lot as shown above. From the numerical results, we can observe that the radiation temperature $T_r$ and the material temperature propagate much slower with a large value of $z$ in the middle region. Besides, the radiation temperature $T_r$ is approaching the material temperature $T$ in the middle region, reaching a thermodynamic equilibrium. In Fig.~\ref{Marshak_hete_2}, we compare the numerical results between the 2nd and 3rd order methods. ``ref'' represents the reference solution obtained by the 1st order explicit RK LDG method on $1024$ elements with a small enough time step. We can find that the solutions of the 3rd order method match the reference better than the 2nd order ones, especially for the radiation temperature. From this example, we can see that high order methods can capture fine structures than corresponding lower order methods.   }
	\end{example}
	
	\begin{example}\label{test2D_order} 
		{(\bf Accuracy test in 2D)} {\rm In this example we test the errors and convergence orders of accuracy at the equilibrium in the 2D case. We consider a 2D radiation diffusion problem with two given source terms $f_1(x,y,t)$ and $f_2(x,y,t)$, which are chosen properly so that exact solutions are available for the following system: }
		\begin{equation}
			\begin{cases}
				\dfrac{\partial E}{\partial t}-\nabla \cdot\left(D_{r} \nabla E\right)=\sigma\left(T^{4}-E\right)+f_1(x,y,t),
				\\[5pt]
				\dfrac{\partial T}{\partial t}-\nabla \cdot\left(D_{t} \nabla T\right)=\sigma\left(E-T^{4}\right)+f_2(x,y,t),
			\end{cases}
		\end{equation}
	\end{example}
	\hspace{-6mm} where $T(x,y,t) = \left(0.8+0.1\sin\left(2\pi\left(x-t\right)\right)\right)\left(0.8+0.1\sin\left(2\pi\left(y-t\right)\right)\right)$ and $E(x,y,t) = T(x,y,t)^4$, $\Omega = [0,1]^2$. For convenience we remove the artificially added term in the energy radiation coefficient, that is $D_r=\frac{1}{3\sigma}$.
	
	We take $\kappa=0.01$ and $z(x)=1$. The problem is run to time $t = 0.5$ using the 1st, 2nd, and 3rd order methods, respectively. We list the numerical $L^2$ errors and orders of accuracy with different time step sizes in Table \ref{table_4}. Similar results as the 1D case are obtained.
	\begin{table}[!htbp]
		\centering
		\caption{The numerical $L^2$ errors and orders of accuracy for the 1st, 2nd and 3rd order schemes with different time steps for Example \ref{test2D_order}. $t=0.5.$}\label{table_4}
		\begin{tabular}{ccccccccc}
			\toprule
			& \multirow{2}{*}{N} & \multirow{2}{*}{$\Delta t$}           & \multicolumn{1}{c}{$L^2$ error} & order & \multicolumn{1}{c}{$L^2$ error} & order & \multicolumn{1}{c}{$L^2$ error} & order \\ \cline{4-9} 
			&                 &                                       & \multicolumn{2}{c}{E}                   & \multicolumn{2}{c}{T}                   & \multicolumn{2}{c}{B}                   \\ \hline
			\multirow{9}{*}{1st} & $16\times 16$ & \multirow{3}{*}{$\frac{h}{2}$} & 1.22e-2     & -                        & 1.02e-2     & -                        & 1.23e-2     & -                        \\
			& $32\times 32$ &                                       & 6.21e-3     & 0.98                     & 5.17e-3     & 0.98                     & 6.23e-3     & 0.98                     \\
			& $64\times 64$ &                                       & 3.13e-3     & 0.99                     & 2.60e-3     & 0.99                     & 3.14e-3     & 0.99                     \\
			& $16\times 16$ & \multirow{3}{*}{$h$}           & 1.85e-2     & -                        & 1.80e-2     & -                        & 2.10e-2     & -                        \\
			& $32\times 32$ &                                       & 9.54e-3     & 0.96                     & 9.35e-3     & 0.95                     & 1.08e-2     & 0.95                     \\
			& $64\times 64$ &                                       & 4.84e-3     & 0.98                     & 4.77e-3     & 0.97                     & 5.50e-3     & 0.98                     \\
			& $16\times 16$ & \multirow{3}{*}{$2h$}          & 3.67e-2     & -                        & 4.07e-2     & -                        & 4.48e-2     & -                        \\
			& $32\times 32$ &                                       & 1.89e-2     & 0.96                     & 2.14e-2     & 0.93                     & 2.33e-2     & 0.94                     \\
			& $64\times 64$ &                                       & 9.65e-3     & 0.97                     & 1.10e-2     & 0.96                     & 1.19e-2     & 0.97                     \\
			\hline
			\multirow{9}{*}{2nd} & $16\times 16$ & \multirow{3}{*}{$\frac{h}{2}$} & 1.11e-3     & -     & 1.04e-3     & -     & 1.14e-3     & -     \\
			& $32\times 32$ &                                       & 2.78e-4     & 2.00  & 2.63e-4     & 1.99  & 2.90e-4     & 1.98  \\
			& $64\times 64$ &                                       & 6.93e-5     & 2.00  & 6.57e-5     & 2.00  & 7.27e-5     & 2.00  \\
			& $16\times 16$ & \multirow{3}{*}{$h$}           & 2.56e-3     & -     & 2.89e-3     & -     & 2.51e-3     & -     \\
			& $32\times 32$ &                                       & 6.33e-4     & 2.01  & 7.45e-4     & 1.95  & 6.40e-4     & 1.97  \\
			& $64\times 64$ &                                       & 1.57e-4     & 2.01  & 1.88e-4     & 1.99  & 1.61e-4     & 1.99  \\
			& $16\times 16$ & \multirow{3}{*}{$2h$}          & 9.68e-3     & -     & 1.01e-2     & -     & 8.91e-3     & -     \\
			& $32\times 32$ &                                       & 2.34e-3     & 2.05  & 2.79e-3     & 1.86  & 2.33e-3     & 1.94  \\
			& $64\times 64$ &                                       & 5.81e-4     & 2.01  & 7.22e-4     & 1.95  & 5.93e-4     & 1.97 \\ \hline
			\multirow{9}{*}{3rd} & $16\times 16$ & \multirow{3}{*}{$\frac{h}{2}$} & 7.63e-5     & -     & 2.40e-5     & -     & 2.59e-5     & -     \\
			& $32\times 32$ &                                       & 1.13e-5     & 2.75  & 3.26e-6     & 2.88  & 3.31e-6     & 2.97  \\
			& $64\times 64$ &                                       & 1.59e-6     & 2.83  & 4.32e-7     & 2.91  & 4.22e-7     & 2.97  \\
			& $16\times 16$ & \multirow{3}{*}{$h$}           & 4.56e-4     & -     & 1.11e-4     & -     & 1.23e-4     & -     \\
			& $32\times 32$ &                                       & 7.17e-5     & 2.67  & 1.62e-5     & 2.78  & 1.50e-5     & 3.04  \\
			& $64\times 64$ &                                       & 1.08e-5     & 2.73  & 2.37e-6     & 2.78  & 1.95e-6     & 2.94  \\
			& $16\times 16$ & \multirow{3}{*}{$2h$}          & 2.80e-3     & -     & 8.52e-4     & -     & 1.11e-3     & -     \\
			& $32\times 32$ &                                       & 4.55e-4     & 2.62  & 1.10e-4     & 2.95  & 1.21e-4     & 3.20  \\
			& $64\times 64$ &                                       & 7.16e-5     & 2.67  & 1.61e-5     & 2.77  & 1.47e-5     & 3.04 
			\\ \bottomrule\end{tabular}
	\end{table}
	
\begin{figure*}[htbp]
		\begin{minipage}[t]{0.32\linewidth}
			\centerline{\includegraphics[scale=0.38]{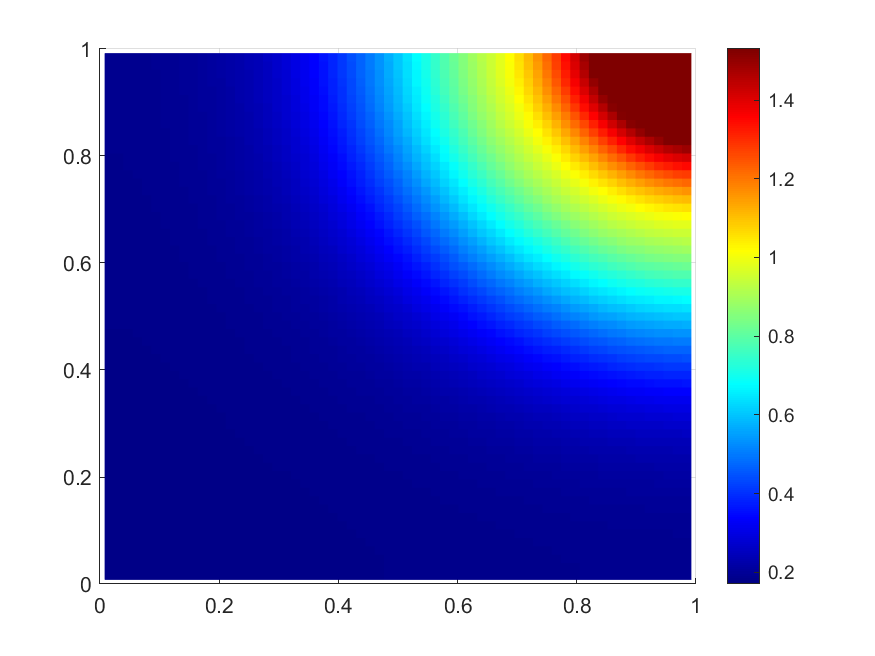}}
		\end{minipage}
		\begin{minipage}[t]{0.32\linewidth}
			\centerline{\includegraphics[scale=0.38]{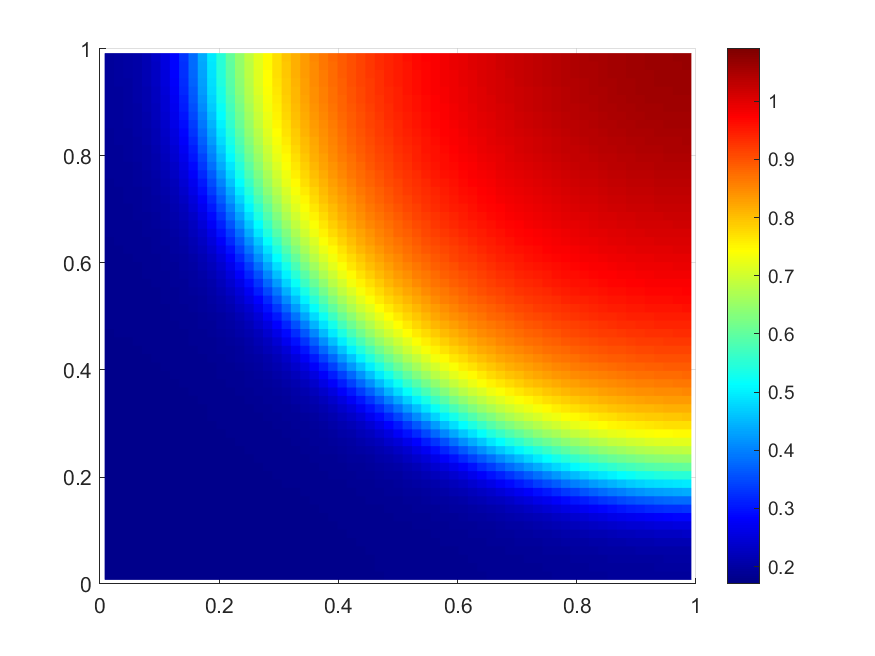}}
		\end{minipage}
		\begin{minipage}[t]{0.32\linewidth}
			\centerline{\includegraphics[scale=0.38]{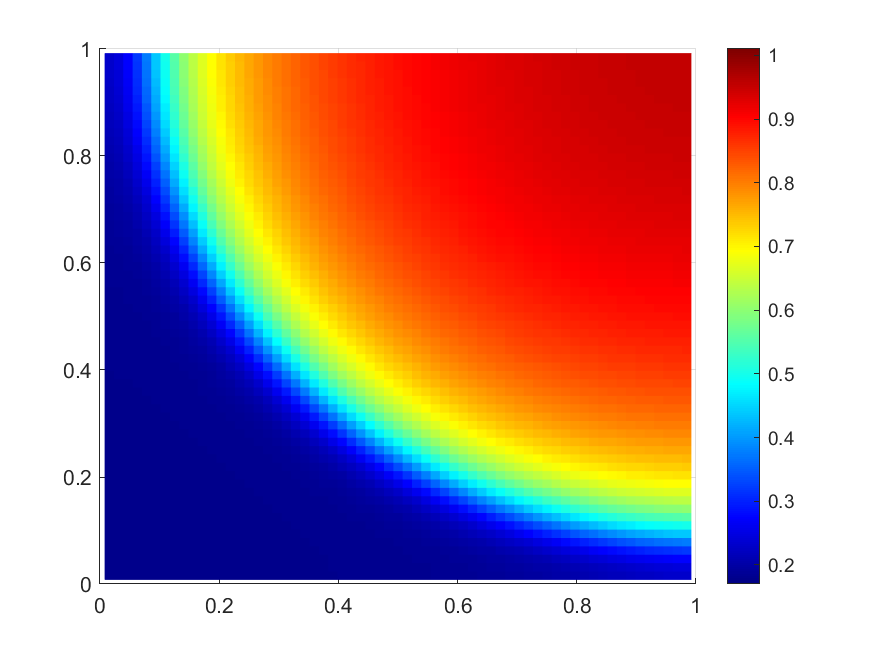}}
		\end{minipage}
		\vfill
		\begin{minipage}[t]{0.32\linewidth}
			\centerline{\includegraphics[scale=0.38]{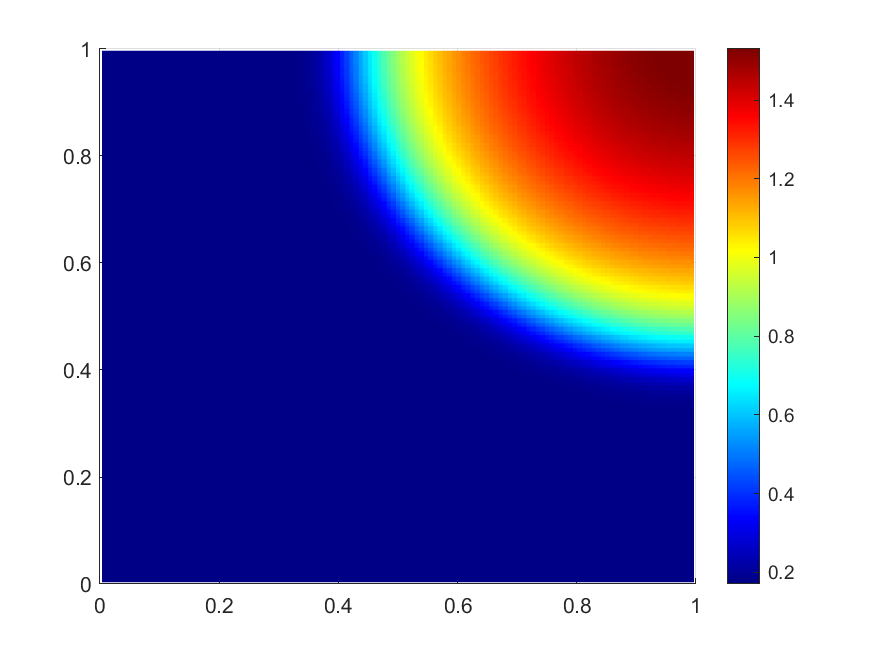}}
		\end{minipage}
		\begin{minipage}[t]{0.32\linewidth}
			\centerline{\includegraphics[scale=0.38]{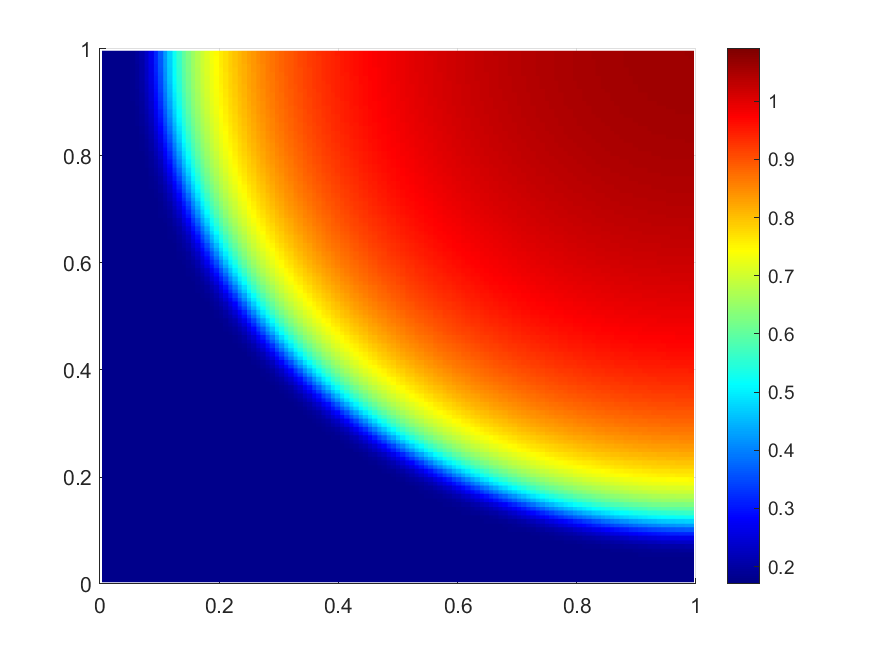}}
		\end{minipage}
		\begin{minipage}[t]{0.32\linewidth}
			\centerline{\includegraphics[scale=0.38]{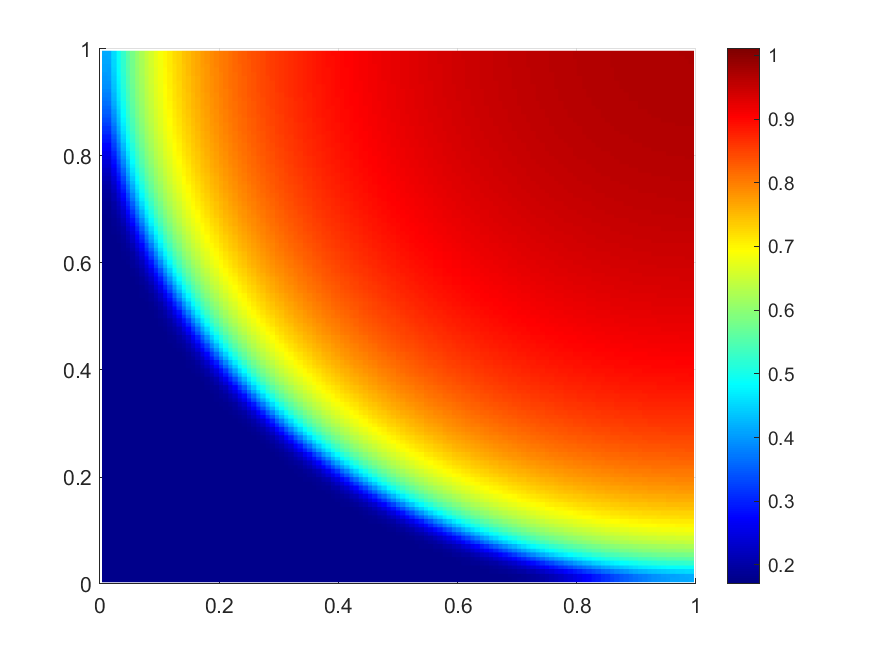}}
		\end{minipage}
		\vfill
		\begin{minipage}[t]{0.32\linewidth}
			\centerline{\includegraphics[scale=0.38]{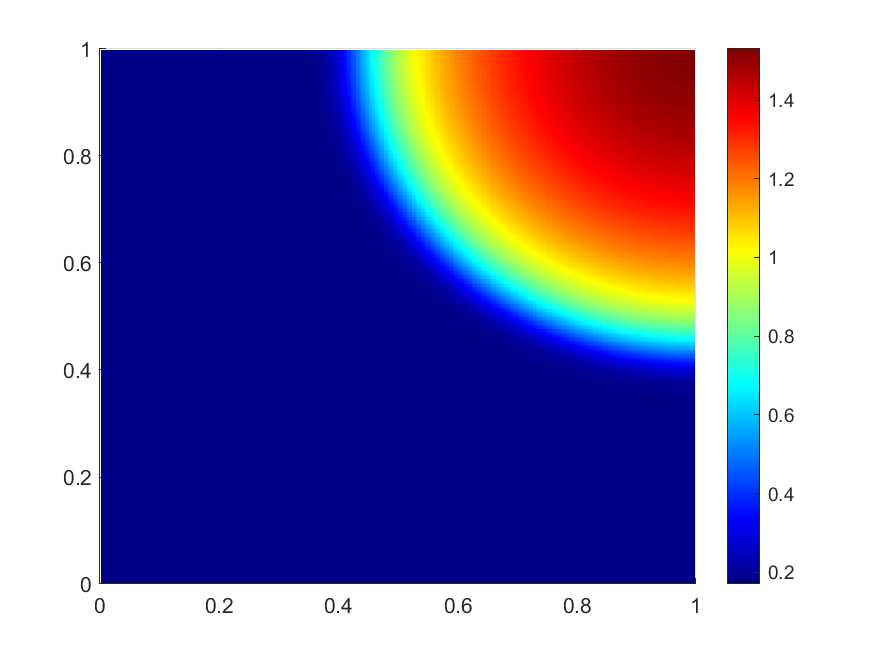}}
		\end{minipage}
		\begin{minipage}[t]{0.32\linewidth}
			\centerline{\includegraphics[scale=0.38]{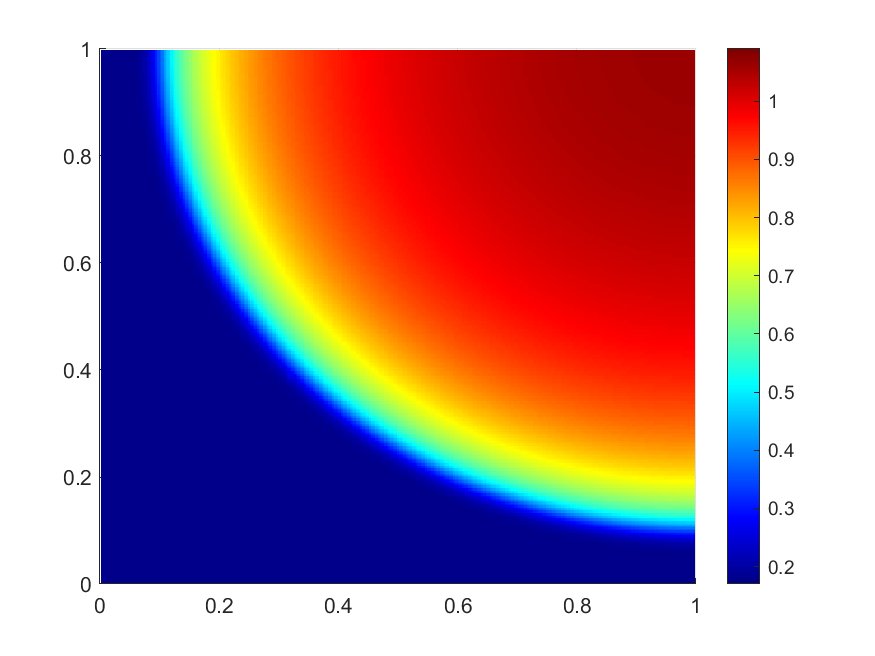}}
		\end{minipage}
		\begin{minipage}[t]{0.32\linewidth}
			\centerline{\includegraphics[scale=0.38]{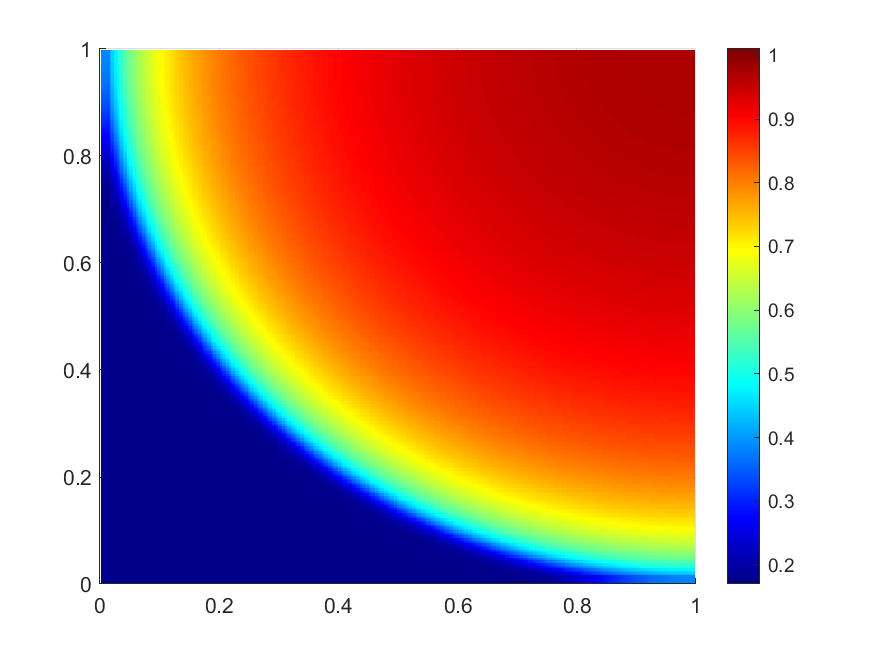}}
		\end{minipage}
		\caption{The numerical results for Example \ref{blast_homo_pb}. From top to bottom: the radiation temperature $T_r$ for 1st, 2nd, and 3rd order schemes. From left to right: time $t=0.5,1.5,2.$ $\Delta t=\frac15 h$.
		}\label{blast_homo_1}
	\end{figure*}

\begin{figure*}[htbp]
		\begin{minipage}[t]{0.32\linewidth}
			\centerline{\includegraphics[scale=0.38]{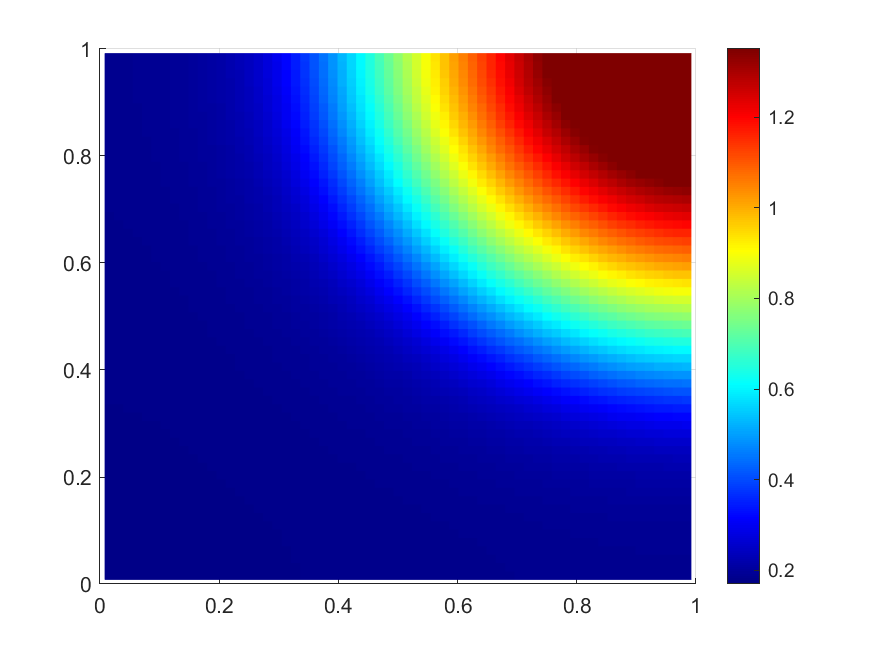}}
		\end{minipage}
		\begin{minipage}[t]{0.32\linewidth}
			\centerline{\includegraphics[scale=0.38]{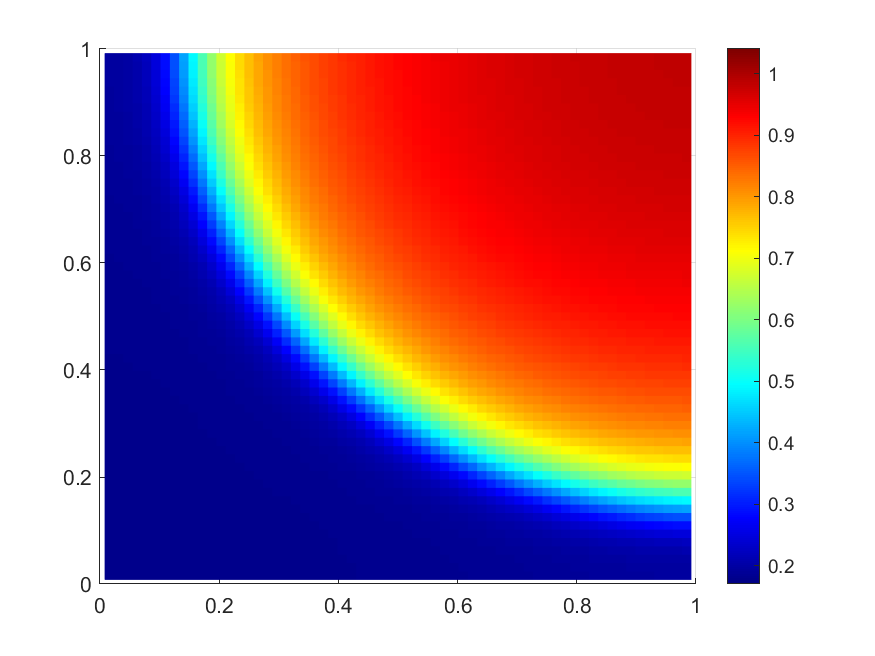}}
		\end{minipage}
		\begin{minipage}[t]{0.32\linewidth}
			\centerline{\includegraphics[scale=0.38]{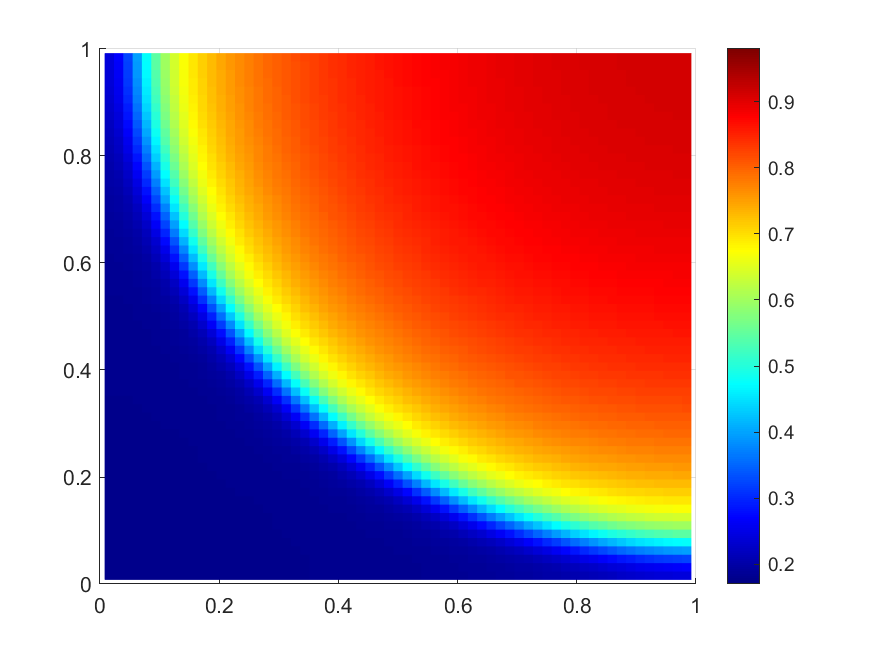}}
		\end{minipage}
		\vfill
		\begin{minipage}[t]{0.32\linewidth}
			\centerline{\includegraphics[scale=0.38]{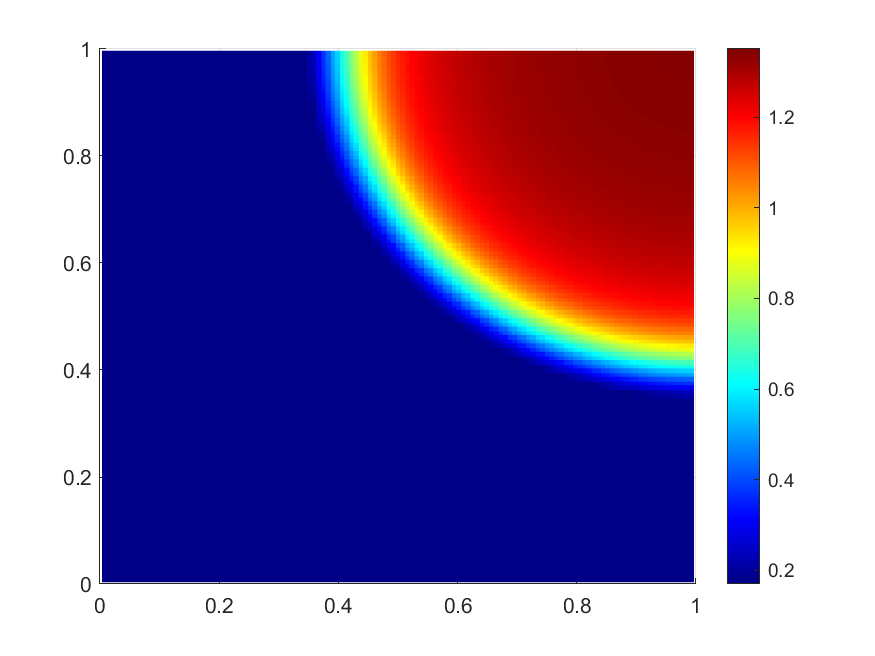}}
		\end{minipage}
		\begin{minipage}[t]{0.32\linewidth}
			\centerline{\includegraphics[scale=0.38]{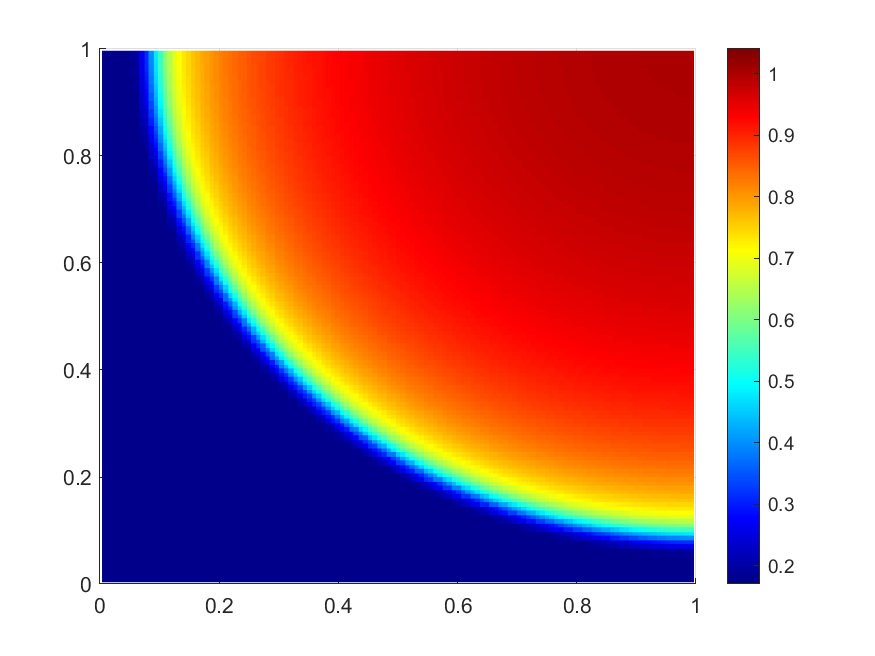}}
		\end{minipage}
		\begin{minipage}[t]{0.32\linewidth}
			\centerline{\includegraphics[scale=0.38]{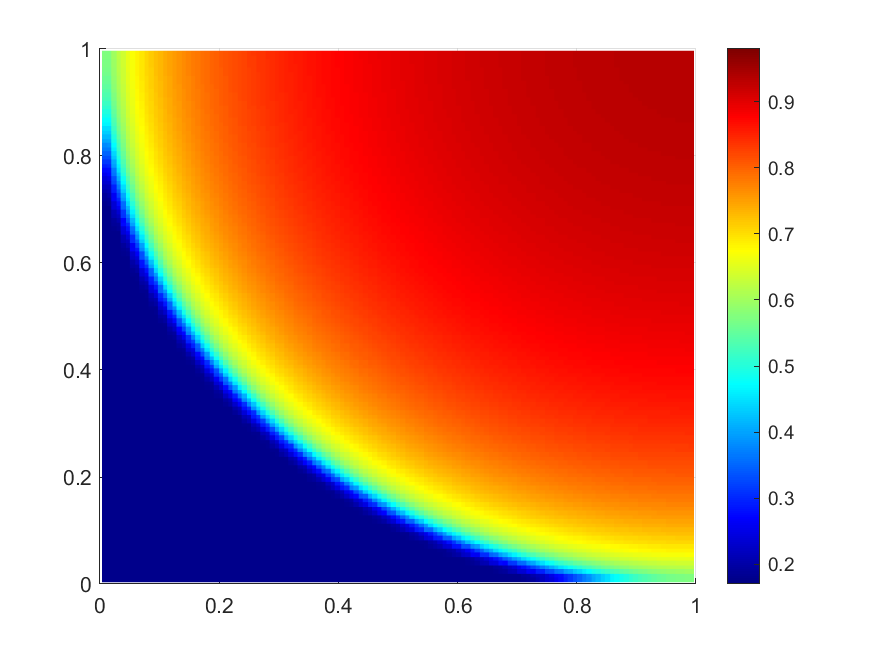}}
		\end{minipage}
		\vfill
		\begin{minipage}[t]{0.32\linewidth}
			\centerline{\includegraphics[scale=0.38]{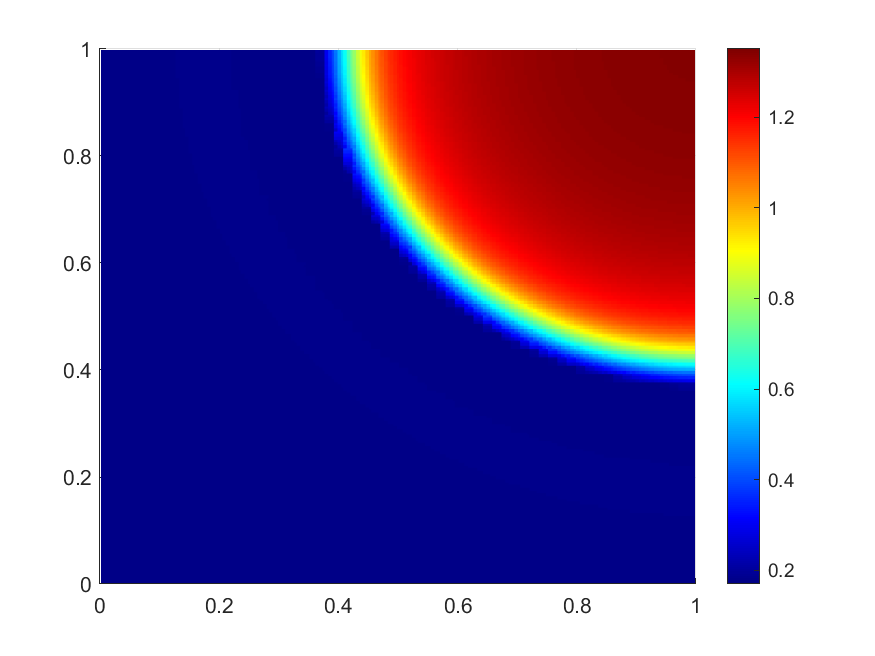}}
		\end{minipage}
		\begin{minipage}[t]{0.32\linewidth}
			\centerline{\includegraphics[scale=0.38]{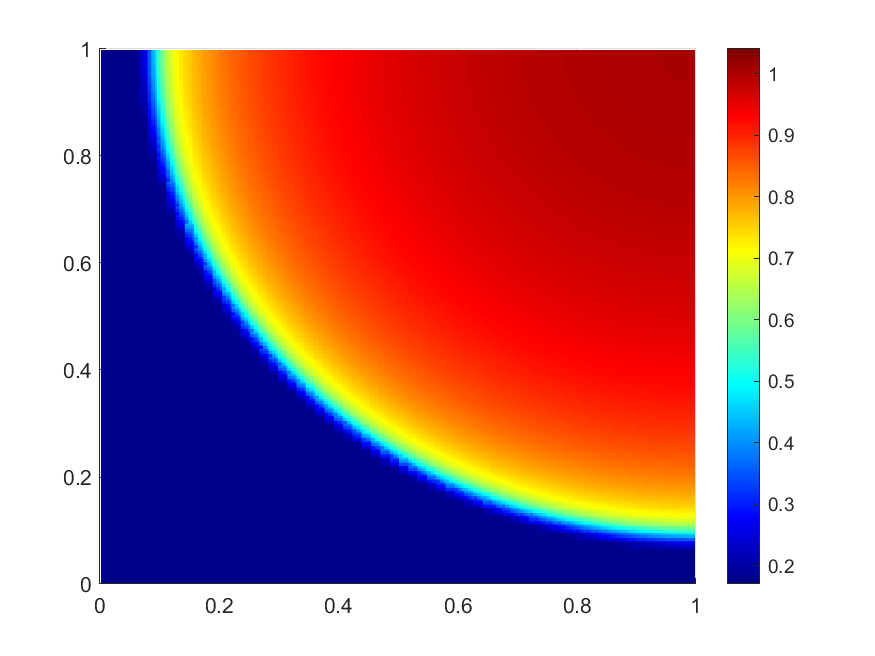}}
		\end{minipage}
		\begin{minipage}[t]{0.32\linewidth}
			\centerline{\includegraphics[scale=0.38]{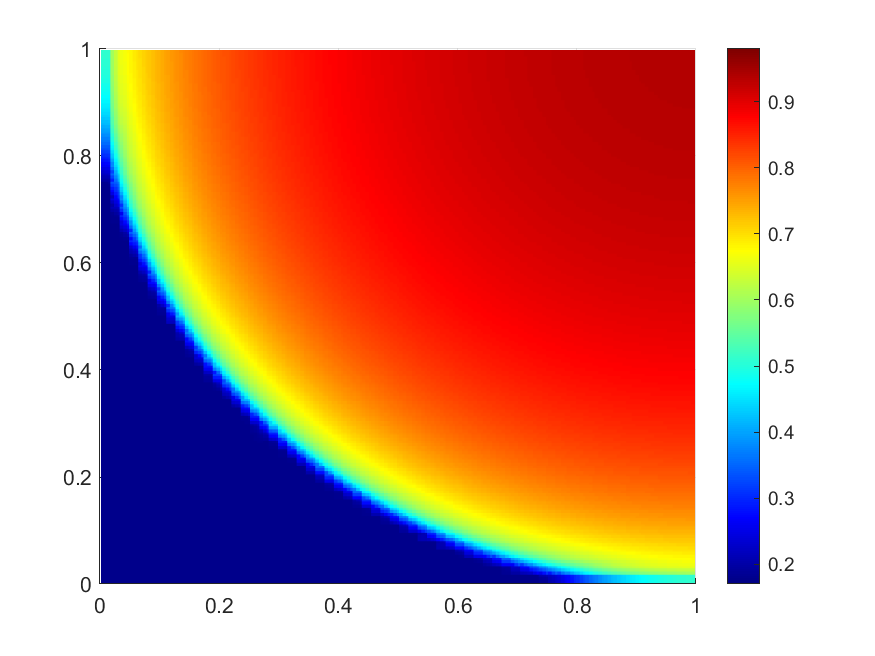}}
		\end{minipage}
		\caption{The numerical results for Example \ref{blast_homo_pb}. From top to bottom: the material temperature $T$ for 1st, 2nd, and 3rd order schemes. From left to right: time $t=0.5,1.5,2.$ $\Delta t=\frac15 h$.
		}\label{blast_homo_2}
	\end{figure*}

\begin{figure*}[!htbp]
		\begin{minipage}[t]{0.49\linewidth}
			\centerline{\includegraphics[scale=0.55]{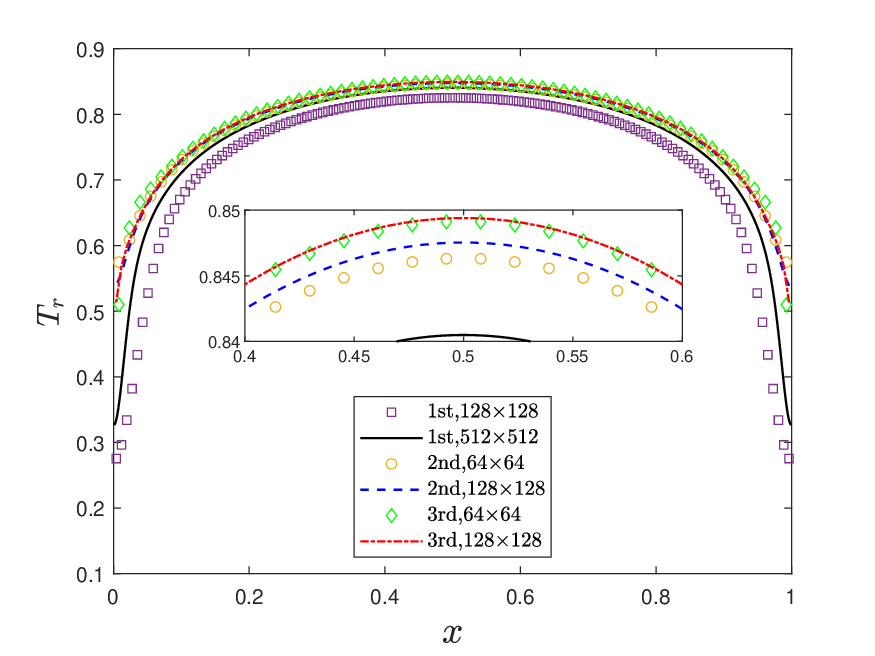}}
		\end{minipage}
		\begin{minipage}[t]{0.49\linewidth}
			\centerline{\includegraphics[scale=0.55]{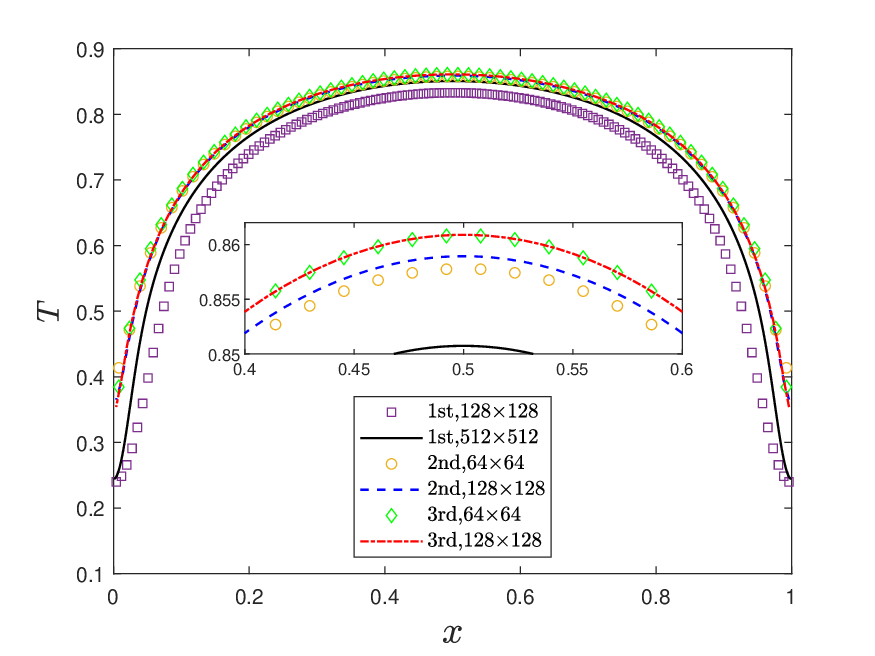}}
		\end{minipage}
		\caption{The numerical results for Example \ref{blast_homo_pb} at time $t=2$. Left: the radiation temperature $T_r$; Right: the material temperature $T$. $\Delta t =\frac{1}{5} h$.
		}\label{blast_homo_comparison}
	\end{figure*}

	\begin{example}\label{blast_homo_pb}
		{\rm Now we consider a 2D blast wave problem \eqref{RDE} with initial and boundary conditions given by}
		\begin{equation}
			\label{RDE_2D_IBC}
			\begin{cases}
				E(x,y,0)=10^{-3}+100e^{-100((x-1)^2+(y-1)^2)},\\
				T(x,y,0)=E(x,y,0)^{\frac{1}{4}},\\
				\frac{\partial E}{\partial \bm{n}}|_{\partial \Omega}=0,\quad \frac{\partial T}{\partial \bm{n}}|_{\partial \Omega}=0,
			\end{cases}
		\end{equation}
		{\rm where $\Omega=[0,1]\times[0,1]$, $\kappa=0.01$.}
		
		{\rm First we consider a homogeneous case with $z(x,y)=1$. In Fig.~\ref{blast_homo_1} and Fig.~\ref{blast_homo_2}, we show the numerical solutions $T_r$ and $T$ at time $t=0.5,1.5,2$ on a mesh of $64\times 64$, respectively. We take a time step $\Delta t=\frac{1}{5}h$. From top to bottom, numerical results obtained by the 1st, 2nd, and 3rd order methods are presented, respectively. We can observe that our methods can capture the sharp fronts of the blast wave well. Moreover, in Fig.~\ref{blast_homo_comparison}, we show the cutting plots for the numerical solutions along $y=x+1$ at time $t=2$ on a mesh of $N\times N$ for $N=64, 128, 512$. We can see that the three methods converge and match each other. 2nd and 3rd order methods are clearly better than the 1st order method. From zoom-in figures, we can see the results of 3rd order method are slightly better than the 2nd order ones.}              
	\end{example}
	
\begin{figure*}[htbp]
		\begin{center}
			\begin{tikzpicture}
				\draw [thick, <->] (0,5) node [left] {$y$}
				-- (0,0) -- (5,0) node [below right] {$x$};
				\draw (0,0) rectangle (4,4);
				\draw [thick] (0,-.1) node[below]{\rm 0} -- (0,0);
				\draw [thick] (4,-.1) node[below]{\rm 1} -- (4,0);
				\draw [thick] (-.1,3/4) node[left]{$\frac{3}{16}$} -- (0,3/4);
				\draw [thick] (-.1,7/4) node[left]{$\frac{7}{16}$} -- (0,7/4);
				\draw [thick] (-.1,9/4) node[left]{$\frac{9}{16}$} -- (0,9/4);
				\draw [thick] (-.1,13/4) node[left]{$\frac{13}{16}$} -- (0,13/4);
				\draw [thick] (3/4,-.1) node[below]{$\frac{3}{16}$} -- (3/4,0);
				\draw [thick] (7/4,-.1) node[below]{$\frac{7}{16}$} -- (7/4,0);
				\draw [thick] (9/4,-.1) node[below]{$\frac{9}{16}$} -- (9/4,0);
				\draw [thick] (13/4,-.1) node[below]{$\frac{13}{16}$} -- (13/4,0);
				\draw [thick] (-.1,4) node[left]{\rm 1} -- (0,4);
				\draw (3/4,9/4) rectangle (7/4,13/4);
				\draw (9/4,3/4) rectangle (13/4,7/4);
				\draw[dashed, black] (3/4,9/4) -- (3/4,0);
				\draw[dashed, black] (7/4,9/4) -- (7/4,0);
				\draw[dashed, black] (9/4,3/4) -- (9/4,0);
				\draw[dashed, black] (13/4,3/4) -- (13/4,0);
				\draw[dashed, black] (3/4,9/4) -- (0,9/4);
				\draw[dashed, black] (3/4,13/4) -- (0,13/4);
				\draw[dashed, black] (9/4,3/4) -- (0,3/4);
				\draw[dashed, black] (9/4,7/4) -- (0,7/4);
				\node[below] (1) at (5/4,3) {\rm 3};
				\node[below] (2) at (11/4,3) {\rm 1};
				\node[below] (3) at (11/4,3/2) {\rm 3};
			\end{tikzpicture}
		\end{center}
		\caption{\rm The atomic mass number $z$ in the heterogeneous case for Example \ref{blast_hete_pb}.}
		\label{2dhete}
	\end{figure*}
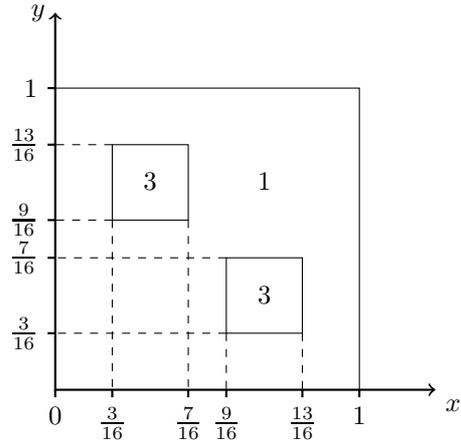

\begin{figure*}[htbp]
		\begin{minipage}[t]{0.32\linewidth}
			\centerline{\includegraphics[scale=0.38]{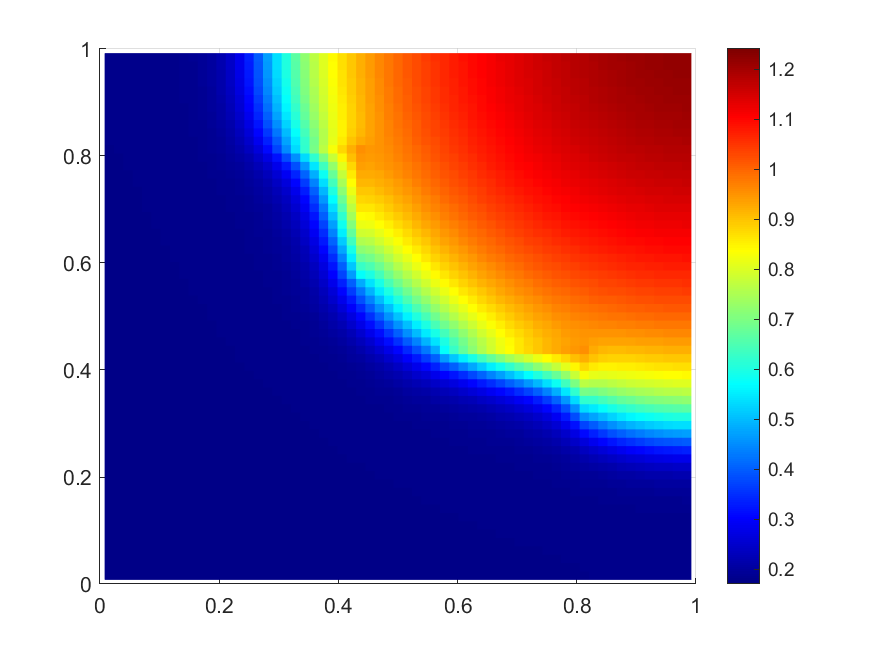}}
		\end{minipage}
		\begin{minipage}[t]{0.32\linewidth}
			\centerline{\includegraphics[scale=0.38]{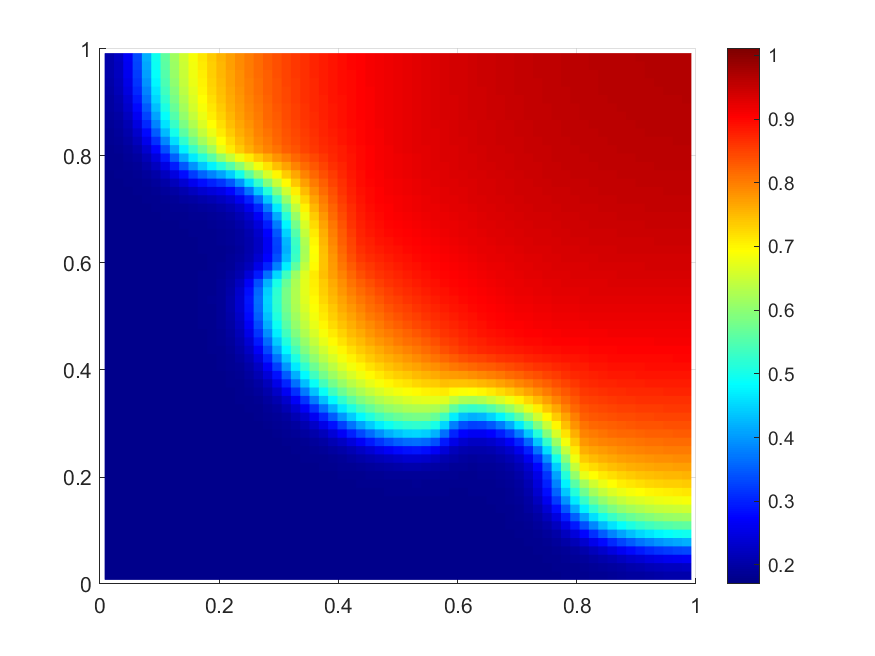}}
		\end{minipage}
		\begin{minipage}[t]{0.32\linewidth}
			\centerline{\includegraphics[scale=0.38]{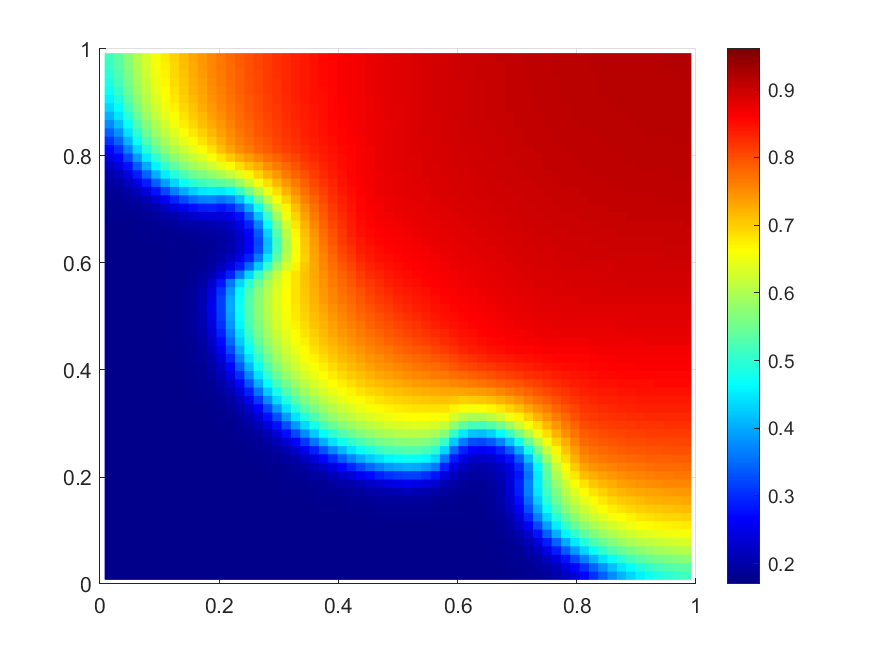}}
		\end{minipage}
		\vfill
		\begin{minipage}[t]{0.32\linewidth}
			\centerline{\includegraphics[scale=0.38]{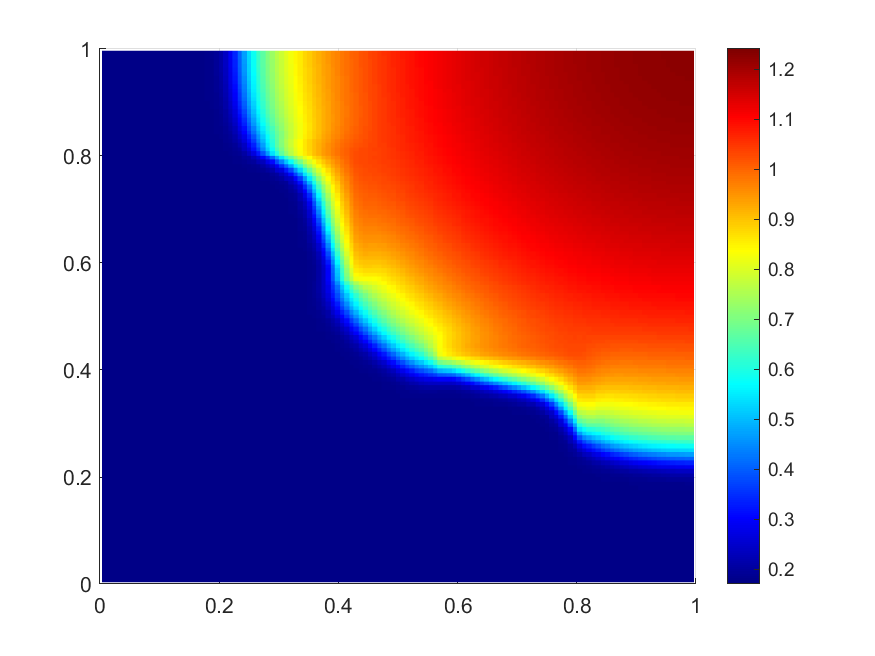}}
		\end{minipage}
		\begin{minipage}[t]{0.32\linewidth}
			\centerline{\includegraphics[scale=0.38]{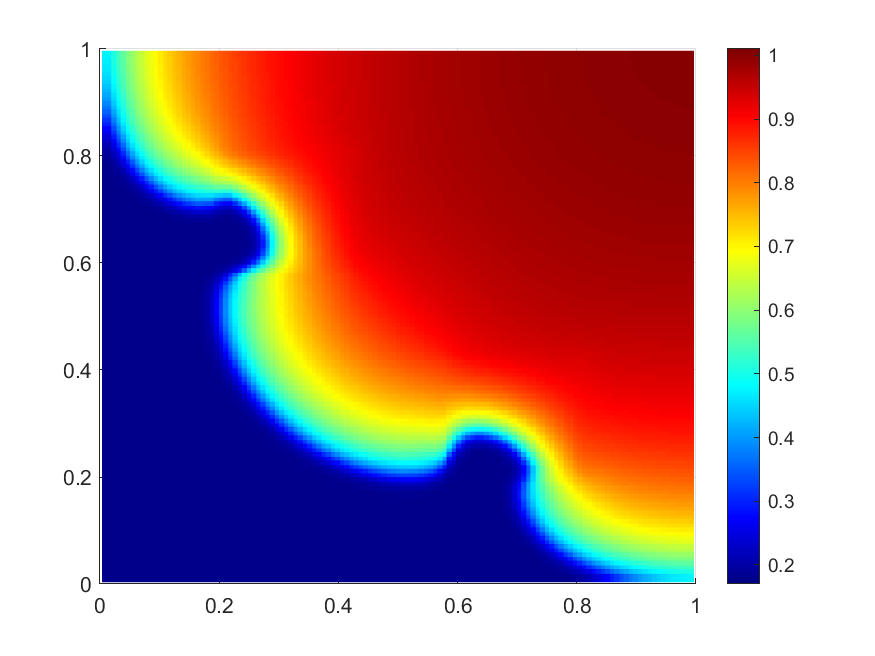}}
		\end{minipage}
		\begin{minipage}[t]{0.32\linewidth}
			\centerline{\includegraphics[scale=0.38]{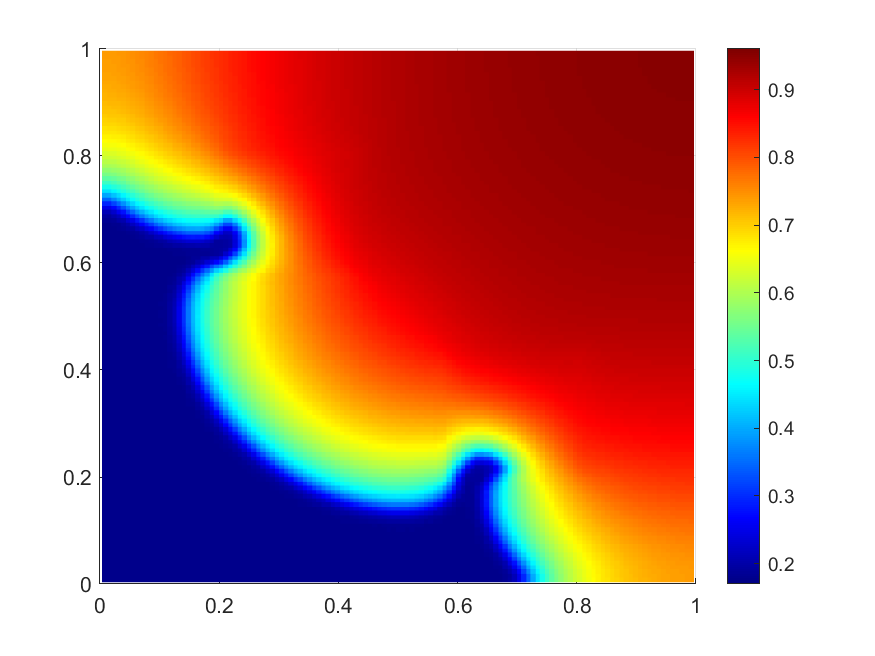}}
		\end{minipage}
		\vfill
		\begin{minipage}[t]{0.32\linewidth}
			\centerline{\includegraphics[scale=0.38]{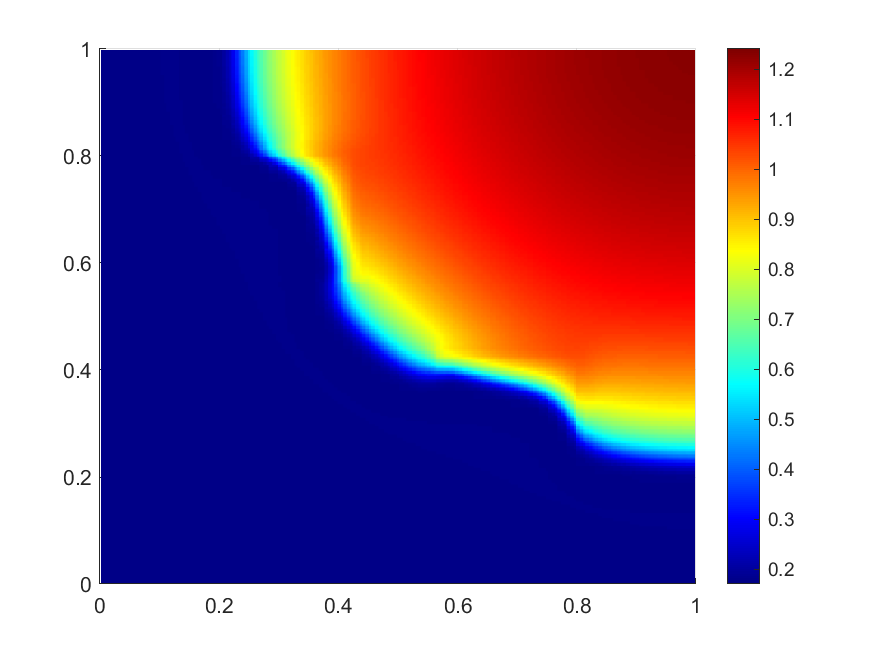}}
		\end{minipage}
		\begin{minipage}[t]{0.32\linewidth}
			\centerline{\includegraphics[scale=0.38]{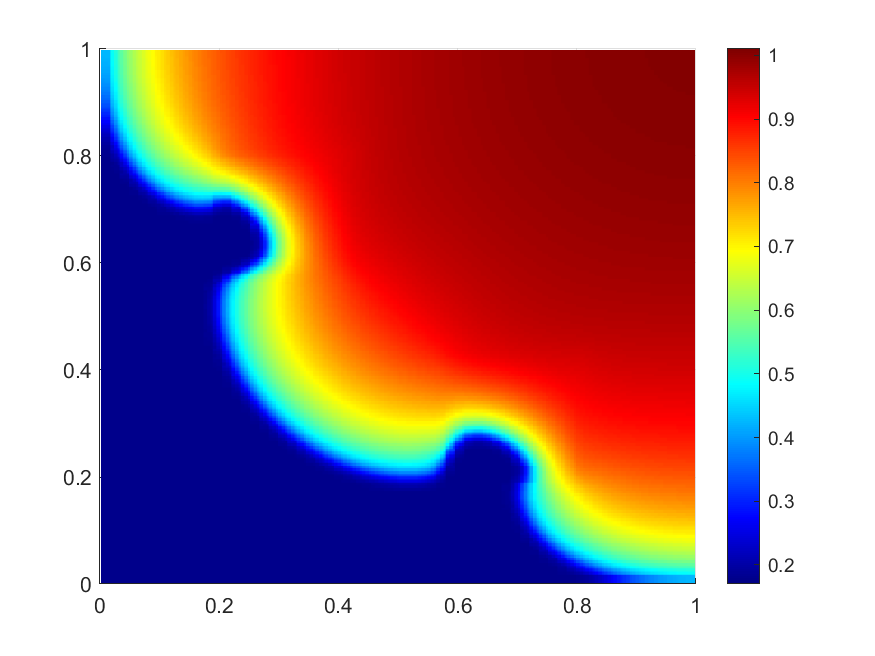}}
		\end{minipage}
		\begin{minipage}[t]{0.32\linewidth}
			\centerline{\includegraphics[scale=0.38]{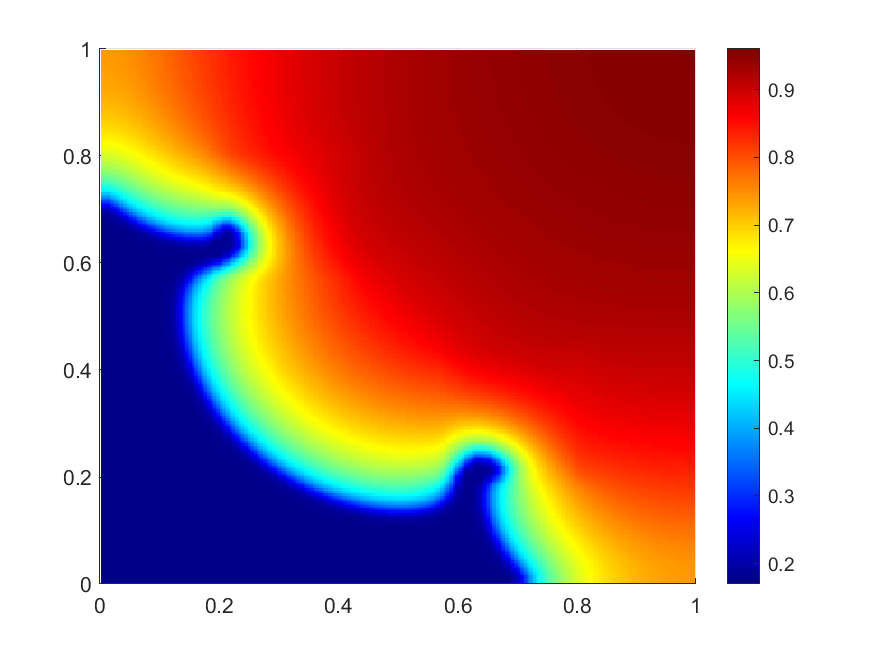}}
		\end{minipage}
		\caption{The numerical results for Example \ref{blast_hete_pb}. From top to bottom: the radiation temperature $T_r$ for 1st, 2nd, and 3rd order schemes. From left to right: time $t=1.0,2.0,2.5.$ $\Delta t =\frac{1}{10} h$.
		}\label{blast_hete_1}
	\end{figure*}

\begin{figure*}[htbp]
		\begin{minipage}[t]{0.32\linewidth}
			\centerline{\includegraphics[scale=0.38]{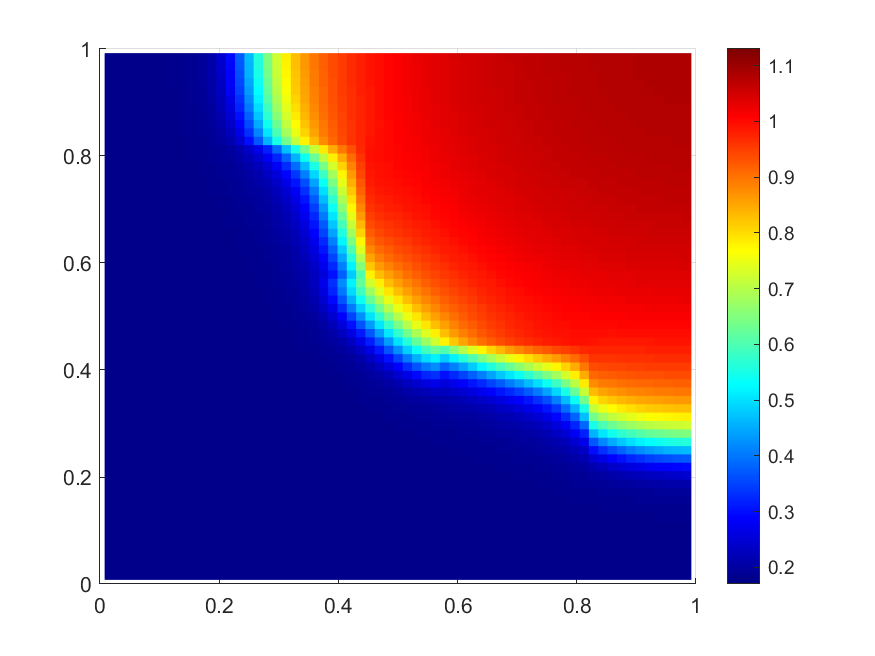}}
		\end{minipage}
		\begin{minipage}[t]{0.32\linewidth}
			\centerline{\includegraphics[scale=0.38]{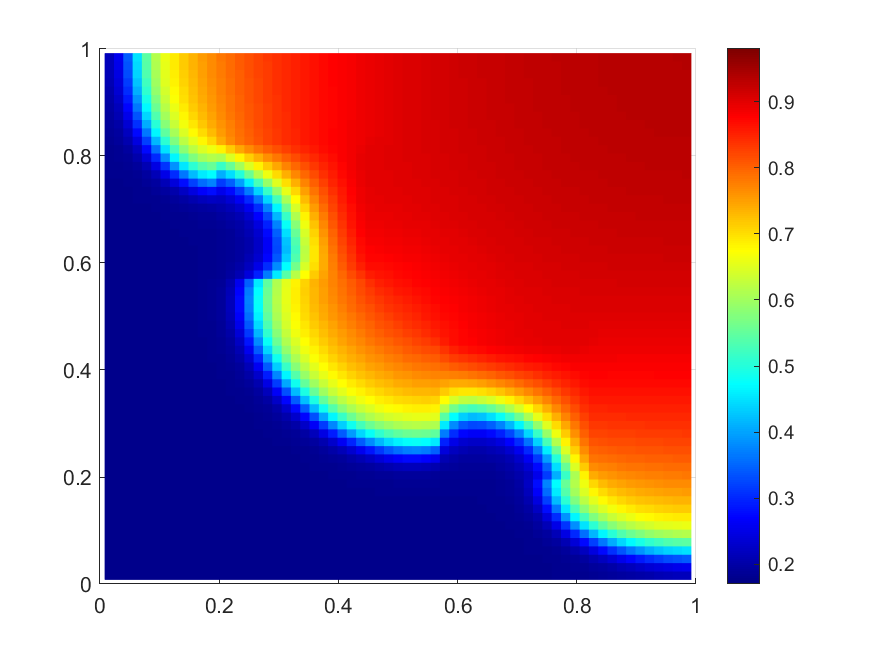}}
		\end{minipage}
		\begin{minipage}[t]{0.32\linewidth}
			\centerline{\includegraphics[scale=0.38]{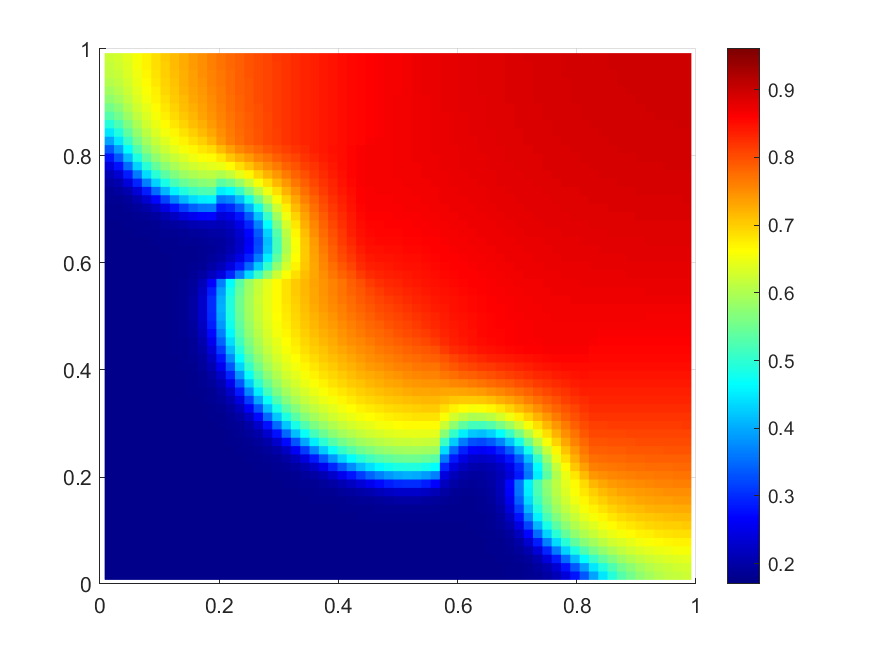}}
		\end{minipage}
		\vfill
		\begin{minipage}[t]{0.32\linewidth}
			\centerline{\includegraphics[scale=0.38]{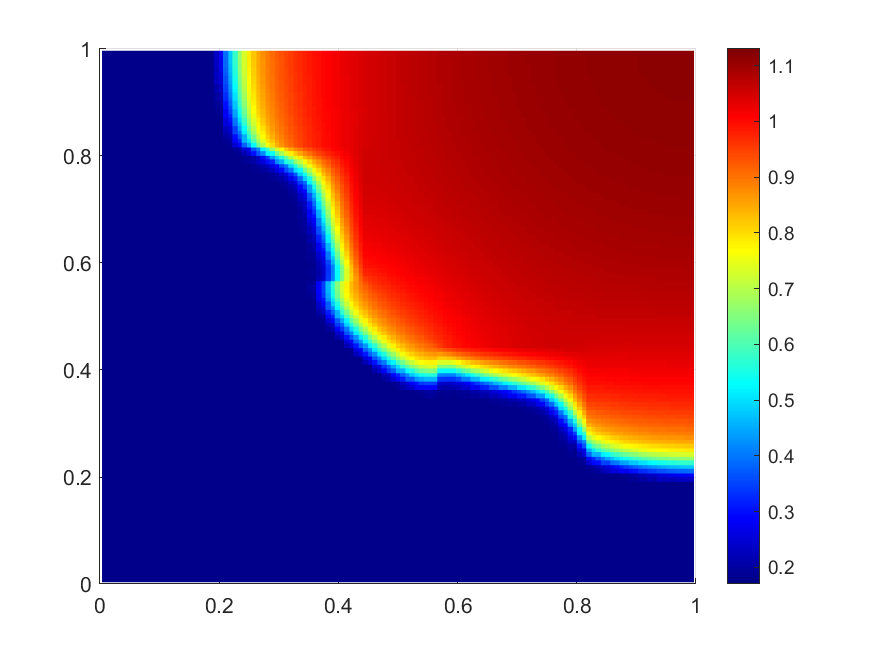}}
		\end{minipage}
		\begin{minipage}[t]{0.32\linewidth}
			\centerline{\includegraphics[scale=0.38]{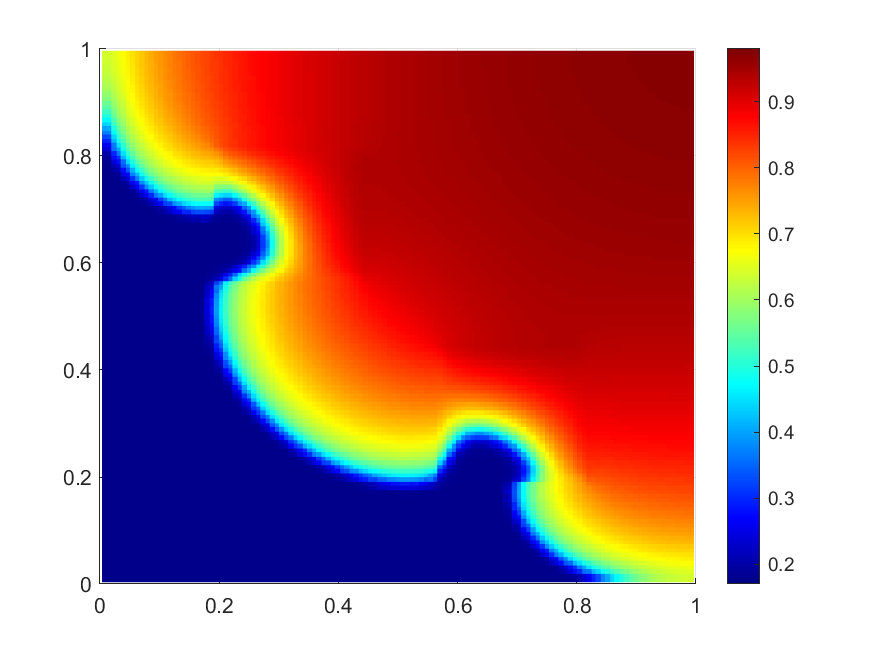}}
		\end{minipage}
		\begin{minipage}[t]{0.32\linewidth}
			\centerline{\includegraphics[scale=0.38]{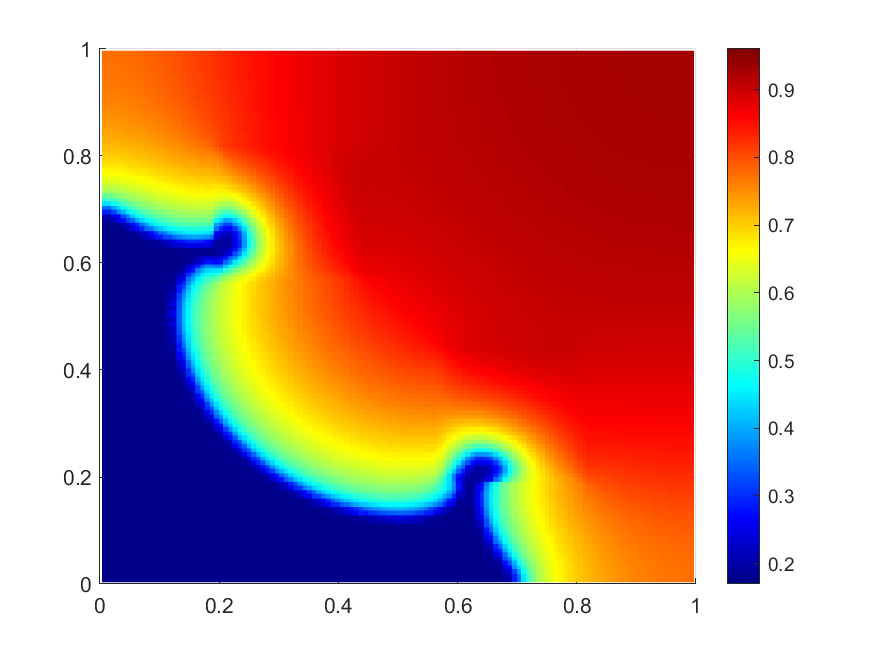}}
		\end{minipage}
		\vfill
		\begin{minipage}[t]{0.32\linewidth}
			\centerline{\includegraphics[scale=0.38]{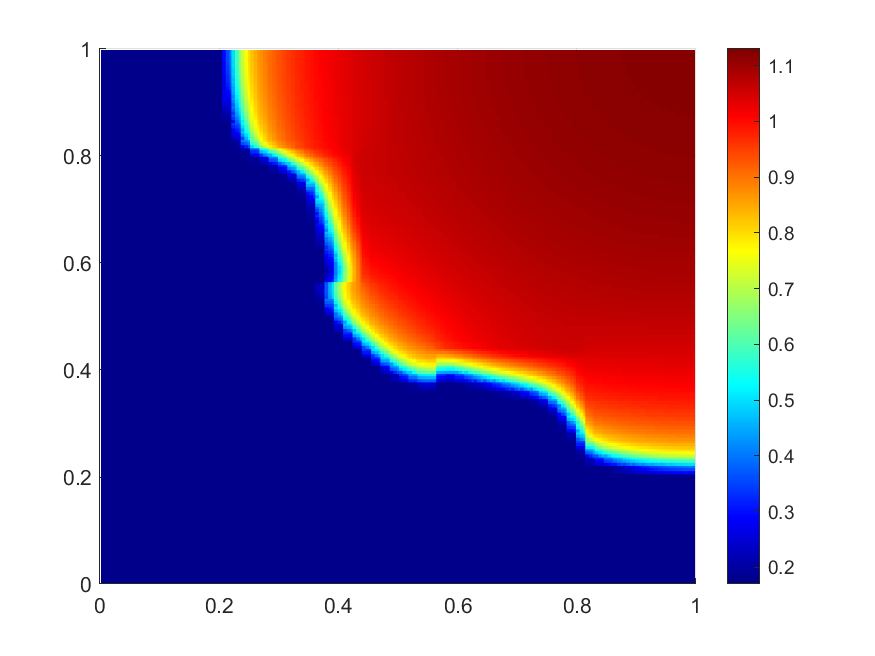}}
		\end{minipage}
		\begin{minipage}[t]{0.32\linewidth}
			\centerline{\includegraphics[scale=0.38]{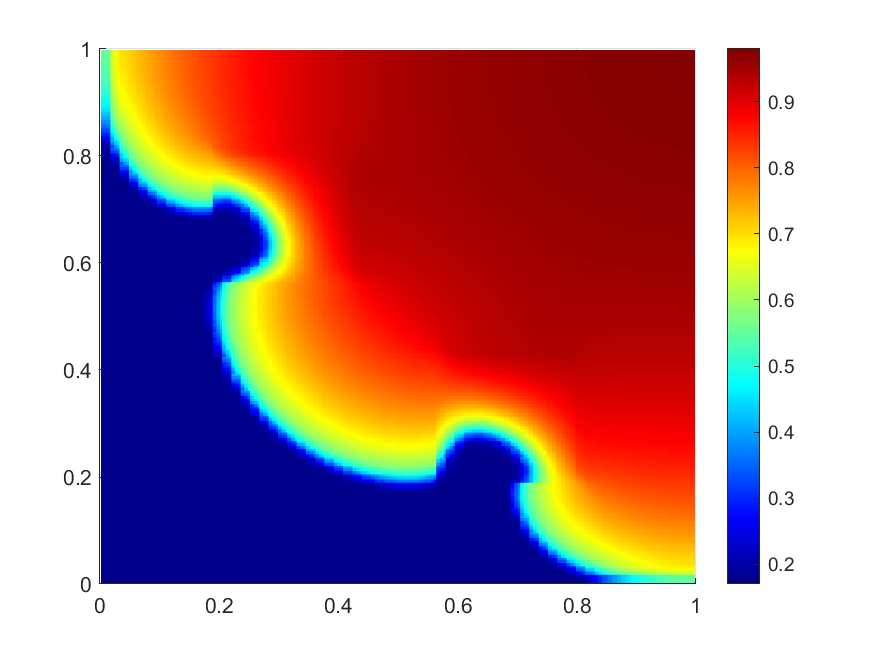}}
		\end{minipage}
		\begin{minipage}[t]{0.32\linewidth}
			\centerline{\includegraphics[scale=0.38]{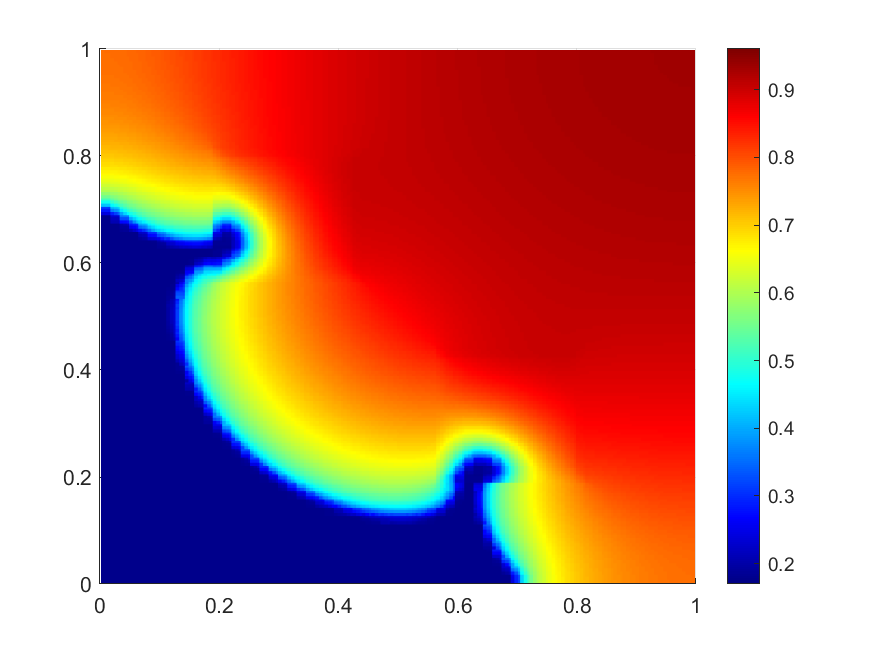}}
		\end{minipage}
		\caption{The numerical results for Example \ref{blast_hete_pb}. From top to bottom: the material temperature $T$ for 1st, 2nd, and 3rd order schemes. From left to right: time $t=1.0,2.0,2.5.$ $\Delta t =\frac{1}{10} h$.
		}\label{blast_hete_2}
	\end{figure*}

\begin{figure*}[!htbp]
		\begin{minipage}[t]{0.49\linewidth}
			\centerline{\includegraphics[scale=0.55]{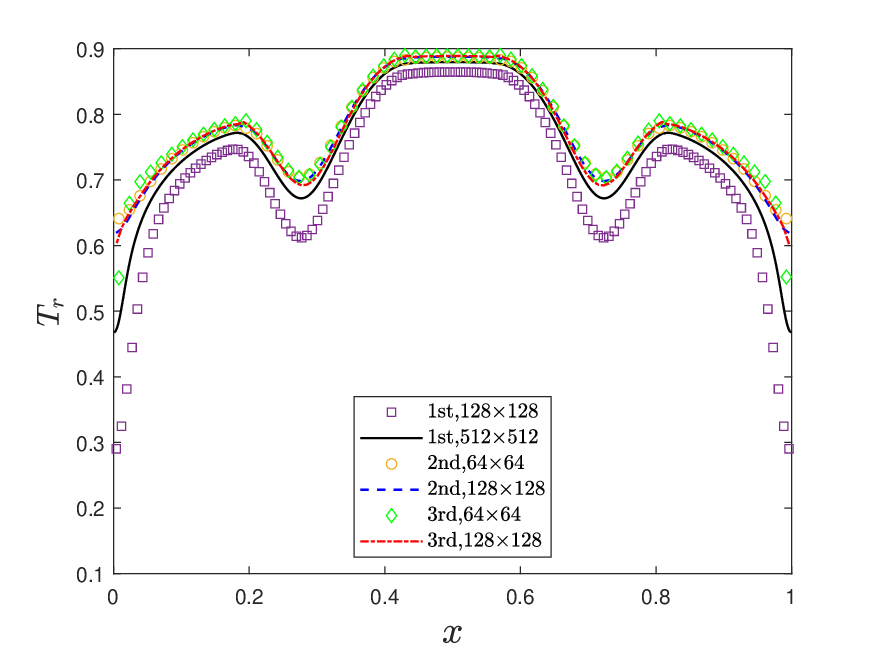}}
		\end{minipage}
		\begin{minipage}[t]{0.49\linewidth}
			\centerline{\includegraphics[scale=0.55]{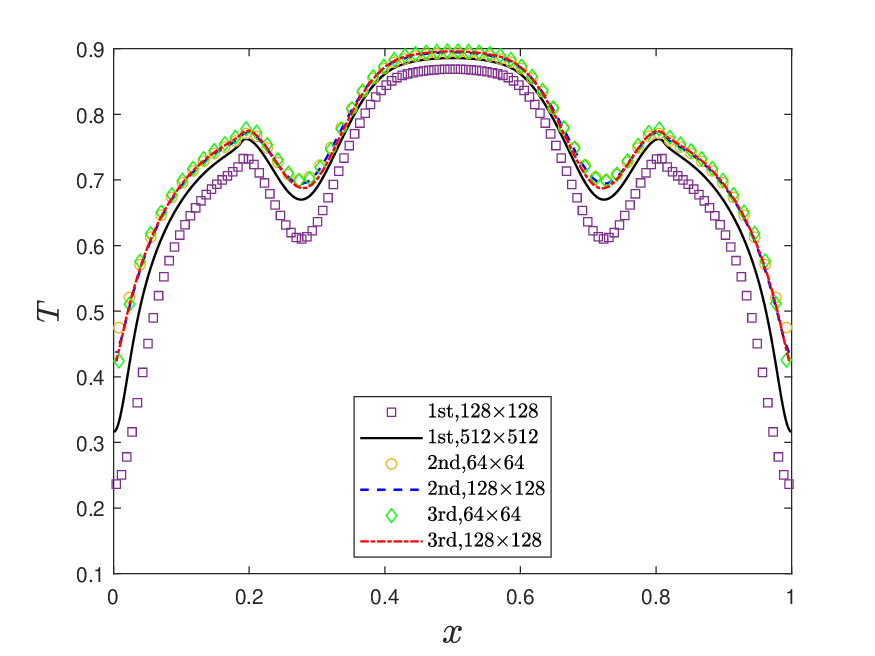}}
		\end{minipage}
		\caption{The numerical results for Example \ref{blast_hete_pb} at time $t=2$. Left: the radiation temperature $T_r$; Right: the material temperature $T$. $\Delta t =\frac{1}{10} h$.
		}\label{blast_hete_comparison}
	\end{figure*}

	\begin{example}\label{blast_hete_pb}
	{\rm Finally we consider the problems \eqref{RDE}-\eqref{RDE_2D_IBC} in a heterogeneous medium. The atomic mass number $z(x,y)$ is $1$ everywhere except in two inner square regions $\frac{3}{16}<x<\frac{7}{16}$, $\frac{9}{16}<y<\frac{13}{16}$ and $\frac{9}{16}<x<\frac{13}{16}$, $\frac{3}{16}<y<\frac{7}{16}$, where the value of $z(x,y)$ is $3$, as shown in Fig. \ref{2dhete}. In Fig.~\ref{blast_hete_1} and Fig.~\ref{blast_hete_2}, we present the numerical solutions $T_r$ and $T$ at time $t=1,2,2.5$ on a mesh of $64\times 64$, respectively. 
		As compared to $\Delta t = \frac{1}{32}h$ in \cite{huang2016monotone,zhao2013finite}, we use a lager time step $\Delta t=\frac{1}{10}h$.
		From top to bottom, numerical results obtained by the 1st, 2nd, and 3rd order methods are presented, respectively. From the results, we find that the results of the 1st order method are very smeared due to numerical viscosities. 2nd and 3rd order methods perform better. In Fig.~\ref{blast_hete_comparison}, we show the cutting plots for the numerical solutions along $y=x+1$ at $t=2$ on a mesh of $N\times N$ for $N=64, 128, 512$. We can observe that by a mesh refinement, the results of 2nd and 3rd methods match each other. In addition, we observe that when the front of the radiation temperature propagates to the interface of two different materials, the process is hindered by a region with dense atoms. Higher energy exchange $\sigma = z(\bx)^3 / T^3$ occurs in these regions, resulting in a well balance between the material and radiation temperatures. The temperature propagates faster in the lower $z(\bx)$ region. Our results agree with those in \cite{guangwei2009progress,yang2019moving}.
		
	}
    \end{example}
	
	\begin{rmk}\label{numberofiterations}
		{\rm We are not able to prove the convergence of the Picard iteration for \eqref{algebraic system} theoretically. For the 2nd and 3rd order methods in space, the matrix $A(U^{(i)})$ may not be diagonally dominant, so that $A(U^{(i)})$ cannot be guaranteed to be an $\mathbf{M}$ matrix. However, in our numerical examples, all results converge quickly under a tolerance $\delta=10^{-8}$. In Table \ref{table_of_iterations}, we show an average number of iterations per stage for the nonlinear algebraic system \eqref{algebraic system} in the heterogeneous medium. We can observe that in such challenging cases, our methods with large time steps converge around four iterations, which show the robustness of our proposed methods.}
	\end{rmk}

	\begin{table}[!htbp]
		\centering
		\caption{Average iteration numbers per stage for some examples in the heterogeneous medium.}\label{table_of_iterations}
		\begin{tabular}{ccccc}
			\toprule
			&     & \multirow{2}{*}{mesh}    & \multirow{2}{*}{\begin{tabular}[c]{@{}c@{}}time \\ step\end{tabular}}             & \multirow{2}{*}{\begin{tabular}[c]{@{}c@{}} an average number \\ of iteration per stage\end{tabular}} \\ 
			&     &                          & \\ \hline
			\multirow{2}{*}{\begin{tabular}[c]{@{}c@{}} 1D heterogeneous Marshak\\ wave(Example \ref{Marshak_hete_problem})\end{tabular}} & 2nd & \multirow{2}{*}{180} & \multirow{2}{*}{$\frac{1}{15}h$} & 4.02                         \\ 
			& 3rd &  &  & 5.05                \\\hline
			\multirow{3}{*}{\begin{tabular}[c]{@{}c@{}} 2D heterogeneous Blast\\ wave(Example \ref{blast_hete_pb})\end{tabular}} & 1st & \multirow{3}{*}{$64\times 64$} & \multirow{3}{*}{$\frac{1}{10}h$} & 4.12                        \\ 
			& 2nd &  &  & 3.99                 \\
			& 3rd &  &  & 4.00                 \\
			\bottomrule
		\end{tabular}
		
	\end{table}
	
	
	
	\section{Conclusion}\label{Conclusion}
	
	In this paper, we have developed a class of high order conservative LDG-IMEX methods for non-equilibrium radiation diffusion problems. The proposed scheme is based on a predictor-corrector approach. During the predictor step, we solve a reformulated system to improve convergence and approach thermodynamic equilibrium. Subsequently, in the corrector step, we solve the original system using nonlinear coefficients and initial estimates obtained from the predictor step. This dual-step process ensures the conservation of total energy and robust convergence. We incorporate linear diffusion terms to circumvent implicit discretization for nonlinear diffusion terms. This approach permits the use of larger time step sizes, in contrast to the restrictive parabolic time step conditions $\Delta t = \mathcal{O}(h^2)$ associated with purely explicit discretizations. Local discontinuous Galerkin finite element approximations in space are employed. Numerical examples in both 1D and 2D illustrate the benefits of high-order conservative methods in accurately capturing steep solution fronts within both homogeneous and heterogeneous media. Our proposed methods exhibit robustness, as confirmed by numerical experiments. However, theoretical analysis to ensure the convergence of such a nonlinear iteration is a challenging task, which we plan to explore in our future work. The extension to a 3D radiation diffusion model \cite{lai2017monotone} or three-temperature (3T) model \cite{su2020vertex,yu2019finite}, and other reaction-diffusion equations \cite{smoller2012shock} will also be investigated.
	
	\newpage
	
	\appendix

	\section{IMEX Butcher tableau}\label{IMEX Butcher tableau}
	In this paper, the double Butcher tableaux we used from \cite{ascher1997implicit} are listed below:
	\begin{tableau}\label{Tableau1}{\rm 2 stage 1st order:}
		\begin{table}[!htbp]
			\centering
			\begin{tabular}{ c|cccc|cc}
				0&0&0&&0&0&0\\
				1&0&1&&1&1&0\\
				\cline{1-3}\cline{5-7}
				&0&1&&&1&0\\
			\end{tabular}.
		\end{table}
	\end{tableau}

	\begin{tableau}\label{Tableau2}{\rm 3 stage 2nd order:}
		\begin{table}[!htbp]
			\centering
			\begin{tabular}{c|ccccc|ccc}
				0&0&0&0&&0&0&0&0\\
				$\gamma$&0&$\gamma$&0&&$\gamma$&$\gamma$&0&0\\
				1&$\gamma$&0&$\gamma$&&1&0&1&0\\
				\cline{1-4}\cline{6-9}
				&$\gamma$&0&$\gamma$&&&0&1&0\\
			\end{tabular},
		\end{table}
	\end{tableau} 
	where $\gamma=\frac12$.
	\begin{tableau}\label{Tableau3}{\rm 5 stage 3rd order:}
		\begin{table}[!htbp]
			\centering
			\begin{tabular}{ c|ccccccc|ccccc}
				0&0&0&0&0& 0&&0&0&0&0&0& 0\\
				1/2&0&1/2&0&0& 0&&1/2&1/2&0&0&0&0\\
				2/3&0&1/6&1/2&0&0&&2/3&11/18&1/18&0&0&0\\
				1/2&0&-1/2&1/2&1/2&0&&1/2&5/6&-5/6&1/2&0& 0\\
				1&0&3/2&-3/2&1/2&1/2&&1&1/4&7/4&3/4&-7/4& 0\\
				\cline{1-6}\cline{8-13}
				&0&3/2&-3/2&1/2&1/2&&&1/4&7/4&3/4&-7/4& 0\\
			\end{tabular}.
		\end{table}
	\end{tableau}

	\bibliographystyle{plain}
	\bibliography{reference}
	
\end{document}